\documentclass{article}
\usepackage{amsmath,amssymb,amsthm,mathtools}
\usepackage{graphicx}
\usepackage{hyperref}
\usepackage{xcolor}
\usepackage{pgfplots}
\usepackage{bbm}
\usepackage{thm-restate}
\usepackage{booktabs}
\usepackage{subcaption}
\usepackage{enumitem}
\usepackage{todonotes}
\usepackage{verbatim}
\usepackage{wrapfig}
\usepackage{siunitx} 
\usepackage[bottom]{footmisc}
\usepackage{iclr2026_conference,times}
\newcommand{\ul}[1]{\underline{#1}}

  \providecommand{\R}{\mathbb{R}} %

  \providecommand{\cB}{\mathcal{B}}

  \providecommand{\cE}{\mathcal{E}}
  \providecommand{\cF}{\mathcal{F}}
  \providecommand{\cG}{\mathcal{G}}
  
  \providecommand{\cI}{\mathcal{I}}

  \providecommand{\cP}{\mathcal{P}}
  \providecommand{\cQ}{\mathcal{Q}}

  \providecommand{\cT}{\mathcal{T}}

  \providecommand{\cX}{\mathcal{X}}
  \providecommand{\cY}{\mathcal{Y}}
  \providecommand{\cW}{\mathcal{W}}
  \providecommand{\cZ}{\mathcal{Z}}

  \usepackage{bm}

\providecommand{\mycomment}[3]{\todo[caption={},size=footnotesize,color=#1!20]{\textbf{#2: }#3}}%
\providecommand{\inlinecomment}[3]{%
  {\color{#1}#2: #3}}%
\newcommand\commenter[2]%
{%
  \expandafter\newcommand\csname i#1\endcsname[1]{\inlinecomment{#2}{#1}{##1}}
  \expandafter\newcommand\csname #1\endcsname[1]{\mycomment{#2}{#1}{##1}}
}

\newtheorem{proposition}{Proposition}
\newtheorem{lemma}{Lemma}
\newtheorem{corollary}[lemma]{Corollary}

\newtheorem{definition}{Definition}
\newtheorem{remark}[lemma]{Remark}

\newtheorem{theorem}[lemma]{Theorem}

\newcommand{\e}{\mathbb{E}}

\newcommand{\vp}{\varphi}
\newcommand{\ve}{\varepsilon}
\newcommand{\II}{\mathbb{I}}

\usetikzlibrary{intersections}
\pgfplotsset{compat=1.18}

\theoremstyle{definition}

\title{Infinitely divisible  privacy and beyond I\\ resolution of the $s^2=2k$ conjecture}

\author{
\textbf{Aaradhya Pandey}\textsuperscript{1},
\textbf{Arian Maleki}\textsuperscript{2}, \textbf{Sanjeev Kulkarni}\textsuperscript{1}
 \\
\textsuperscript{1}Princeton University\quad \textsuperscript{2} Columbia University \\
}

\iclrfinalcopy

\usepackage{float}
\begin{document}

\maketitle

\begin{abstract} 
Differential privacy is increasingly formalized through the lens of hypothesis testing via the \textit{robust} and \textit{interpretable}$f$-DP framework, where privacy guarantees are encoded by a baseline \textit{Blackwell} trade-off function $f_{\infty}= T(P_{\infty},Q_{\infty})$ involving a pair of distributions $(P_{\infty}, Q_{\infty})$ \cite{dong2022gdp}. The problem of  `\textit{choosing the right privacy metric in practice}' \cite{Cummings_HDSR_2024_AdvancingDP} gives rise to a central question in this framework: \textit{what is a statistically appropriate baseline $f_{\infty}$ given some prior modeling assumptions}? The special case of Gaussian differential privacy (GDP) \cite{dong2022gdp} showed that, under compositions of nearly perfect DP mechanisms, these baseline trade-off functions exhibit a central limit behavior: $f_n=T(P_n^{\otimes n}, Q_n^{\otimes n}) \to f_{\infty}= T(P_{\infty},Q_{\infty})$, where the limiting baseline trade-off function $f_{\infty}$ involves two shifted Gaussians $P_{\infty}, Q_{\infty}$. Inspired by Le Cam's theory of limits of statistical experiments \cite{LeCam1986}, \textit{we answer this question in full generality in its natural infinitely divisible setting}. 

We show that the sequence of composition experiments $(P_n^{\otimes n}, Q_n^{\otimes n})$ `\textit{converge}' in general to a binary limit experiment $(P_{\infty}, Q_{\infty})$ whose log-likelihood ratio  $\log \frac{dQ_{\infty}}{dP_{\infty}}$ is an infinitely divisible distribution under the limiting null $P_{\infty}$. So, any such limiting trade-off function $f_{\infty}$ involves an infinitely divisible law $P_{\infty}$ characterized by its \textit{Lévy–Khintchine} triplet,  and its \textit{Esscher} tilt $dQ_{\infty}(x) = e^{x} dP_{\infty}(x)$. As a consequence, we determine all the limiting baseline trade-off functions $f_{\infty}$ under the composition of nearly perfect differentially private procedures. Our framework recovers GDP as the purely Gaussian case,  yields explicit many non-Gaussian limits including Poisson, and many others. Consequently, we \textit{positively} resolve the $s^2=2k$  conjecture proposed in the GDP paper \cite{dong2022gdp}, which curiously observed that in the Gaussian examples $P_{\infty} \overset{d}{=} -k+ |s|Z$ with $Z \sim N(0,1)$,  we always have $s^2=2k$. Then, we describe situations where one can goes beyond the infinite divisibility framework, such as $\log  \frac{dQ_{\infty}}{dP_{\infty}}$ being a Gaussian mixture under $P_{\infty}$. Finally, we describe a mechanism for count statistics that optimally achieves (asymmetric) Poisson differentially privacy, going beyond the natural noise-adding mechanisms used in practice. 
\end{abstract}

\section{Introduction,  previous work, and our contributions}

Modern data analysis routinely draws on information about people: their location traces, search histories, media consumption, social graphs, alongside long-standing, highly sensitive records from hospitals, schools, and national censuses. This necessarily demands a framework of data analysis that preserves the privacy of individuals. Early \textit{privacy} practices leaned on ad-hoc anonymization, but a series of high-profile re-identifications made it clear that simple de-identification is brittle at scale \cite{narayanan2008netflix}. This motivated a rigorous  foundation for privacy in statistical analysis. Differential privacy (DP) provides exactly that: a principled, mathematically verifiable guarantee against disclosure risks that persists under arbitrary post-processing and composition. Since its introduction, differential privacy has grown into a practical standard, with widespread industrial adoption \cite{apple2017dp} and deployment in official statistics \cite{abowd2018census}.

\subsection{From $(\varepsilon,\delta)$ and divergence based DP to  Hypothesis testing based DP}

In its classical form,   differential privacy formalizes protection via two real-valued parameters $\varepsilon\ge 0$,   $0\le \delta \le 1$, which jointly bound the advantage of any adversary attempting to distinguish between neighboring datasets from a mechanism’s output  \cite{dwork2006eurocrypt, dwork2006tcc}.
The \((\varepsilon,\delta)\) formulation of differential privacy became the standard due to its robustness to post-processing, group privacy, and composability, but sharp accounting under composition proved to be technically delicate \cite{dwork_roth_2014}. A line of \emph{divergence-based} relaxations was subsequently developed to tighten composition analyses and enable modular accounting: concentrated DP (CDP and zCDP)~\cite{dworkrothblum2016cdp} \cite{bun2016concentrated}, and Rényi DP (RDP)~\cite{mironov2017renyi}, with follow-ups on amplification by subsampling and analytic accountants~\cite{wang2019subsampled,abadi2016deep}. 

In parallel, a practical, robust, and  interpretable \emph{hypothesis-testing} viewpoint of differential privacy was formalized, starting with \cite{wasserman2010statistical} and  \cite{kairouz2015composition}. It describes privacy as the hypothesis testing problem of distinguishing between two neighboring datasets  $S,S'$ via the output distributions $M(S), M(S')$ of a private mechanism $M$. This was formalized in the framework of $f$-differential privacy ($f$-DP)  \cite{dong2022gdp}, which describes privacy of $M$ by requiring that the optimal type~I vs type~II error curve for the distinguishing test for $M(S), M(S')$, for any two neighboring datasets,  $S$ and $S'$,  stays above a baseline \textit{trade-off function} (TOF) $f$ .

Immediately, the $f$-DP framework is statistically interpretable via hypothesis testing. It recovers the classical \((\varepsilon,\delta)\)-DP formulation for a certain choice of baseline trade-off function $f_{\ve, \delta}$ \cite{wasserman2010statistical}\footnote{The expression of $f_{\varepsilon,\delta}$ is given by $f_{\varepsilon, \delta}(\alpha) = \max\left(0, 1-\delta-e^{\varepsilon}\alpha, e^{-\varepsilon}(1-\delta- \alpha)\right)$ for $\alpha \in [0,1]$.}, and it is a robust framework because it uniformly preserves privacy at all scales of type $I$ error $0\leq \alpha \leq 1$.  Moreover, it is achievable by appropriate noise injection mechanisms \cite{awan2023canonical}\footnote{Although very general, the noise adding mechanisms described here are not quite naturally implementable.}, and it is closed and tight under composition, thereby allowing for lossless reasoning about composition \cite{dong2022gdp}. Finally, \cite{dong2022gdp} introduced Gaussian differential privacy (GDP) as a one-parameter family of baseline trade-off curves $G_{\mu}$ that require distinguishing between a pair of shifted normal distributions $N\left(-\frac{\mu^2}{2}, \mu^2\right)$ and   $N\left(\frac{\mu^2}{2}, \mu^2\right)$ with the same variance $\mu^2$ for a real-valued parameter $\mu \in \mathbb{R}$\footnote{From the basic symmetries of Gaussian trade-off functions $G_{\mu}$ \cite{Pandey2025GaussianCertifiedUnlearning},  our description of $G_{\mu}$ is equivalent to the usual description of itself involving the pair of distributions $N(0,1)$ and $N(\mu,1)$. However, our formulation will make the resolution of the $s^2= 2k$ conjecture of \cite{dong2022gdp} appear naturally.} and showed that it optimally captures the widely used Gaussian mechanism.  \cite{dong2022gdp} further showed that under compositions of nearly perfect DP mechanisms $M_n\otimes \cdots \otimes  M_1$, these baseline trade-off functions $f_n \otimes \cdots \otimes f_1$ converge to a limiting baseline trade-off function $G_{\mu}$ involving two shifted Gaussians as above.

\subsection{Our contribution: In search of a baseline trade-off function} 

\cite{dong2022gdp} established a central limit behavior within differential privacy, indicating that the baseline trade off function $f$ should be chosen as a trade off function $T(P_{\infty}, Q_{\infty})$ (see definition \ref{def:TOF}) for a shifted Gaussian pair $P_{\infty}=N\left(-\frac{\mu^2}{2}, \mu^2\right)$ and   $Q_{\infty}=N\left(\frac{\mu^2}{2}, \mu^2\right)$, for some $\mu\in \mathbb{R}$ under a larger number of (nearly perfect) differentially private outputs based on queries on the same dataset. 
\begin{itemize}
    \item \textbf{An infinitely divisible framework of privacy.} In this paper, we generalize this result in its natural infinitely divisible setting (see Theorem \ref{thm:IDP}) and establish that under a large number of (nearly perfect) differentially private operations, if the deteriorated trade-off curves converge to a limiting trade-off curve, then the limiting trade-off curve has an infinitely divisible form (see Equation \eqref{eq:IDP}), encompassing both Gaussian and Poisson components\footnote{Inspired by   \cite{LeCamYang2000} `The gaussian ones are
used everywhere because of mathematical tractability. Poisson experiments are less tractable and less studied, but they are probably more important. They will loom large in the new century.', and given that we are in the next century, we show how the Poisson baselines are extremely relevant for differentially private analysis on graph based relational datasets.}.
\item \textbf{Resolution of the $s^2=2k$ conjecture.} As an application of  Theorem \ref{thm:IDP}, we resolve the $s^2=2k$ conjecture (see Theorem \ref{thm:GDP}) proposed in  \cite{dong2022gdp}, which  observed that in the Gaussian case $P_{\infty} \overset{d}{=} -k+ |s|Z$ with $Z \sim N(0,1)$,  we always have $s^2=2k$.

\item \textbf{Random number of operations and beyond infinite divisibility.} In practice, we do not know the number of differentially private operations one is going the implement on the same dataset. In section \ref{ssec:practice} (see Equation \eqref{eq:Levylimit}) we derive that under a large but random number of (nearly perfect) differentially private operations, the limiting trade-off functions can go beyond the class of infinitely divisible trade-off functions. Inspired by \cite{LeCam1986}, we demonstrate this with a locally asymptotically mixed normal example \ref{eq:GMM}.

 \item \textbf{Achieving Poisson differential privacy in practice.} In Theorem  \ref{thm:PDP} , we describe a mechanism that optimally achieves (asymmetric) Poisson differential privacy real-valued statistics, including count statistics with baseline privacy $f= T(P(\lambda_1), P(\lambda_2))$ going beyond the popular noise -adding mechanisms.

\item \textbf{Applications to graph-based relational datasets.} We demonstrate how Poisson differential privacy becomes extremely relevant in tightly capturing privacy on datasets that are relational or graph-based and only allow released statistics to be inherently discretized  \ref{sec:BIDP}.

 \item \textbf{A coarsened Neyman-Pearson lemma.} Inspired  by  a Bayesian perspective of  differential privacy  \cite{StrackYang2024PrivacyPreserving}, we prove (see Theorem \ref{thm:NP-coarsen}) exactly how the trade-off functions $T(P,Q)$ on $(\Omega, \cF)$ degrade when we shrink the underlying the $\sigma$-algebra of information from $\cF$ to $\cG$, and demonstrate with a Gaussian example\footnote{The results extend with appropriate modifications for any symmetric log concave shift family on $\mathbb{R}$.} (see Proposition \ref{prop:GTOFC}). 
\end{itemize}
 
\section{Hypothesis Testing-Based Privacy:  Definitions and Motivations}

We describe objects of the hypothesis-testing-based privacy framework and motivate our problem. 

\subsection{Basics of trade-off functions and Blackwell ordering}

 \begin{definition}[Trade-off function \cite{dong2022gdp}]    \label{def:TOF} 
Given two probability distributions $P, Q$ on a measurable space $(\cW, \mathcal{F}_{\cW})$, we define the  \emph{trade-off function} as  the map $= T(P,Q):[0,1]\to[0,1]$   
\begin{equation}\label{eq:TOF}
 T(P,Q)(\alpha):=
  \inf_{\varphi}\Bigl\{
        \beta_{\vp} := \e_{Q}[1-\vp] 
        \,\Bigm\vert\,
        \alpha_{\vp}:= \e_{P}[\vp]\le \alpha,\;
        \vp:(\cW, \cF_{\cW})\!\to[0,1]\text{ measurable}
      \Bigr\}.
\end{equation}
\textbf{Intuition, the Neyman-Pearson optimizer.} For  any type I error $\alpha \in [0,1]$, the trade-off function (TOF) returns the smallest value of type II error $\beta_{\vp}$ over all possible Borel measurable test functions $\vp\colon \cW\to[0,1]$. As a consequence of Neyman-Pearson lemma  \cite{polyanskiy_wu_2025}[Thm 14.11] $T(P,Q)(\alpha)$ is achieved by the following \text{likelihood ratio test} for every $\alpha \in  [0,1]$
        \begin{equation}\label{eq:LRT}
            \varphi^*(w) :=\mathbf{1}\left( \frac{dQ}{d\mu} >\tau^* \frac{dP}{d\mu} \right) + \lambda^* \mathbf{1}\left(\frac{dQ}{d\mu} =\tau^* \frac{dP}{d\mu}\right), \quad \frac{dQ}{dP}: = \frac{dQ}{d\mu}/\frac{dP}{d\mu}:(\Omega, \cF) \to [0, \infty]
        \end{equation}
  for $2\mu=P+Q$, are the ratio of Radon–Nikodym derivatives, and $\tau^*, \lambda^*$ are determined  as follows. Choose $\tau^* \in [0,\infty]$ as the unique number\footnote{Such a $\tau^*$ uniquely exists  as the CDF $F:\mathbb{R}\to [0,1],$ defined as $ F(t):= P\left(\frac{dQ}{dP} \leq t\right)$ is monotone increasing as a function of $t$ satisfying $(t_1\leq t_2\rightarrow F(t_1)\leq F(t_2))$, right continuous ($\lim_{t\downarrow t_0}F(t)=F(t_0)$) with $F(0^-)= P_{0}\left(\frac{dQ}{dP} < 0\right)=0$ (lower limit of $\alpha$) and $\lim_{t \to \infty}F(t)=1$ (upper limit of $\alpha$).} such that $P\left(\frac{dQ}{dP} \geq \tau^{*}\right) \geq \alpha \geq  P\left(\frac{dQ}{dP} > \tau^{*}\right)$, and 
  \[\text{define }\lambda^*:= \frac{\alpha - P\left(\frac{dQ}{dP} > \tau^{*}\right) }{P\left(\frac{dQ}{dP} \geq \tau^{*}\right)  - P\left(\frac{dQ}{dP} > \tau^{*}\right) } \mathbf{1}\left(\alpha - P\left(\frac{dQ}{dP} > \tau^{*}\right) >0\right) \text{ to have } \mathbb{E}_{P}[\varphi^*]=\alpha.\] 

\textbf{Characterization of TOFS and their inverse.} \cite{dong2022gdp}[Prop 1] proved that the collection of all trade-off functions $\cT:=\{T(P,Q): (P,Q) \text{ on some measurable } (\cW, \cF_{\cW})\}$ is the same as\footnote{As argued in \cite{Torgersen1991Comparison} and is also immediate that the set $\cT$ says the same, if we keep all the conditions intact, except one, where we replace the condition $f(x)\leq 1-x$ for all $x\in [0,1]$ by just $f(1)=0$.} 
\begin{equation} \label{def:TOFs}
\cT= \{f: [0,1] \to [0,1] \text{ convex, continuous, decreasing and } f(\alpha)\leq 1-\alpha \text{ for all } \alpha\in [0,1]\}
\end{equation}
Moreover,  any such $f\in \cT$ can be realized as $T(P,Q)$ on  $\cW=[0,1], \cF_{\cW}= \cB_{[0,1]}$ by taking $P$ to be the standard uniform probability measure $U$ on $([0,1],\cB_{[0,1]})$ and $Q=Q_f$ to be a Borel probability measure on $[0,1]$ with cumulative distribution function $Q([0,x]) = f(1-x)$ for $x\in [0,1)$ and having an atom at $x=1$ of mass $Q(\{1\}) = 1- f(0)$ whenever $f(0) <1$. Finally, if $f= T(P,Q)$ then its generalized inverse is a trade off function and $f^{-1} =T(Q,P)$ \cite{dong2022gdp}[Lem A.2]
\begin{equation} \label{eq:geninv}
\text{ where } 
f^{-1} : [0,1] \to [0,1]
\text{ defined as } 
f^{-1}(\alpha)= \inf \{t \in [0,1]: f(t) \leq \alpha\} 
\text{ for }
\alpha \in [0,1].
\end{equation}

 \textbf{Functional ordering and indistinguishability (privacy).} One can define a \textit{functional ordering} 
 \begin{equation} \label{eq:Blackwellordering}
 \textbf{Blackwell ordering:}  
 \quad
 T(P',Q')(\alpha) \geq T(P,Q)(\alpha) \text{  for all } \alpha \in [0,1]
 \end{equation}
 capturing the intuition that the pair of distributions $(P',Q')$ defined on some space $(\cW', \cF_{\cW'})$ is more difficult to distinguish from a sample than the pair $(P,Q)$ defined on another space $(\cW, \cF_{\cW})$, uniformly across every scale $0\leq \alpha\leq 1$. We refer to  \cite{Torgersen1991Comparison} for a comprehensive treatment of equivalent conditions for this ordering also known as Blackwell-Le-Cam equivalence.
 
\textbf{Blackwell ordering.} Further, the functional ordering $f\geq g$ for $f,g \in \cT$ (a partial ordering on $\cT$ since not all pairs are comparable) described above is equivalent to the Blackwell ordering \cite{dong2022gdp}[Thm 2]. More precisely, $T(P',Q')(\alpha) \geq T(P,Q)(\alpha)$ for all $\alpha \in [0,1]$ if and only if there exists a Markov kernel $R: \cW \to  \cW'$ \cite{Kallenberg2021FMP3} such that $(P',Q')= (R(P),R(Q))$. 
\end{definition}

\subsection{Basics of $f$-differential privacy, motivation  of our work}

In this subsection, we introduce the f-differential privacy framework of \cite{dong2022gdp} based on Blackwell ordering and specialize to Gaussian differential privacy to motivate our generalization. Along the way, inspired by \cite{Torgersen1991Comparison}, we also state an equivalent way of describing the $f$-differential privacy of \cite{dong2022gdp} through their minimum Bayes risks.

\begin{definition}[$f$-differential privacy \cite{dong2022gdp}] \label{def:fDP}
    For $k \in \mathbb{Z}_{\geq 1}$, given an input space $(\cX, \cF_{\cX})$ and an output space $(\cY, \cF_{\cY})$, a  Markov kernel (randomized mechanism) $M:\mathcal{X}^k \to \mathcal{Y}$ is said to satisfy  $f$-DP  for some trade off function $f\in \cT$ if minimizing over  neighboring datasets $S,S' \in \cX^k$\footnote{Instead of the product metric space $(\cX^k, d_H)$ with the Hamming metric $d_{H}(\ul{x}, \ul{y}) = \sum_{i=1}^k \mathbf{1}(x_i\neq y_i)$, one can consider any metric space $(\cX, d)$ with the appropriate change in definitions. See \ref{thm:GDPdong} for more.}
\begin{equation}
    \inf_{d_H(S,S') \leq 1} T\big(M(S),M(S'))\big)\;\ge\; f \quad \text{pointwise on }[0,1].
\end{equation}

\end{definition}

\textbf{$k$ dependence, intuition, post processing.} For different values of $k$, one has different Markov kernels $M(k)$ and can allow different baseline trade off functions $f_k \in \cT$.  Let  $f=T(P, Q)$ for probability measures $P, Q$ on a space $(\cW, \cF_{\cW})$. Then the above definition captures the intuition that distinguishing the outputs of $M$ for any neighboring \(S, S' \in \cX^k\) is at least as hard as distinguishing between a sample of $P$ or $Q$. Moreover, if \(M:\cX^k \to \cY\) is \(f\)-DP and \(R:\cY \to \cZ\) is a Markov kernel independent of the data given \(M\), then \(R\circ M:\cX^k \to \cZ\) is also \(f\)-DP \cite{dong2022gdp}. 

\textbf{Convex duality and connections to $(\ve,\delta)$-DP and others.}  It follows from definition   that if a Markov kernel $M$ satisfy  $f_i$ DP for all $i\in I$, then it satisfy $\sup_{i\in I} f_i$ DP. Moreover,  \cite{wasserman2010statistical} established that satisfying $(\ve,\delta)$-DP in the classical sense is equivalent to satisfying $f_{\ve, \delta}$ -DP in the sense above, with $f_{\ve, \delta} (\alpha)= \max (0, 1- \delta  - e^{\ve} \alpha, e^{-\ve}(1-\delta -\alpha))$ for  $\alpha \in [0,1]$. As a consequence of the two above and convex duality\footnote{A convex, continuous  $f:[0,1] \to [0,1]$ is the pointwise supremum of all affine functions lying below it.},  for a \textit{symmetric} trade off function $f=f^{-1}$ we have a Markov kernel $M$ satisfy $f$-DP if and only if it satisfy $(\ve, \delta(\ve))$-DP in the classical sense, for all $\ve \geq 0$ with $\delta(\ve)= 1 + f^{*}(-e^{\ve})$, where $f^*(y)= \sup_{\alpha \in [0,1]} (\alpha y - f(\alpha))$ is the convex conjugate of $f$. We refer to \cite{dong2022gdp}[Appendix B] for the conversion of $f$-DP to divergence-based DP. Inspired by \cite{Torgersen1991Comparison}, one can derive a duality relation between the trade-off curve 
\begin{equation*}
f(\alpha)= T(P,Q)(\alpha) = \inf_{\varphi}\Bigl\{
        \e_{Q}[1-\vp] 
        \,\Bigm\vert\,
        \e_{P}[\vp]\le \alpha,\;
        \vp:(\cW, \cF_{\cW})\!\to[0,1] \text{ measurable}
      \Bigr\}, 
\end{equation*}

and the minimum Bayes risk for the testing problem $(P,Q)$ with prior $(1-\lambda, \lambda)$ for $\lambda \in [0,1]$
\begin{equation}
    b(\lambda)= B(P,Q)(\alpha)= \min_{\varphi} 
    \Bigl\{
         (1-\lambda) \e_P [\varphi] + \lambda \e_Q [1-\varphi]
        \,\Bigm\vert\
        \vp:(\cW, \cF_{\cW})\!\to[0,1] \text{ measurable}
      \Bigr\}.
\end{equation}
The precise duality relation for the pair  of functions $f,b:[0,1] \to [0,1]$ is given by\footnote{This representation writes out $f$ as a suprema over linear functions $f(\alpha)= \sup_{0 < \lambda \leq 1} \Big[ \frac{b(\lambda)}{\lambda}- \frac{(1-\lambda)}{\lambda} \alpha \Big], $} 
\begin{equation}
    f(\alpha)= \sup_{0 < \lambda \leq 1} \Bigg[ \frac{b(\lambda)- (1-\lambda)\alpha}{\lambda}\Bigg],
    \quad 
     b(\lambda)= \inf_{0\leq \alpha \leq 1} \Bigg[ (1-\lambda) \alpha + \lambda f(\alpha)\Bigg].
\end{equation}
 The equivalence between $f\leftrightarrow b$ gives rise to the an enumeration of all such functions $\{B(P,Q): [0,1] \to [0,1]: (P,Q) \text{ probability measures on some } (\Omega, \cF)\}$ as the following class.
\begin{equation}
    \cB:= \{ b: [0,1] \to [0,1] \quad | \quad  b \text{ concave and } b(\lambda) \leq \min(\lambda, 1-\lambda) \text{ for all }\lambda \in [0,1] \}
\end{equation}
Moreover, Blackwell ordering $f_1\geq f_2$ for $f_1, f_2 \in \cT$ can be equivalently stated in terms of their corresponding minimum Bayes risks  $b_1 \geq b_2$, pointwise on $[0,1]$. So, the entire framework of $f$-differential privacy proposed in \cite{dong2022gdp} can also be written equivalently in terms of the minimum Bayes risk function $b$. However, working with $b$ has the advantage that $b$ is affine under mixtures of experiments. More precisely, given pairs of binary (Blackwell-Le-Cam) experiments $\cE_1=(\Omega_1,\cF_1, P_1,Q_1)$ and $\cE_2=(\Omega_2,\cF_2, P_2,Q_2)$, consider the mixture experiment\footnote{The intuition is that in practice the mixture experiment is performed by deciding with probabilities $(\delta, 1-\delta)$ whether to perform a sample from $\cE_1$ or $\cE_2$ and then a sample is drawn. The output of the experiment has the form $(I,X_I)$, where $I$ is the index \{1, 2\} of which experiment is performed and $X_I$ is the sample itself.}
\begin{equation}
\delta \cE_1+ (1-\delta) \cE_2: = (\Omega_1 \bigsqcup \Omega_2, \sigma (\cF_1\sqcup \cF_2), \delta P_1 + (1-\delta) P_2,\delta Q_1 + (1-\delta) Q_2).
\end{equation}
\begin{equation}
\text{Then } 
B(\delta P_1 + (1-\delta) P_2, \delta P_1 + (1-\delta) P_2)\equiv \delta B(P_1, Q_1) + (1-\delta) B(P_2, Q_2).
\end{equation}This affineness of the functional $B$ is of importance in the sense that one can restrict the analysis of experiments that are only supported on $(\Omega, \cF) \equiv (\{0,1\}, 2^{\{0,1\}})$. It would be interesting to see the consequences of this in the practice of formalizing differential privacy through the $b$ function.

\textbf{$G_{\mu}$-DP and motivation of our work.} A central question in this hypothesis testing-based $f$-DP framework is to choose a statistically appropriate baseline trade-off function $f$ depending on the situation at hand. As a consequence of a central limit behavior, the focus of the paper \cite{dong2022gdp} was on the special case of a mechanism $M$ satisfying  $G_{\mu}$-DP by letting the baseline trade off function $f= G_{\mu}=T(P,Q)$ where $P= N\left(-\frac{\mu^2}{2}, \mu^2\right)$ and $Q= N\left(\frac{\mu^2}{2}, \mu^2\right)$ for some $\mu \in \mathbb{R}$. Our motivation is to search for a universal list of baseline TOFs $f=T(P,Q)$ beyond the Gaussian one.

\subsection{Composition and beyond the universal central limit behaviour}

Now, our interests are in quantifying how privacy degrades when we compose a sequence of mechanisms $M_1, \cdots, M_n$ on a (fixed) private data set $S\in \cX^k$ in which each analysis is informed by prior analyses on the same data set, and more and more analysis is publicly released with increasing $n$. For simplicity, we keep the $\sigma$-algebras and the measurability conditions implicit in definitions. 

\begin{definition}[Composition of private mechanisms \cite{dong2022gdp}, \cite{Kallenberg2021FMP3}]
Consider a Markov kernel $M_1:\cX^k \to \cY_1$, and conditional on the following pair (data, first output) = \((S,M_1(S))\), let $M_2: \cX^k \times \cY_1 \to \cY_2$ be another Markov kernel. The joint Markov kernel (randomized mechanism) $M= (M_1,M_2): \cX^k \to \cY_1 \times \cY_2$ is defined is $M(S)= (M_1(S), M_2(S, M_1(S)))$.
\end{definition}
\textbf{Intuition and $n$ step generalization.} Given the dataset $S$, the distribution of output of the joint mechanism $M(S)$ is constructed from the marginal distribution of $M_1(S)$ on $Y_1$ and the conditional distribution of $M_2(S, y_1)$ on $Y_2$ given $M_1(S) = y_1$. More generally, given a sequence of Markov kernels $M_i: \cX^k \times \cY_1 \times \cdots \times \cY_{i-1} \to \cY_i$, we can recursively define the joint mechanism as their composition $M = (M_1, \cdots, M_n): \cX^k \to \cY_1 \times \cdots \times \cY_n$. Now, with increasing $n$, we publicly release more and more analysis of the dataset $S$, first  $M_1(S)$ then $M_2(S, M_1(S))$, then $M_3(S, M_1(S), M_2(S, M_1(S)))$, and continue. The key question of  how the privacy degrades as the number of analyses increases, namely under composition leads us to the following definition.
\begin{definition}[Tensor product \cite{dong2022gdp}] \label{def:TPoTOF} Consider trade off functions $f= T(P,Q)$ for probability measures $(P,Q)$ on some space $(\cW, \cF_{\cW})$ and $g= T(P',Q')$ for probability measures $(P',Q')$ on some space $(\cW', \cF_{\cW'})$. We define the tensor product  trade off function as $f\otimes g := T(P\otimes P', Q\otimes Q')$, where $(P\otimes P', Q \otimes Q')$ are the product measures on $(\cW\times \cW', \cF_{\cW} \otimes \cF_{\cW'})$.
\end{definition}
\textbf{Well-definedness and basic properties.} By definition $f\otimes g$ is a trade-off function. \cite{dong2022gdp}[Lem C.2] showed that $f\otimes g$ is well defined in the sense that if $f= T(P,Q)= T(\Tilde{P}, \Tilde{Q})$, then $T(P\otimes P', Q\otimes Q')= T(\Tilde{P}\otimes P', \Tilde{Q}\otimes Q')$. Moreover, we have  \cite{dong2022gdp}[Prop D.1, D.2] 
\textit{monotonicity:} $g_1 \geq g_2$ implies $f\otimes g_1 \geq f\otimes g_2$,
\textit{commutativity:} $f\otimes g = g\otimes f$, \textit{associativity:} $(f\otimes g) \otimes h =  f \otimes (g\otimes h)$, \textit{identity:} $f\otimes I = I\otimes  f = f $ with $I(\alpha)=1- \alpha$ for $\alpha \in [0,1]$, \textit{inverse:} $(f\otimes g)^{-1} = f^{-1}\otimes g^{-1}$, with a \textit{trivial limit:} for any trade off function $\cT \ni f \neq I$, we have $\lim_{n\to \infty} f^{\otimes n}(\alpha)= 0$ for all $\alpha \in (0,1]$ as a consequence of \cite{polyanskiy_wu_2025}[Rem 7.6]. Now, we state the celebrated composition theorem that tightly quantifies  the privacy degradation of sequential mechanisms with the tensor product construction of trade off functions.

\begin{theorem}[Tightest composition theorem] \cite{dong2022gdp}[Thm 4]
    For $1\leq i \leq n$ consider Markov kernels $M_i$ as above so that  $M_i (\cdot, y_1, \cdots, y_{i-1}):\cX^k \to \cY_i$ satisfy $f_i$-DP for all $y_1 \in \cY_1, \cdots, y_{i-1} \in \cY_{i-1}$. Then the joint mechanism $M=(M_1,\cdots, M_n)$ satisfy  \(f_1\otimes \cdots \otimes f_n\)-DP.
\end{theorem}

Before stating a central limit behavior of nearly perfect private mechanisms as the number of compositions $n \uparrow \infty$, we need to define a few moment functionals of the log-likelihood
ratio (LLR) of $P$ and $Q$ such that $f = T(P,Q)$. The functionals $\text{kl}(f)$ $:=$ $ - \int_0^1 \log |f'(x)| dx,
\kappa_2(f):=$ $ \int_0^1 \log^2 |f'(x)| dx,
\kappa_3(f):=$ $ \int_0^1 \big|\log |f'(x)|\big|^3  dx,
\overline{\kappa}_3(f):=$ $ \int_0^1 \big|\log |f'(x)| + \text{kl}(f)\big|^3  dx$ (this one  requires $\text{kl}(f) <\infty$) take values in $[0,+\infty]$. Now, consider a triangular array of Markov kernels (randomized mechanisms)
$\{ M_{n1}, \ldots, M_{nn} \}_{n=1}^\infty$ where $M_{ni}$ satisfy  $f_{ni}$-DP for $1 \le i \le n$ \ref{def:fDP}.

\begin{theorem}\cite{dong2022gdp}[Thm 6] \label{thm:clt}
Let $\{ f_{ni} : 1 \le i \le n \}_{n=1}^\infty$ be a triangular array of
symmetric $(f=f^{-1})$ trade-off functions and assume for 
constants $k \ge 0$ and $s > 0$ as $n \to \infty$ we have 
$\sum_{i=1}^n \text{kl}(f_{ni}) \to k,  \max_{1 \le i \le n} \text{kl}(f_{ni}) \to 0, \sum_{i=1}^n \kappa_2(f_{ni}) \to s^2,\sum_{i=1}^n \kappa_3(f_{ni}) \to 0.$ Then, 
\begin{equation}
\lim_{n \to \infty}
f_{n1} \otimes  \cdots \otimes f_{nn}(\alpha)= G_{2k/s}(\alpha) \text{ uniformly for all } \alpha \in [0,1].
\end{equation}
\end{theorem}
\textbf{Limitation of $G_{\mu}$-DP and a non-Gaussian (Poisson) limit illustration.} \label{ex:poisson} The above theorem shows (under its assumptions) that the composition of nearly perfect DP mechanisms $M_{n1} \otimes \cdots \otimes  M_{nn}$  asymptotically satisfy $G_{2k/s}$-DP, where $G_{\mu}= T\left(N\left(-\frac{\mu^2}{2}, \mu^2\right), N\left(\frac{\mu^2}{2}, \mu^2\right)\right)$. But, consider the baseline  TOFs $ f_{ni}= T\left(\text{Ber}\left(\frac{\lambda_1}{n}\right), \text{Ber}\left(\frac{\lambda_2}{n}\right)\right)$ for $i \in [n]$, then it follows from the Neyman-Pearson lemma \eqref{eq:LRT} and \textit{law of small numbers} that uniformly for all  $\alpha \in [0,1]$ (see Corollary  \ref{cor:Poissonlimit})
\begin{align} \label{Poissonex}
&
f_{n1} \otimes \cdots \otimes f_{nn} (\alpha)
=
T\left(\text{Ber}\left(\frac{\lambda_1}{n}\right)^{\otimes n}, \text{Ber}\left(\frac{\lambda_2}{n}\right)^{\otimes n}\right) (\alpha)
\\
&
=
T\left(\text{Bin}\left(n,\frac{\lambda_1}{n}\right), \text{Bin}\left(n,\frac{\lambda_2}{n}\right)\right)(\alpha)
\to
T(P(\lambda_1), P(\lambda_2))(\alpha),
\end{align} 
where $P(\lambda_1), P(\lambda_2)$  denotes Poisson distributions with mean $\lambda_1$ and $\lambda_2$  for some $\lambda_1, \lambda_2 >0$. This example is clearly beyond the reach of the framework of Gaussian differential privacy, and depicts that based on the situation at hand, the baseline trade-off functions could be non-Gaussian. We unify Gaussian and Poisson limits under the framework of infinite divisibility \cite{JanssenMilbrodtStrasser1985}.

\textbf{On the $s^2=2k$ conjecture of \cite{dong2022gdp}.} It was remarked in \cite{dong2022gdp} that `in all examples of the application of Theorem \ref{thm:clt}' it is observed that $s^2= 2k$. We resolve this conjecture positively in our infinitely divisible setting as a consequence of a contiguity requirement.

\subsection{Infinitely divisible distributions and Lévy-Khintchine formula}
Now, we define infinitely divisible distributions, give some examples, state basic properties including closure under scaling and convolution, and weak limits, describe their appearances as distributional limits of sums of triangular arrays of (row-wise) IID random variables \cite{Kallenberg2021FMP3}.
\begin{definition}[Infinitely divisible distribution (IDD)]\label{def:IDD}
    A random variable $X$ on $\mathbb{R}$ is called infinitely divisible if $\forall$ $n \in \mathbb{Z}_{\geq 1}$, there exists  IID random variables $X_{n1}, \cdots, X_{nn}$ such that $X\overset{d}{=} \sum_{i=1}^n X_{ni}$. 
\end{definition}
\textbf{Equivalent definitions.} A probability measure $\mu \in \cP(\mathbb{R}, \cB_{\mathbb{R}})$ or its (equivalent) Fourier transform  $ \hat{\mu}(t)= \mu(e^{i t x})$ for $t\in \mathbb{R}$  is called infinitely divisible if for every $n\in \mathbb{Z}_{\geq 1}$, there exists a probability measure (which is also unique) $\mu_n \in \cP(\mathbb{R}, \cB_{\mathbb{R}}) $ such that $\mu= \mu_n^{*n}$ or equivalently $ \hat{\mu}(t)= \hat{\mu}_n(t)^n$ for $t\in \mathbb{R}$, where $\int_{\mathbb{R}} f(x) d\mu * \nu (x):= \int_{\mathbb{R}^2} f(x+y) d\mu(x) d\nu(y)$ defines the convolution of measures $\mu *\nu$ by varying $f$. We denote with $\cI$ the collection of all infinitely divisible distributions on $\mathbb{R}$. 

\textbf{Intuition, Gaussian and Compound Poisson examples in $\cI$.} The notion of infinitely divisible distribution captures the intuition that the random variable $X$ can be broken into smaller IID pieces of arbitrary order. Elements of $\cI$ include Normal distribution $\mu_1 \overset{d}{=} m + \sigma Z  \leftrightarrow \hat{\mu}_1(t)= e^{itm - \frac{t^2\sigma^2}{2}}, t\in \mathbb{R}$ for  some $m, \sigma \in \mathbb{R}$, Compound Poisson distribution $\mu_2 \overset{d}{=} \sum_{i=1}^{N}X_i$, with IID $X_i \sim \frac{\nu}{\|\nu\|}$   for $i \in \mathbb{Z}_{\geq 1}$ and independent $N \sim P(\|\nu\|) \perp \{X_i\}_{i\geq 1}$,   where $\nu$ is a (uniquely determined from $\mu_2$) positive finite Borel measure on $\mathbb{R}\setminus \{0\}$ of total mass $\|\nu\| = \nu(\mathbb{R})$ $\leftrightarrow \log \hat{\mu}_2(t)= \int_{\mathbb{R}} (e^{itx} - 1) d\nu(x)$. One can also represent the measure $\mu_2(\nu) = \sum_{k\in \mathbb{Z}_{\geq0}} e^{-\nu(\mathbb{R})} \frac{\nu^{*k}}{k!}$ as a Poisson mixture of measures.

\textbf{Basic closure properties of $\cI$.}  $\cI$ is closed under some algebraic properties including \textit{affine change:} $\cI \ni X$ implies  $ c_1X+c_2 \in \cI$ for $c_1, c_2 \in \mathbb{R}$, \textit{convolution:} $X_1 \perp X_2 \in \cI$ implies $ X=X_1 +X_2 \in \cI$. As a consequence we have an independent sum of Gaussian and Compound Poisson $X = m + \sigma Z + \sum_{i=1}^N X_i$ as described above, is also infinitely divisible.  $\cI$ is also closed under analytic properties including \textit{weak convergence:} $\cI \ni \mu_n \overset{d}{\to} \mu $ implies $\mu \in \cI$. Consequently, Poisson, Geometric, Negative Binomial, as well as Chi-square, Exponential, Gamma, and Cauchy are all elements of $\cI$. Moreover, a probability measure $\mu$ on $(\mathbb{R},\cB_{\mathbb{R}})$ is infinitely divisible if and only if there is a sequence of positive finite Borel measures $\{\nu_n\}_{n\in \mathbb{Z}_{\geq 1}}$ on $\mathbb{R} \setminus \{0\}$ so that $\mu_2(\nu_n) \overset{d}{\to} \mu$  \cite{Kallenberg2021FMP3}.

\textbf{Generalization of CLT and parametrizing limits as elements of IDD.} A remarkable generalization of the central limit theorem says  if $\mu \in \cP(\mathbb{R}, \cB_{\mathbb{R}})$ is such that there exists a triangular array of (row-wise) IID random variables $\{X_{ni}: 1 \leq i \leq n\}_{n \geq 1}$ such that $\sum_{i=1}^n X_{ni} \overset{d}{\to}\mu$, then $\mu \in \cI$ \cite{Durrett2019}, \cite{BoseDasGuptaRubin2002}. In fact, something more general holds, where one can drop the identically distributed condition \footnote{ We still keep the independence of elements of the row $\{X_{ni}: 1\leq i\leq n\}$ for every row.} on the triangular array, and replace it with a weaker condition on the triangular array also known as the \textit{uniformly asymptotically negligible} condition given by 
\begin{equation} \label{eq:UAN}
    \lim_{n\to \infty} \max_{1\leq i\leq n}\mathbb{P}[|X_{ni}| \geq \varepsilon ]= 0 \text{ for all } \varepsilon >0
\end{equation}
Even under such weaker row-wise independent triangular array, $\{X_{ni}: 1\leq i\leq n\}_{n\in \mathbb{Z}_{\geq 1}}$, if $\sum_{i=1}^nX_{ni} \overset{d}{\to} \mu$, then $\mu \in \cI$ \cite{Klenke2008Probability}. This enumeration of all possible limit distributions as elements of $\cI$ is also the key technical ingredient of our enumeration of all possible baseline trade-off functions.  We will establish our results again in two such stages like the above.

\textbf{Lévy-Khintchine representation of IDD.} Every IDD   $\mu$ on $\mathbb{R}$ has  a unique Lévy-Khintchine parametrization $(m, \sigma, \nu)$ of their Fourier transform with Gaussian and Compound Poisson parts, with  $\nu$ a positive Borel measure on $\mathbb{R}\setminus \{0\}$ satisfying $\int_\mathbb{R} (x^2\wedge 1) d\nu(x) < \infty$ \cite{Kallenberg2021FMP3}
\begin{equation}
    \log \hat{\mu}(t)=   i t m - \frac{1}{2} t^2 \sigma^2 + \int_{\mathbb{R}} (e^{itx} - 1 - itx \mathbf{1}(|x| \leq 1) ) d\nu(x).
\end{equation}

\textbf{Kolmogorov series representation of IDD.} Consider $m\in \mathbb{R}, \sigma \in \mathbb{R}_{\geq 0}$. Then one can define 
\begin{equation}
X= m + \sigma Z + X_0 +  \sum_{k \in \mathbb{Z}_{\geq 1}} (X_k - \mathbb{E}[X_k])
\end{equation}
a random series with independent random variables $Z\overset{d}{=} N(0,1)$,  $X_k \overset{d}{=}\mu_2(\nu_k)$ on a common probability space $(\Omega, \cF)$ where the positive finite Borel measures $\nu_k$ is supported on $I_k := \left( -\frac{1}{k}, -\frac{1}{k+1}\right ] \cup \left[\frac{1}{k+1}, \frac{1}{k}\right)$\footnote{We follow the convention that $\frac{1}{0}= \infty$ here to also include the case $k=0$ for $\nu_0$.}. As a consequence of Kolmogorov three series theorem, it is necessary and sufficient that $\sum_{k\in \mathbb{Z}_{\geq 1}} \mathbb{E}\big[\left(X_k -\mathbb{E}[X_k]\right)^2\big] < \infty$, for the infinite series representation to converge almost surely. Moreover, $\mathbb{E}[X_k]= \int x\nu_k (dx)$, $\mathbb{E}\big[\left(X_k -\mathbb{E}[X_k]\right)^2\big] = \int x^2 \nu_k(dx)$, and $\sum_{k\in \mathbb{Z}_{\geq 1}} \mathbb{E}\big[\left(X_k -\mathbb{E}[X_k]\right)^2\big] < \infty  \leftrightarrow \int x^{2} d \sum_{k\in \mathbb{Z}_{\geq 1}} \nu_k(x) <\infty \leftrightarrow \int (x^2 \wedge 1) d\nu(x) < \infty$  where $\nu = \sum_{k\in \mathbb{Z}_{\geq 0}} \nu_k$ is the canonical measure  of the infinitely divisible distribution  $\mu$ \cite{Klenke2008Probability}.

\subsection{Resolution of the $s^2=2k$ conjecture}

One of the major conclusions of the \cite{dong2022gdp} paper was that  under the limit of large number of compositions of private mechanisms, a natural baseline trade-off function is the Gaussian trade-off function $f=T\left(N\left(-\frac{\mu^2}{2}, \mu^2\right), N\left(\frac{\mu^2}{2}, \mu^2\right)\right)$ for some $\mu \in \mathbb{R}$. But, as we saw earlier in \ref{ex:poisson} that this framework does not capture the Poisson example $f= T(P(\lambda_1), P(\lambda_2))$. We show that under the limit of large number of compositions of private mechanisms, an enumeration of all the baseline trade-off functions (encompassing both Gaussian and Poisson case) is given by the following.
\begin{equation} \label{eq:IDP}
    \cI_T= \{f= T(P,Q): P \in \cI \text{ and } dQ(x) = e^{x}dP(x)\}
\end{equation}
 
\begin{theorem} \label{thm:GDP}
    Consider the infinitely divisible distribution $P \overset{d}{=} -k + |s|Z$ with only Gaussian component $Z\sim N(0,1),k,s\in \mathbb{R}$. If $dQ(x) = e^x dP(x)$ is a probability measure, then $s^2=2k$. More precisely, for a Gaussian $P \in \cI$ so that its Esscher tilt $dQ(x)= e^{x} dP(x)$ is a probability measure and $T(P,Q) \in \cI_{f}$  if and only if $P=N\left(-\frac{\mu^2}{2}, \mu^2\right)$.
\end{theorem} 
\begin{proof}
    The proof follows immediately from a computation of the moment-generating function of a Gaussian distribution, saying that $1=Q(\mathbb{R})=\mathbb{E}_P[\exp(x)] = \exp\left({-k + \frac{s^2}{2}}\right) \rightarrow s^2= 2k$.
\end{proof}
Now, in the following lemma, we describe the basic properties of trade-off functions in general and show that the class of infinitely divisible trade-off functions $\cI_{\cT}$ that captures both the Gaussian $T\left(N(0,1), N(\mu,1)\right)$  and the Poisson $T(P(\lambda_1), P(\lambda_2))$ trade-off functions.
    \begin{lemma} Closure of trade-off functions under bijections and likelihood-ratios. \label{lem: TOFB}
        \begin{enumerate}
            \item For two random variables defined $X,Y$ on some probability space $(\Omega, \cF, \mathbb{P})$ and constants $ 0 \neq c_1 , c_2 \in \mathbb{R}$ we have  $T(X,Y) = T(c_1X +c_2, c_1Y +c_2)$, where we denote for any random elements $(\Omega, \cF, \mathbb{P}) \overset{A,B}{\to} (\cX, \cF_{\cX})$ trade-off function  $T(A,B):= T(P_A, P_B)$ with $P_A = \mathbb{P} \circ  A^{-1}$, $P_B= \mathbb{P} \circ B^{-1}$ are the corresponding distributions. More generally, for any bijective measurable map $(\cX, \cF_{\cX}) \overset{h}{\to} (\cY, \cF_{\cY})$ we have $T(A,B)= T(h(A), h(B))$.
            \item For any two probability distributions $P,Q$ on some measurable space $(\cX, \cF_{\cX})$ we have $T(P,Q)= T(L_P,L_Q)$, where $L_P := P \circ \left(\log \frac{dQ}{dP}\right)^{-1}, L_Q:= Q \circ \left(\log \frac{dQ}{dP}\right)^{-1}$ are the distributions of the likelihood ratio  $\log \frac{dQ}{dP}$ under their respective probability measures.
            \item \textbf{Gaussian} We observe that if $P \overset{d}{=} N\left(-\frac{\mu^2}{2}, \mu^2\right)$, then  $dQ(x)= e^{x}dP(x) \overset{d}{=}  N\left(\frac{\mu^2}{2}, \mu^2\right)$. So, the Gaussian trade-off curve as observed in \cite{dong2022gdp} belongs to $\cI_T$.
            \item \textbf{Poisson.} The Poisson curve $T(P(\lambda_1), P(\lambda_2)) \in \cI_T$, since $T(P,Q)=T(P(\lambda_1), P(\lambda_2)$, where $P\overset{d}{=} \lambda_1- \lambda_2 + N \log \left( \frac{\lambda_2}{\lambda_1}\right)$, with $N\sim P(\lambda_1)$ and $dQ(x)= e^{x} dP(x)$.
        \end{enumerate}
    \end{lemma}
    \begin{proof}
    \begin{enumerate}
        \item Observe that from Blackwell's theorem \cite{dong2022gdp} post-processing always increases the trade-off function. Therefore, for any measurable functions $(\cX, \cF_{\cX}) \overset{h}{\to} (\cY, \cF_{\cY})$, we have $T(h(A),h(B)) \geq T(A,B)$. Moreover, if $h$ is bijective, reversing the process above yields $T(h(A), h(B)) \geq T(h^{-1}(h(A)), h^{-1}(h(B))) = T(A,B)$. Now, for scalar (or even vector) valued random variables, scaling $(x \to c_1 x, c_1\neq 0)$ and translations  $(x \to x + c_2)$ are measurable bijections, and hence the equality.
        \item  Let $X\sim P$, $Y\sim Q$ for random variables $X,Y$ on some probability space $(\Omega, \cF, \mathbb{P})$. Let $L = \log \frac{dQ}{dP}$ the likelihood ratio, considered as a random variable $(\cX, \cF_{\cX}) \to \mathbb{R}$. Again, by post–processing, $T\big(L(X),L(Y)\big) \geq T(X,Y) =T(P,Q)$. Now, for the reverse inequality,  we apply the Neyman–Pearson likelihood ratio test \eqref{eq:LRT} for the pair $(P,Q)$ of size \(\alpha\).

        \begin{equation*}
            \varphi^*(w) :=\mathbf{1}\left( \log \left(\frac{dQ}{dP}\right) > \log \tau^* \right) + \lambda^* \mathbf{1}\left(\log \left(\frac{dQ}{dP}\right)=\log \tau^*\right), L =  \log \left(\frac{dQ}{dP}\right)
        \end{equation*}
 with $\frac{dQ}{dP}: = \frac{dQ}{d\mu}/\frac{dP}{d\mu}:(\cX, \cF_{\cX}) \to [0, \infty]$, where the thresholds $\tau^*\in [0,\infty], \lambda^* \in [0,1]$ are chosen as in \eqref{eq:LRT} so that
\(\mathbb{E}_{P}[\varphi^*]=\alpha\). Type I and type II errors are given by ($t^*= \log \tau^*$)
\begin{equation}
\alpha= \mathbb{E}_{P}[\varphi^*]
= \mathbb{E}_{P}\!\big[\mathbf{1}\{L>t^*\}+\lambda^*\,\mathbf{1}\{L=t^*\}\big]
= P\!\big[L>t^*\big]
  + \lambda^*P \!\big[L=t^*\big],
\end{equation}
\begin{equation}
T(X,Y)(\alpha)=T(P,Q)(\alpha)= 1-\mathbb{E}_{Q}[\varphi^*]
= 1-{Q}\!\big[L>t^*\big] - \lambda^*Q\big[L=t^*\}\big].
\end{equation}
Now, consider the trade-off function  between  \(P_X \overset{d}{=}L(X)\)  and \(P_Y\overset{d}{=}L(Y)\) with the  test  $\varphi^*(x)= \mathbf{1}\{x>t^*\}+\lambda^*\,\mathbf{1}\{x=t^*\}$ to  have the type $I$ error $\mathbb{E}_{P_{X}}[ \varphi^*]= \alpha$ and type II error $T(L(X), L(Y))(\alpha) \leq 1- \mathbb{E}_{P_Y}[\varphi^*]= T(X,Y)(\alpha)$. Since $\alpha \in [0,1]$ was arbitrary, we have $T(X,Y) =T(L(X), L(Y))$. Combining both inequalities proves the claim.
\item Observe that a shifted and scaled Gaussian $P= N\left(-\frac{\mu^2}{2}, \mu^2\right) \overset{d}{=}  -\frac{\mu^2}{2} + |\mu|Z$ is indeed infinitely divisible and its Esscher tilt $\mathbb{Q}[\mathbb{R}]=\mathbb{E}_{P}[\exp(x)]= \exp\left( -\frac{\mu^2}{2} + \frac{\mu^2}{2}\right)=1$, and therefore $Q$ is a probability measure. Moreover, the moment generating function  (MGF)
\begin{align*}
&
\mathbb{E}_{Q}[\exp(tx)]
=
\mathbb{E}_{P}[\exp((t+1)x)]
= 
\exp\left( (t+1)\left(-\frac{\mu^2}{2} + \frac{1}{2} (t+1) \mu^2\right)\right) 
\\
&
=
\exp\left( \frac{t\mu^2}{2} + \frac{t^2\mu^2}{2} \right)
= 
\mathbb{E}\big[\exp\left(t\left( \frac{\mu^2}{2} + |\mu|Z\right)\right)\big] \text{ for } t\in \mathbb{R}.
\end{align*}

Since the moment generating function of $Q$ matches that of a Gaussian, we have $Q\overset{d}{=} \frac{\mu^2}{2} + |\mu| Z$. Now the equality $T(P,Q)= T(Z, Z+\mu)=T(N(0,1),N(\mu,1))$ for $\mu>0$ follows from the first part of the lemma about the invariance of trade-off curves under (common) translations and scaling, and additional translation symmetries of the Gaussian tradeoff functions $T(Z+c_1, Z+c_2)= T(Z, Z+c_2-c_1)$ for $c_1,c_2 \in \mathbb{R}$.
\begin{equation}
    T(P,Q) = T\left(-\frac{\mu^2}{2} + |\mu|Z, \frac{\mu^2}{2} + |\mu|Z\right) = T\left(-\frac{\mu}{2} + Z, \frac{\mu}{2} + Z\right) = T(Z, Z+\mu). 
\end{equation}
\item Consider a shifted and scaled Poisson $P\overset{d}{=} \lambda_1- \lambda_2 + N_1\log \left( \frac{\lambda_2}{\lambda_1}\right)$, with $N_1\sim P(\lambda_1)$, and  therefore infinitely divisible and $Q[\mathbb{R}]= \mathbb{E}_{P}[\exp(x)]= \exp\left(\lambda_1 - \lambda_2 + \lambda_1\left(\frac{\lambda_2}{\lambda_1}-1\right)\right)  = 1$, therefore  $Q$ is a probability measure. Moreover, the moment generating function  (MGF) 
\begin{align*}
&
\mathbb{E}_Q[\exp(tx)]
= \mathbb{E}_{P}[\exp((t+1)x))] 
= \exp\left((t+1) (\lambda_1 - \lambda_2)  + \lambda_1\left(\left(\frac{\lambda_2}{\lambda_1}\right)^{t+1}-1\right)\right) 
\\
&
=
\exp\left(t (\lambda_1 - \lambda_2)  + \lambda_2\left(\left(\frac{\lambda_2}{\lambda_1}\right)^{t}-1\right)\right) 
=
\mathbb{E}\big[ \exp\left(t \left(\lambda_1- \lambda_2 + N_2\log \left( \frac{\lambda_2}{\lambda_1}\right)\right)\right)\big]
\end{align*}
where $N_2\sim P(\lambda_2)$. Since the moment generating function of $Q$ matches that of  a scaled and shifted Poisson, we have $Q\overset{d}{=} \lambda_1- \lambda_2 + N_2\log \left( \frac{\lambda_2}{\lambda_1}\right)$.  Now the equality $T(P,Q)= T(P(\lambda_1),P(\lambda_2))$ for $\lambda_1, \lambda_2>0$ follows from the first part of the lemma about the invariance of trade-off curves under (common) translations and scaling\footnote{Observe that equality holds if $\lambda_1= \lambda_2 >0$, and we do not need to apply any invariance properties.}.
\begin{equation*}
    T(P,Q)
    =
    T\left(\lambda_1- \lambda_2 + N_1\log \left( \frac{\lambda_2}{\lambda_1}\right), \lambda_1- \lambda_2 + N_2\log \left( \frac{\lambda_2}{\lambda_1}\right)\right) 
    =
    T(N_1, N_2) 
\end{equation*}
\end{enumerate}
    \end{proof}
    Now, in the following lemma, we describe the basic properties of infinitely divisible trade-off functions, from which a closure  of $\cI_{\cT}$ under pointwise limits $\cI_{\cT} \supset f_n \to f \in \cI_T$ can also be inferred.
    \begin{lemma} \label{lem:TOF1}
    Basic properties of the infinitely divisible trade-off functions. 
    \begin{enumerate}
        \item \textbf{Tensor product.} Consider $P_1,P_2 \in \cI$ with their corresponding Esscher tilts $Q_1,Q_2$. Then the tensor product $f= T(P_1\otimes P_2, Q_1\otimes Q_2) =T(P,Q) \in \cI_T$ with $P= P_1 * P_2$ and $dQ= e^x dP$. This immediately extends to $(P_i,Q_i)_{i=1}^n$ for $P_i\in \cI$ for any $n\in \mathbb{Z}_{\geq 1}$ as well. 
        \item \textbf{Infinite divisibility.} For any trade off function $f=T(P,Q) \in \cI_T$ and any $n \in \mathbb{Z}_{\geq 1}$ there exists a trade off function $f_n = T(P_n,Q_n) \in \cI_T$ such that $f= f_n^{\otimes n}$.
    \end{enumerate}
\end{lemma}
    \begin{proof}
    \begin{enumerate}
        \item 
    Observe that  the likelihood ratio of the tensor product pairs of distributions satisfies
\begin{equation}L_{P_1\otimes P_2}:= (P_1 \otimes P_2) \circ \left(\log \frac{d Q_1\otimes Q_2}{d P_1\otimes P_2}\right)^{-1} =P_1*P_2 =:P, \text{ and}
\end{equation}
\begin{equation}
 dL_{Q_1\otimes Q_2} =d(Q_1 \otimes Q_2) \circ \log \left(\frac{d Q_1\otimes Q_2}{d P_1\otimes dP_2}\right)^{-1} = e^x d P_1*P_2 = e^x dP =:dQ.
\end{equation}

Now, using the lemma above, we have $T(P_1\otimes P_2, Q_1\otimes Q_2) = T(L_{P_1\otimes P_2}, L_{Q_1\otimes Q_2}) = T(P,Q)$. Finally, observe that $P_1*P_2$ is infinitely divisible if $P_1$ and $P_2$ are \cite{Kallenberg2021FMP3} since the addition of two independent infinitely divisible random variables is also an infinitely divisible random variable. Moreover, $\mathbb{E}_{P_1*P_2} [\exp(x)]= \mathbb{E}_{P_1} [\exp(x)] \mathbb{E}_{P_2} [\exp(x)] =1$ satisfies the normalization requirement.
\item 
Observe that by the definition of infinite divisibility of $P \overset{d}{=} X$  satisfying $\mathbb{E}_{P}[\exp(x)] = \mathbb{E}[\exp(X)]=1$, for every $n \in \mathbb{Z}_{\geq 1}$ we have $P= P_n^{*n}$ for some probability measure $P_n$ on $(\mathbb{R}, \cB_{\mathbb{R}})$ such that $1=\mathbb{E}_P[\exp(x)] = \left(\mathbb{E}_{\mathbb{P}_n}[\exp(x)]\right)^n$. Then, from the positivity of the quantity $\mathbb{E}_{\mathbb{P}_n}[\exp(x)] >0$, we have $\mathbb{E}_{\mathbb{P}_n}[\exp(x)]=1$, and so $dQ_n = e^x dP_n$ is a probability measure.  Moreover, from the Levy-Khintchine characterization $(m, \sigma, \nu)$ of the characteristic function $\hat{P}(t)$ of an infinitely divisible distribution $P$, it follows that $P_n$ itself is an infinitely divisible distribution with its Levy-Khintchine triplet $\left(\frac{m}{n}, \frac{\sigma}{\sqrt{n}}, \frac{\nu}{n}\right)$. Now, $f_n^{\otimes n} = T(P_n^{\otimes n}, Q_n^{\otimes n})= T(P,Q)$, where $P=P_n^{*n}$, and $dQ= e^{x}dP$.
\end{enumerate}
\end{proof}
\section{Technical results}

In this section, we prove the main technical result of the paper, establishing why is the collection of infinitely divisible trade-off functions $\cI_T$ \eqref{eq:IDP} appears under the limit of a large number of compositions of (nearly perfect) differentially private mechanisms, as a converse of the result \ref{lem:TOF1}.
\begin{theorem} \label{thm:IDP}
    Consider a sequence of trade off functions $\{f_n\}_{n\in \mathbb{Z}_{\geq 1}}$ such that $f_n(0)=1$ for all $n\in \mathbb{Z}_{\geq 1}$ and $f_n^{\otimes n}(\alpha) \to f_{\infty}(\alpha)$ pointwise on $[0,1]$ for some trade-off function $f_{\infty}$ \eqref{eq:TOF}, then $f_{\infty}$ is an infinitely divisible trade off function. More precisely $f_{\infty}\in \cI_T$.
\end{theorem}
\begin{proof}
We prove it in a few steps. First, recall that any trade-off function $f\in \cT$ can be realized as $f=T(P,Q)$ on  the measurable space $(\cW,  \cF_{\cW})=$ $([0,1],\cB_{[0,1]})$ by taking $P$ to be the standard uniform probability measure $U$ on $([0,1],\cB_{[0,1]})$ and $Q=Q_f$ to be a Borel probability measure on $[0,1]$ with a cumulative distribution function $Q([0,x]) = f(1-x)$ for $x\in [0,1)$ and having an atom at $x=1$ of mass $Q(\{1\}) = 1- f(0)$ whenever $f(0) <1$\footnote{This choice is called the canonical pairs of probability measures $(U,Q_f)$ for $f$ \cite{Torgersen1991Comparison}.}. So, for tradeoff functions $f_n$ with $f_n(0)=1$ we have $f_n= T(P, Q_n)$, where $P = U$ as above and $Q_n := Q_{f_n} \ll P$\footnote{For two probability measures $Q \ll P$ on $(\cX, \cF_{\cX})$ means that for all $A\in \cF_{\cX}$ if $P(A)= 0 \rightarrow Q(A)=0$ or equivalently for all measurable functions  $\varphi(\cX,\cF_{\cX}) \to [0,1]$ we have $\mathbb{E}_{P}(\varphi) = 0 \rightarrow \mathbb{E}_Q(\varphi)=0$.} since $f_n(0)=1$\footnote{This follows from the fact that by definition of $f_n \in \cT$ \eqref{eq:TOF} we have  $Q_n$ is absolutely continuous with respect to $U$ on $[0,1)$ with Radon-Nikodym density $\frac{dQ_n}{dP}(x)= -f'_n(1-x)$, $P$ almost everywhere on  $[0,1)$, ($f_n$ is convex and therefore Lebesgue almost differentiable \cite{Rockafellar1970ConvexAnalysis}).  Now $f_n(0)=1$ implies that this absolute continuity $Q_n\ll U$ extends to the entire closed interval $[0,1]$.}. More precisely,  $f_n(0)=1$ is equivalent to the fact  $Q_n \ll P$\footnote{$Q_n \ll P\rightarrow $ `for measurable $(\cX, \cF_{\cX}) \overset{\varphi}{\to}[0,1]$ we have $\mathbb{E}_P(\varphi)= 0 \rightarrow \mathbb{E}_{Q_n}(\varphi) = 0 \rightarrow f_n(0)=1$.}. Therefore, $f_n(0)=1$ implies $Q_n$ has no atom or, more precisely, no singular component with respect to $P=U$.
    
Second, $f_n^{\otimes n} = T\left(U^{\otimes n}, Q_n^{\otimes n}\right)$ converges to $f_{\infty}= T(P_{\infty},Q_{\infty})$ pointwise (hence uniformly \cite{dong2022gdp},\cite{lehmann1986testing})\footnote{The intuition is that if a collection of monotone functions $f_n :[0,1] \to [0,1] $ converges pointwise to a continuous monotone function $f: [0,1] \to [0,1]$, then the convergence is uniform. Moreover, the domain (in this case $[0,1]$) need not be compact, only the range  (in this case $[0,1]$) needs to be. The version, where trade-off functions $\{(f_n)_{n\in \mathbb{Z}_{\geq 1}}, f\}$ are replaced by cumulative distributive function $\{(F_n)_{n\in \mathbb{Z}_{\geq 1}}, F\}$ with $F$ continuous, this is known as the Glivenko-Cantelli lemma or the fundamental theorem of statistics, whose finite sample concentration  is  known as the Dvoretzky–Kiefer–Wolfowitz inequality
 \cite{Kallenberg2021FMP3}.} for some pairs of probability measures $(P_{\infty}, Q_{\infty})$ on  a measurable space $(\Omega, \cF)$\footnote{Observe that we are not taking the canonical choice for $f$, because we believe the infinite divisible parametrization $\cI_{\cT}$ is more useful for expressing limiting trade-off functions and practically more informative in the sense that will be once we demonstrate the procedures on how to achieve such privacy requirements.}. Now, $f_n(0)=1$ for all $n \in \mathbb{Z}_{\geq 1} \rightarrow f_{\infty}(0)=1$ which is equivalent to  $Q_{\infty}\ll P_{\infty}$. This is an equivalent reformulation of  contiguity  for  $(U^{\otimes n}, Q_n^{\otimes n})$ \cite{LeCam1986}.
 
Third, the pointwise convergence of trade-off functions for the sequence  $\left(U^{\otimes n}, Q_{n}^{\otimes n}\right)$ to the trade-off function for the pair $ (P_{\infty},Q_{\infty})$ is equivalent to weak  convergence of the sequence on $\mathbb{R}$ \footnote{This equivalence of weak convergence of `Blackwell experiments' for the two point set $\Theta = \{0,1\}$ and different equivalences has been established in detail in \cite{LeCam1986}, \cite{Torgersen1991Comparison}. The version with weak convergence of the likelihood ratio (under the null) to the likelihood ratio under the limiting null is ubiquitous in statistics, and the pointwise convergence of the trade-off functions to a limiting trade-off function can be exactly written as `Le Cam distance' between $(U^{\otimes n}, Q_n^{\otimes n})$ and $(P,Q)$ going to zero.} 
\begin{equation}
U^{\otimes n} \circ  \left(\log \frac{d Q_n^{\otimes n}}{d U^{\otimes n}}\right)^{-1}   \overset{d}{\to} P_{\infty} \circ \left(\log \frac{dQ_{\infty}}{dP_{\infty}}\right)^{-1} \overset{d}{=} L(P_{\infty})
\end{equation}

Now, we observe that on the canonical probability space $\left([0,1]^{n}, \cB_{[0,1]^n},U^{\otimes n} \right)$ the random variable $\log \left(\frac{dQ_n^{\otimes n}}{dU^{\otimes n}}(x_1, \cdots, x_n)\right) = \sum_{i=1}^n \log \left( \frac{dQ_n}{dU}(x_i)\right)$ is a sum of independent and identically distributed random variables. Moreover, by assumption of the theorem, this IID sum converges in distribution to a probability distribution $L(P_{\infty})$. Now, by the fundamental theorem of infinitely divisible distributions\footnote{The fundamental theorem of infinitely divisible distributions \cite{Durrett2019} says that if a probability measure $\mu$ on $(\mathbb{R}, \cB_{\mathbb{R}})$ is such that there exists a sequence of probability measures $\{\mu_n\}_{n\in \mathbb{Z}_{\geq 1}}$ on $(\mathbb{R}, \cB_{\mathbb{R}})$ satisfying weak convergence under convolutions $ \mu_n^{*n} \overset{d}{\to}\mu$, then $\mu$ is infinitely divisible.}  $L(P_{\infty})$ is infinitely divisible, and moreover $T(P_{\infty},Q_{\infty})= T(L(P_{\infty}), L(Q_{\infty}))$ where $d L(Q_{\infty})(x)= e^{x} dL(P_{\infty})(x)$. Hence, $f_{\infty}=T(P_{\infty},Q_{\infty})\in \cI_T$ is infinitely divisible.
\end{proof}
\textbf{A Gaussian example.} Consider the case $f_n(\cdot)= T(N(0,1), N(\mu_n,1))(\cdot)= \Phi(\Phi^{-1}(1-\cdot)- \mu_n)$ for $\mu_n >0$.  Then $f_n(0)=1$, and by the dimension-freeness of Gaussian trade-off functions\ref{lem: GTOF}, we have $f_n^{\otimes n}(\cdot)= T(N(0,1), N(\mu_n \sqrt{n}, 1)) =\Phi(\Phi^{-1}(1-\cdot)- \mu_n\sqrt{n})$. Therefore,  the pointwise convergence of $f_n^{\otimes n}$ to a limiting trade-off function $f\in \cT$ immediately implies that the sequence  $\{\mu_n \sqrt{n}: n\in \mathbb{Z}_{\geq 1}\}$ has to converge, and moreover, the $\lim_{n\to \infty} \mu_n \sqrt{n} = \mu_{\infty}$   has to be finite, because the  trade-off function $f_{\infty}$ has to have the form $f_{\infty}(\cdot)= \Phi(\Phi^{-1}(1-\cdot)- \mu_{\infty})$, and if $\mu_{\infty}=\infty$, then $f_{\infty}\equiv 0$\footnote{Observe that for any $\alpha >0$ we have $f_{\infty}(\alpha)=\Phi(\Phi^{-1}(1-\alpha)- \mu_{\infty})=0$ whenever $\mu_{\infty}=\infty$, and therefore the continuity requirement on $f$ throughout  the closed interval on $[0,1]$ implies that $f_{\infty}(0)=0.$}, a contradiction. We compute non-Gaussian examples in section  \ref{ssec:practice}

\textbf{The Contiguity assumption in Theorem \ref{thm:IDP}.} The assumption of  $f_n(0)=1$ for all $n\in \mathbb{Z}_{\geq 1}$ in Theorem \ref{thm:IDP} is maintained so that the canonical pair of probability measures $P=U$ and $Q_n:= Q_{f_n}$ satisfy absolute continuity $Q_n \ll U$ for all $n\in \mathbb{Z}_{\geq 1}$. Consequently, our results  are not directly applicable to the convergence of trade-off functions involving $f_{\varepsilon_n, \delta_n}$ with $\delta_n >0$ for all $n\in \mathbb{Z}_{\geq 1}$. However, the tensor product decomposition $f_{\varepsilon, \delta} =f_{\varepsilon, 0} \otimes f_{0,\delta}$ as observed in \cite{dong2022gdp} immediately allows us to apply our results directly to the $f_n=f_{\varepsilon_n,0}$ part and treat the $f_n=f_{0, \delta_n}$ part separately, as has been done in \cite{dong2022gdp} using the identity $f_{0,\delta_1} \otimes f_{0,\delta_2} = f_{0, 1- (1-\delta_1)(1-\delta_2)}$. Moreover, as  \cite{dong2022gdp} established the limits of such tensor products 
\begin{equation}
f_{0,\frac{\delta}{n}}^{\otimes n} \to  f_{\infty}=f_{0, 1- e^{-\delta}} \text{ where } f_{0, 1-e^{-\delta}}(0) = 1-\delta \neq 1 \text{ whenever }\delta >0.
\end{equation}
So $f_{\infty}=f_{0, 1- e^{-\delta}}$ is not infinitely divisible according to our definition of $\cI_{\cT}$. On a different note, we believe that our conclusions of Theorem \ref{thm:IDP} will continue to hold under the weaker assumption that $f_n(0)$ need not to be equal to $1$, but we require $f_{n}^{\otimes n}(0) \to 1$ as $n\to \infty$. 

\textbf{The Convergence assumption in Theorem \ref{thm:IDP} .}
The assumption $g_n(\alpha):=f_n^{\otimes n}(\alpha) \to f_{\infty}(\alpha)$ was made pointwise on $[0,1]$ for some trade-off function $f$ \eqref{eq:TOF} in the statement of Theorem \ref{thm:IDP}. This is natural since our primary goal was to  enlarge  the class of Gaussian trade-off functions obtained in \cite{dong2022gdp} and to see what class of trade-off functions can be obtained in general under the composition of a large number of nearly perfect differentially private operations, under the assumption that these degradation curves $g_n:= f_n^{\otimes n}$ do converge to a limiting trade-off curve. However, in practice of differential privacy through the $f$-differential privacy framework, one has to equivalently show (under the assumptions of our Theorem \ref{thm:IDP}) that if $f_n=T(P_n,Q_n)$ with $f_n(0)=1$, then $f_n^{\otimes n} \to f_{\infty}$ for some limiting trade-off function, $f_{\infty}= T(P_\infty, Q_{\infty})$   is the same as the following (possibly more tractable) distributional convergence of the likelihood ratio under the null
\begin{equation} \label{eq:practice}
\text{ under } P_{n}  \quad \log \frac{dQ_n^{\otimes n}}{dP_n^{\otimes n}} \overset{d}{\to} X\overset{d}{=} \log \frac{dQ_{\infty}}{dP_{\infty}} \text{ under } P_{\infty}.
\end{equation}

In section \ref{ssec:practice}, we show the distributional convergence of $\log \frac{dQ_n^{\otimes n}}{dP_n^{\otimes n}}$  under $P_{n}$ for various statistical models, including the locally asymptotically normal, Poisson and Compound Poisson families.

\textbf{Blackwell-Le-Cam equivalence.} Observe that our conditions of Theorem \ref{thm:IDP} only involve trade-off functions $\{f_n: n\in \mathbb{Z}_{\geq 1}\}, f$, and it is deliberate since the definition of $f$-differential privacy of a mechanism \ref{def:fDP} does not involve any details of the binary experiment $\cE= (\Omega, \cF, P,Q)$, except that it works for any such equivalent representation of the baseline trade-off function $f$ through $f= T(P,Q)$. This relates to the Blackwell-Le Cam equivalence of experiments, which we will briefly  describe for binary experiments with a parameter set $\Theta=\{0,1\}$ consisting of two elements. 

Given two binary experiments $\cE_1= (\Omega_1, \cF_1, P_1,Q_1)$, $\cE_2= (\Omega_2, \cF_2, P_2,Q_2)$, a notion of distance between $\cE_1, \cE_2$ (better known as \textit{Le-Cam deficiency}) was defined in \cite{LeCam1986} that is given as 
\begin{equation}  \label{eq:Lecamdist}
\textbf{Le-Cam distance: }
\Delta(\cE_1, \cE_2)= \max(\delta(\cE_1, \cE_2), \delta(\cE_2,\cE_1)),   \text{ where}
\end{equation}
\begin{equation} \label{eq:Lecamdef}
    \textbf{Le-Cam deficiency: } \delta(\cE_1, \cE_2) = \inf_{R: (\Omega_1, \cF_1) \to (\Omega_2, \cF_2)} \max(\|RP_1- P_2\|_1, \|RQ_1- Q_2\|_1), 
\end{equation}
 where the infimum is taken over all (appropriate)  Markov kernels $R$ and for a pair of probability measures $\mu_1, \mu_2$ on $(\Omega, \cF_2)$, $\|\mu_1- \mu_2\|_1:= \sup_{A \in \cF_2} |\mu_1(A)- \mu_2(A)|$ is the total variation distance. 

 One of the key results of \cite{LeCam1986}, \cite{Torgersen1991Comparison} is that $\Delta(\cE_1, \cE_2) =0$ if and only if $f_1:= T(P_1,Q_1) \equiv T(P_2, Q_2)=:f_2$. Moreover, $f_2 \geq f_1$ if and only if $\delta(\cE_1, \cE_2) =0$ if and only if there exists a Markov kernel $R :(\Omega_1, \cF_1) \to (\Omega_2, \cF_2)$ such that $RP_1=P_2$ and  $RQ_1= Q_2.$ 
 Furthermore, the correspondence $\cE = (\Omega, \cF, P, \cQ) \leftrightarrow f=T(P,Q)$ preserves the Le-Cam distance in the following sense and therefore is quantitative, since \cite{Torgersen1991Comparison}[Corollary 9.3.27] states
 \begin{equation}  \label{LC=TOFD}
     \sqrt{2}\Delta(\cE_1, \cE_2)= \Lambda(Q^{(1)}, Q^{(2)}), \text{ where the }
 \end{equation}
  cumulative distribution functions $F^{(1)}, F^{(2)}$ of probability measures $Q^{(1)},Q^{(2)}$ on $[0,1]$ are given
 \begin{equation}
 \text{by }Q^{(i)}[0,\alpha]= 1-f_i(\alpha) \text{ for } \alpha \in (0,1]  \text{ and }  Q^{(i)}[\{0\}] =(1-f_i(0)), 
  \text{ for }i \in \{1,2\}
  \end{equation}
 and $\Lambda(Q^{(1)},Q^{(2)})$ denotes the Levy distance between the probability measures $Q^{(1)},Q^{(2)}$ through their cumulative distribution functions  $F^{(1)}, F^{(2)}$ on $[0,1]$ as \cite{Kallenberg2021FMP3} 
\begin{equation}
\Lambda(F^{(1)}, F^{(2)})= \inf \{\varepsilon >0: F^{(2)}(x-\varepsilon) - \varepsilon \leq F^{(1)}(x) \leq F^{(2)}(x+\varepsilon)+ \varepsilon \text{ for all }x \in [\varepsilon,1-\varepsilon]\}.
\end{equation}
The key observations about the Levy-distance $\Lambda(\cdot, \cdot)$ (hence $\Delta$ for binary $\Theta$) is that it precisely captures convergence in distribution. More precisely, our assumption in Theorem \ref{thm:IDP} about the  pointwise convergence of $f_n^{\otimes n} \to f_{\infty}$  on $[0,1]$ \footnote{The assumption $f_{n}(0)=1$ (hence $f_{\infty}(0)=1$) is there to ensure that $Q_{\infty}\ll P_{\infty}$ for $f_{\infty}=T(P_{\infty}, Q_{\infty}).$} is the equivalent to the condition that $\Delta(1- f_n^{\otimes n}, 1-f_{\infty}) \to 0$, and therefore captures the meaning of convergence of experiments $\Delta(\cE_n^{\otimes n}, \cE_{\infty}) \to 0$ where the correspondence $\cE_n \equiv (P_n,Q_n) \to f_n= T(P_n, Q_n)$ and $\cE_{\infty} \equiv (P_n,Q_n)\to f_{\infty} \equiv T(P_{\infty},Q_{\infty})$.

From the definition of $\Delta(\cE_1, \cE_2)$, given a fixed parameter set $\Theta= \{0,1\}$ (in our case) one can (consistently) define a equivalence class  of experiments $[\cE]:=\{\cE_1: \Delta(\cE,\cE_1)=0\}$  with the same trade-off curve $f \leftrightarrow [\cE]$. Now, in the $f$-differential privacy framework of \cite{dong2022gdp}, it only depends on the equivalence classes of experiments and not the exact representation, such as $f=T(P,Q)$. However, this explicit representation $f=T(P,Q)$ is important in practice since,  one has to implement a Markov kernel (a generalized noise-adding mechanism) such as adding Gaussian noise to achieve Gaussian differential privacy $f= T(N(0,1) N(\mu,1))$.

This is the reason our explicit parametrization of infinitely divisible trade-off functions $f_{\infty}=T(P_{\infty}, Q_{\infty})$ of $\cI_{\cT}$  is important, and we will demonstrate in section \ref{sec:BIDP} how to
 achieve infinitely divisible privacy for a large class of infinitely divisible trade-off functions $f_{\infty}=T(P_{\infty}, Q_{\infty})$, including the practically relevant Poisson case $f_{\infty}= T(P(\lambda_1), P(\lambda_2))$.

 \textbf{Large deviation regime versus infinitely divisible regime.}  
Instead of a sequence $\{f_n: n\in \mathbb{Z}_{\geq 1}\}$ of trade-off functions and considering their tensor product $f_{n}^{\otimes n}$, if one takes a single  trade-off function $f$ (independent of $n$) and considers $f^{\otimes n}$, one has a large deviation principle 
\begin{equation} \label{eq:LDP}
\lim_{n \to \infty}\frac{1}{n} \log f^{\otimes n} (\alpha) = \int_{0}^{z_f} \log |f'(x)| dx, \text{ for all } \alpha \in (0,1)
\end{equation}

where $z_{f}= f^{-1}(0)$ is the first zero of $f$, where $f^{-1}$ is the generalized inverse of $f$ \eqref{eq:geninv} (see\cite{dong2022gdp}[Prop B.4], \cite{polyanskiy_wu_2025}, \cite{Torgersen1991Comparison}). Although our focus is on the infinitely divisible (central limit) regime, where there is some non-trivial privacy left even after a large number of nearly perfect private operations. The large deviation regime is also of interest and has been used in the theoretical development of watermark detection for LLMs \cite{LiRuanWangLongSu2025AOS}.

\textbf{Some consequences in practice.} One can observe that, as a consequence of the continuous mapping theorem\footnote{To be mathematically precise, we need the \textit{uniform integrability} condition on the sequence of random variables $\{\left(\frac{dQ_n}{dU}\right)\}_{n \in \mathbb{Z}_{\geq 1}}$  under $U^{\otimes n}$ \cite{Kallenberg2021FMP3} so that we can use $X_n \overset{d}{\to}X \implies \e[X_n] \to \e[X]$.} the mean of the distribution $L(P)$ can be expressed in terms of the sequence $\{f_n\}_{n \in \mathbb{Z}_{\geq 1}}$.
\begin{equation}
\mathbb{E}_{L(P)}[x]= \mathbb{E}_{P}\big[\log \left( \frac{dQ}{dP}\right)\big]
=
\lim_{n\to \infty} \sum_{i=1}^n \mathbb{E}_{U} \big[ \log \left(\frac{dQ_n}{dU}\right)\big]
=
\lim_{n\to \infty} n \int_{0}^{1} \left( \log \left(- f_n'(1-x)\right)\right) dx.
\end{equation}
Moreover, the cumulant generating function (CGF) of the distribution $L(P)$ for $t\in [0,1]$\footnote{As pointed out in \cite{Torgersen1991Comparison} the CGF might not be well defined outside $t\in [0,1]$.} can be expressed through the trade-off function sequence as well $\{f_n\}_{n \in \mathbb{Z}_{\geq 1}}$\footnote{$(P,Q) \to \mathbb{E}_{P}\big[\left(\frac{dQ}{dP}\right)^t\big]$ is also known as the Renyi-divergence or Hellinger transform of order $t$ \cite{polyanskiy_wu_2025}. In Le-cam's theory of weak convergence of experiments for the two-point set $\Theta= \{0,1\}$, convergence of the Hellinger transforms of all orders $t\in [0,1]$ for the pair $(U^{\otimes n}, Q_{n}^{\otimes n}) $ to that of the limiting pair $(P,Q)$ is another equivalent conditions along with the weak convergence of the likelihood ratios under the null, and the Le-cam deficiency metric $\Delta(\cE_n, \cE)$ between $\cE_n=(U^{\otimes n}, Q_{n}^{\otimes n})$ and $\cE=(P,Q)$ going to zero.} .
\begin{equation}
\log \mathbb{E}_{L(P)}[\exp(tx)]
= 
\log \mathbb{E}_{P}\big[\left(\frac{dQ}{dP}\right)^t\big]
=
\lim_{n\to \infty} n \log \left(\int_{0}^{1} \left(-f_{n}'(1-x)\right)^t dx\right)\big]
\end{equation}

As a consequence, one can differentiate both sides with respect to $t$ and derive an expression for moments of $L(P)$ in terms of the sequence of $\{f_n\}_{n\in \mathbb{Z}_{\geq 1}}$ \cite{NicaSpeicher2006FreeProbability}.

\textbf{Triangular systems of differentially private procedures.} Now, we rephrase the implications of the above statement in the \(f\)-differential privacy framework, where we consider triangular arrays of trade-off functions $\{f_{ni}: 1\leq i\leq n\}_{n\in \mathbb{Z}_{\geq 1}}$ with $f_{ni}(0)=1$ for all $1\leq i\leq n \in \mathbb{Z}_{\geq 1}$. We consider the tensor product  $f_{n1} \otimes \cdots \otimes f_{nn}= T(U^{\otimes n}, Q_{n1}\otimes \cdots \otimes  Q_{nn})$, where $U$ is the uniform probability measure on $([0,1], \cB_{[0,1]})$ and $Q_{ni}[0,x]= f_{ni}(1-x)$ satisfying $Q_{ni}\ll U$.  Now, if we assume that $f_{n1}\otimes \cdots \otimes f_{nn}$ converges pointwise to $f= T(P,Q)$ for some trade-off function $f$\footnote{Pointwise convergence of the sequence of trade-off functions $\{f_n\}$ with $f_n(0)=1$ requires $f(0)=1$ and therefore $Q\ll P$, a contiguity statement in disguise for the sequence $\left( U^{\otimes n}, Q_{n1}\otimes \cdots \otimes  Q_{nn}\right)$. In our case, $f_n = f_{n1}\otimes \cdots f_{nn} =T(U^{\otimes n}, Q_{n1}\otimes \cdots \otimes Q_{nn})$, and because $f_{ni}(0)=1 \leftrightarrow Q_{ni}\ll U \rightarrow Q_{n1}\otimes \cdots \otimes Q_{nn} \ll U^{\otimes n}$ (absolute continuity is preserved under tensor  products \cite{Folland1999RealAnalysis}) $\leftrightarrow f_n(0)=1$.}, by \cite{LeCam1986, Torgersen1991Comparison} , it becomes equivalent to the weak convergence of the likelihood ratio 
\begin{equation}
    U^{\otimes n} \circ  \left(\log \frac{d Q_{n1}\otimes  \cdots \otimes dQ_{nn}}{d U^{\otimes n}}\right)^{-1}   \overset{d}{\to} P \circ \left(\log \frac{dQ}{dP}\right)^{-1} \overset{d}{=} L(P).
\end{equation}
Now, the likelihood ratio  $\log \left(\frac{d Q_{n1}\otimes  \cdots \otimes dQ_{nn}}{d U \otimes \cdots \otimes dU}(x_1, \cdots, x_n)\right) = \sum_{i=1}^n \log \left( \frac{dQ}{dU}(x_i)\right)$ on the probability space $\left([0,1]^{n}, \cB_{[0,1]^n},U^{\otimes n} \right)$ is a random variable that is a sum of independent random variables $\{\log \left( \frac{dQ_{ni}}{dU}\right)\}_{i=1}^n$ under $U^{\otimes n}$. This holds true for every $n \in \mathbb{Z}_{\geq 1}$.  Now, under the technical assumption of a uniformly asymptotically negligible condition,
\begin{equation} \label{eq:uanllr}
\lim_{n\to \infty} \max_{1 \leq i\leq n} U\left(x \in [0,1] :| \log \left( \frac{dQ_{ni}}{dU}(x)\right)| >\varepsilon\right) = 0 \text{ for all } \varepsilon >0,
\end{equation}
the limiting distribution $L(P)$ is still an IDD with its Esscher tilt distribution $dL(Q)= e^xdL(P)(x)$. One can observe that, as a consequence of the continuous mapping theorem\footnote{To be mathematically precise we need the \textit{uniform integrability} conditions analogous to the IID case.}, the mean $\mathbb{E}_{L(P)}[x]$ of the distribution $L(P)$ can be expressed in terms of the sequence $\{f_{ni}: 1\leq i \leq n\}_{n \in \mathbb{Z}_{\geq 1}}$.
\begin{equation}
\mathbb{E}_{P}\big[\log \left( \frac{dQ}{dP}\right)\big]
=
\lim_{n\to \infty} \sum_{i=1}^n \mathbb{E}_{U} \big[ \log \left(\frac{dQ_{ni}}{dU}\right)\big]
=
\lim_{n\to \infty}  \sum_{i=1}^n \int_{0}^{1} \left( \log \left(- f_{ni}'(1-x)\right)\right) dx.
\end{equation}
Moreover, the cumulant generating function (CGF) of the distribution $L(P)$ for $t\in [0,1]$ can be expressed through the trade-off function sequence as well $\{f_{ni}: 1\leq i\leq n\}_{n \in \mathbb{Z}_{\geq 1}}$\footnote{$(P,Q) \to \mathbb{E}_{P}\big[\left(\frac{dQ}{dP}\right)^t\big]$ is also known as the Renyi-divergence or Hellinger transform of order $t$ \cite{polyanskiy_wu_2025}. In Le-cam's theory of weak convergence of experiments for the two-point set $\Theta= \{0,1\}$, convergence of the Hellinger transforms of all orders $t\in [0,1]$ for the pair $(U^{\otimes n}, Q_{n}^{\otimes n}) $ to that of the limiting pair $(P,Q)$ is another equivalent conditions along with the weak convergence of the likelihood ratios under the null, and the Le-cam deficiency metric $\Delta(\cE_n, \cE)$ between $\cE_n=(U^{\otimes n}, Q_{n}^{\otimes n})$ and $\cE=(P,Q)$ going to zero.} .
\begin{equation}
\log \mathbb{E}_{L(P)}[\exp(tx)]
= 
\log \mathbb{E}_{P}\big[\left(\frac{dQ}{dP}\right)^t\big]
=
\lim_{n\to \infty} \sum_{i=1}^n \log \left(\int_{0}^{1} \left(-f_{ni}'(1-x)\right)^t dx\right)\big]
\end{equation}

As a consequence, one can differentiate both sides with respect to $t$ and derive an expression for moments (cumulants) of $L(P)$ in terms of the sequence of $\{f_n\}_{n\in \mathbb{Z}_{\geq 1}}$ \cite{NicaSpeicher2006FreeProbability}.
Moreover, in the above condition  \eqref{eq:uanllr} if one inserts the maximum inside the probability, and demands $ \lim_{n\to  \infty} U\left(x\in [0,1]| \max_{1\leq i\leq n} | \log \left( \frac{dQ_{ni}}{dU}(x)\right)| >\varepsilon\right) = 0 \text{ for all } \varepsilon >0,$ then this condition is known as a version of the Lindeberg-type condition and causes the limit distribution $L(P)$ to be Gaussian, in our case, of the form $L(P)\overset{d}{=} - \frac{s^2}{2} + |s|Z$ for $Z \sim N(0,1)$ \cite{LeCamYang2000}.

\textbf{Connections to high dimensional hypothesis testing.}
Observe that in the $f$-differential privacy framework of \cite{dong2022gdp},  while sequentially applying more than one  differentially private operation:  first $M_1$ (satisfying $f_1$ DP), then $M_2$ (satisfying $f_2$ DP)\footnote{This includes adaptive operations, meaning the action of $M_2$ not only depends on the original dataset $S$, but also on the output of $M_1(S)$. So, the joint output is $M(S)= (M_1(S), M_2(S,M_1(S)))$  on a dataset $S$.} it has been captured through tensor products of experiments $\cE_1 = (\Omega_1, \cF_1, P_1, Q_1)$ with $f_1= T(P_1, Q_1)$ and $\cE_2= (\Omega_2, \cF_2, P_2, Q_2)$ with $f_2= T(P_2,Q_2)$ so that $M= (M_1, M_2)$ satisfy $f=f_1\otimes f_2= T(P_1\otimes P_2, Q_2\otimes Q_2)$ DP. 

However, an inspection of the proof of the composition theorem \cite{dong2022gdp}[Lemma C.3] reveals that the condition on $M_2$ satisfying $f_2$ DP independent of the exact outcome of $M_1$ on the dataset $S$ gives rise to the tensor product structure\footnote{Since this essentially assumes the worst case privacy loss for the second mechanism $M_2$ as $f_2$ uniformly over all outputs  of the first mechanism $M_1$, it does not allow $f_2$ to depend on the exact output $M_1(S)$.}.  In principle, it is not unreasonable to consider a more flexible situation  so that the degradation of privacy upon applying multiple procedures $M^{(n)}= (M_1,M_2, \cdots, M_n)$ (on the same dataset $S$) might occur in the more general form where each $M^{(n)}$ satisfy $g^{(n)}$-DP where $g^{(n)}= T(P^{(n)},Q^{(n)})$ with the experiment $\cE^{(n)} = (\Omega^{(n)} , \cF^{(n)}, P^{(n)}, Q^{(n)})$.  Observe that, in the previous worst-case analysis, these experiments had the product structure 
\begin{equation} \label{eq:Prodexp}
\cE^{(n)}
=
\cE_1 \times \cdots \times \cE_n
=
(\Omega_1 \times \cdots \times \Omega_n, \cF_1 \otimes \cdots \otimes \cF_n, P_1\otimes \cdots  \otimes P_n, Q_1\otimes \cdots \otimes  Q_n).
\end{equation}
This product structure is precisely what led us to the infinitely divisible structure of our limiting trade-off functions $\cI_{\cT}$. However, for a general such experiment  $\cE^{(n)} = (\Omega^{(n)} , \cF^{(n)}, P^{(n)}, Q^{(n)})$ capturing privacy degradation for $M^{(n)}= (M_1, \cdots, M_n)$ more adaptively than the tensor product framework of \cite{dong2022gdp}, it is difficult in practice to establish the pointwise convergence of the sequence of trade-off functions $g^{(n)}= T(P^{(n)},Q^{(n)}) \to g_{\infty}= T(P_{\infty}, Q_{\infty})$ and to be able to guess  the structure of the class of such limiting trade-off functions $g_{\infty}= T(P_{\infty}, Q_{\infty})$ one might expect in these more general situations. Equivalently, it is indeed more difficult to prove distributional convergence for the likelihood ratio under the null hypothesis. 
\begin{equation} \label{:HDT}
 \text{ under } P^{(n)}  \quad  \log \frac{dQ^{(n)}}{dP^{(n)}} \overset{d}{\to} X \overset{d}{=} \log \frac{dQ_{\infty}}{dP_{\infty}} \text{ under } P_{\infty}.
\end{equation}
Proving the above distributional convergence as $n \to \infty$ is a central task in the problem of high dimensional hypothesis testing, which has turned out to be of increasing interest in recent times \cite{HanJiangShen2023Contiguity},  \cite{JohnstoneOnatski2020SpikedModels}, \cite{OnatskiMoreiraHallin2013SphericityPower}. Among other results, distributional convergence of likelihood ratio under the null to a contiguous pair of Gaussians has been shown in some restricted cases, such as below the celebrated BBP threshold \cite{BaikBenArousPeche2005BBP} for spiked covariance models. However, it will be extremely interesting to analyze what kind of limiting binary experiments or trade-off curves one obtains, as we believe the limiting trade-off functions might belong to a different universality class beyond the infinitely divisible class of trade-off function, with the arise of new universality classes such as the Tracy-Widom distributions \cite{ElKaroui2007TWLimitSampleCovariance}\footnote{More precisely, Tracy-Widom distributions appear as distributional limits of the largest eigenvalues (after appropriate translation and scaling) of a spiked covariance model. as a test statistic, but not as the limit of the likelihood ratio, so it still does not rule out an infinitely divisible limit for the likelihood ratio under the null.}.

\subsection{Some consequences within and beyond infinite divisibility} \label{ssec:practice}
In our Theorem \ref{thm:IDP}, we used the results of \cite{LeCam1986}, \cite{Torgersen1991Comparison} to show if for a sequence of functions $f_n\in \cT$, $f_n(0)=1$, we have  $f_n^{\otimes n} \to f_{\infty}$ converge pointwise on $[0,1]$ for some function $f_{\infty} \in \cT$, then the corresponding canonical likelihood ratio $\log \frac{dQ_n^{\otimes n}}{dU^{\otimes n}}$  under $U^{\otimes n}$  converges in distribution to a limiting ratio $\log \frac{dQ_{\infty}}{dP_{\infty}}$ under $P_{\infty}$ for some $Q_{\infty} \ll P_{\infty}$ so that $f_{\infty}= T(P_{\infty}, Q_{\infty})$. By the fundamental theorem of infinite divisibility, the ratio under $P_{\infty}$ is infinitely divisible. 

\textbf{Application for statistical models.} Consider a statistical model $\{P_{\theta}:\theta \in \Theta\}$, a collection of probability measures $\{P_{\theta}:\theta \in \Theta\}$ on some measurable space $(\Omega, \cF)$ with $f_{n}= T(P_{\theta_n}, P_{\theta_n+\delta_n})$\footnote{Assume the base points $\theta_n \in \Theta $ is such that the local perturbation $\theta_n + \delta_n \in \Theta$ for all $n\in \mathbb{Z}_{\geq 1}$.} with $f_n(0)=1 \leftrightarrow $ $P_{\theta_n +\delta_n} \ll P_{\theta_n}$, and consider $f_n^{\otimes n} = T(P_{\theta_n}^{\otimes n}, P_{\theta_n+\delta_n}^{\otimes n})$. Now, again, by \cite{LeCam1986}, \cite{Torgersen1991Comparison}  to establish the pointwise convergence of $f_n^{\otimes n} \to f_{\infty} \in \cT$ it is enough to establish that  a distributional limit of the likelihood ratio under the null $P_{\theta_n}^{n}$  (on $\Omega^{\times n}, \cF^{\otimes n}$)

\begin{equation} \label{eq:statmodel}
\log\frac{dP_{\theta_0 +\delta_n}^{\otimes n}}{dP_{\theta_n}}(w_1, \cdots, w_n) = \sum_{i=1}^n  \log \frac{dP_{\theta_n +\delta_n}}{dP_{\theta_n}}(w_i) \overset{d}{\to} X  \text{ satisfying } \e[\exp(X)] =1
\end{equation}
for some random variable $X:(\Omega, \cF, \mathbb{P}) \to \mathbb{R}$. Then one defines on $(\mathbb{R}, \cB_{\mathbb{R}}$) the distributions  $P_{\infty} \overset{d}{=}X$\footnote{Since $X$ is a distributional limit of IID sums of random variables $\{\log \frac{dP_{\theta_n +\delta_n}}{dP_{\theta_n}}\}$, it is infinitely divisible.} and $Q_{\infty}(h) = \e[h\exp (X)]$ for bounded measurable $h: (\mathbb{R}, \cB_{\mathbb{R}}) \to (\mathbb{R}, \cB_{\mathbb{R}})$ or equivalently $dQ_{\infty} (x) =e^{x}dP_{\infty}(x)$. Observe that the distributional limit assumption along with its normalization is precisely Le-Cam's first contiguity lemma for $P_{\theta_{n}+\delta_n}^{\otimes n} \triangleleft P_{\theta_n}^{\otimes n}$ \cite{vandervaart1998asymptotic}[Lemma 6.4].

The fact that $f_{\infty}= T(P_{\infty}, Q_{\infty})$ follows from the following version of Le-Cam's result, combined with the fact that for binary experiments, the convergence of Le-Cam distance for binary experiments is the same as the pointwise convergence of their corresponding trade-off functions \eqref{LC=TOFD}.

\begin{lemma}[\cite{pollard2011lecam}[Lemma 1] \label{lem:pollard}
Consider  binary experiments $\cE_n = (\Omega^{(n)},\cF^{(n)},P^{(n)}, Q^{(n)})$ with $Q^{(n)} \ll P^{(n)}$ and likelihood ratio $X_{n} =dQ^{(n)}/dP^{(n)}$, and let
$\cE_{\infty}=(\Omega_{\infty},\cF_{\infty},P_{\infty}, Q_{\infty})$ with $Q_{\infty} \ll P_{\infty}$ and likelihood ratio $X_{\infty} =dQ_{\infty}/dP_{\infty}$ be another. Suppose the random variable
\begin{equation} \label{eq:pollard1}
(\text{under } P^{(n)}) \quad  X_n= dQ^{(n)}/dP^{(n)}\overset{d}{\to} Y= dQ_{\infty}/dP_{\infty} (\text{ under } P_{\infty}). \text{ Then }
\end{equation} 
 $\lim_{n\to \infty} \Delta(\cE^{(n)}, \cE_{\infty}) =0 $  or equivalently $ f_{n}:= T(P^{(n)},Q^{(n)}) \to  f_{\infty} = T(P_{\infty}, Q_{\infty})$ pointwise.
\end{lemma}

\textbf{A dominating extension.} The conclusion $f_n=T(P^{(n)},Q^{(n)}) \to T(P_{\infty}, Q_{\infty})$ in the lemma above  continues to hold when we keep the assumption $Q^{(n)} \ll P^{(n)}$, and there exists a $\sigma$-finite\footnote{A positive measure $\mu$ on $(\Omega, \cF)$ is called $\sigma$-finite if there is a sequence  $\{S_k: k\in \mathbb{Z}_{\geq 1}\}$ of  increasing $S_k \subset S_{k+1}$ measurable sets $S_k \in \cF$ of finite measure $\mu(S_k) <\infty$ for all $k\in \mathbb{Z}_{\geq 1}$ so that $\cup_{k\in \mathbb{Z}_{\geq 1} S_k = \Omega}$ Our examples of $\mu$ will be counting measures on at most countable set \cite{Kallenberg2021FMP3}} measure $\mu^{(n)}$ on $(\Omega^{(n)}, \cF^{(n)})$ dominating both $P^{(n)},Q^{(n)}\ll \mu^{(n)}$ with densities $p^{(n)}, q^{(n)}$  respectively, and a $\sigma$-finite measure  $\mu_{\infty}$ on $(\Omega_{\infty}, \cF_{\infty})$ dominating $P_{\infty}, Q_{\infty} \ll \mu_{\infty}$ with density $p_{\infty}, q_{\infty}$ with $Q_{\infty} \ll P_{\infty}$, but instead of  \eqref{eq:pollard1} we suppose the following distributional convergence of the densities

\begin{equation} \label{eq:pollard2}
(\text{under } \mu^{(n)}) \quad (p^{(n)},q^{(n)})\overset{d}{\to} (p_{\infty}, q_{\infty}) \quad (\text{ under } \mu_{\infty}).
\end{equation}
 
As a consequence of the above result, we immediately obtain results for the celebrated Locally Asymptotically Normal family, the proofs of which can be found in \cite{vandervaart1998asymptotic}.

\begin{definition}[\cite{vandervaart1998asymptotic}]
A statistical model $(P_\theta : \theta \in \Theta \subset \mathbb{R}^k)$ is called
\textit{differentiable in quadratic mean at $\theta_0$} if there exists a
measurable vector-valued function $\dot\ell_{\theta_0}$ such that, as
$\theta \to \theta_0$,
\begin{equation} \label{eq:DQM}
  \int \Big[
    \sqrt{p_\theta}
    - \sqrt{p_{\theta_0}}
    - \tfrac{1}{2} (\theta - \theta_0)^\top \dot\ell_{\theta_0}
      \sqrt{p_{\theta_0}}
  \Big]^2 \, d\mu
  = o\big(\|\theta - \theta_0\|^2\big), 
\end{equation}
where we assume $p_{\theta}= \frac{dP_{\theta}}{d\mu}$ are the densities with-respect to a dominating measure $\mu$ on the $(\Omega, \cF)$.
\end{definition}

\begin{theorem}\cite{vandervaart1998asymptotic}[Theorem 7.2] \label{cor:Gausslimit}
Suppose that $\Theta\subset \mathbb{R}^k$ open and  the
model $(P_\theta : \theta \in \Theta)$ is differentiable in quadratic mean
at $\theta$. Then $P_\theta \dot\ell_\theta = 0$ and the Fisher information
matrix $I_\theta = P_\theta \dot\ell_\theta \dot\ell_\theta^{\mathsf T}$ exists. Furthermore, for every converging sequence $h_n \to h$\footnote{We are specializing to the case of our more general model \eqref{eq:statmodel} with $\theta_n=\theta$ and $\delta_n = h_n/\sqrt{n}$.}, as
$n \to \infty$,
\begin{equation}
 \text{under } P_{\theta}^{\otimes n} 
 \quad  \log \prod_{i=1}^n
  \frac{p_{\theta + h_n / \sqrt{n}}(X_i)}{p_\theta(X_i)}
  \overset{d}{\to}
  X
  \overset{d}{=}
  N\left(-\frac{1}{2} h^{\mathsf T} I_\theta h, h^{\mathsf T} I_\theta h \right)
  =P_{\infty}
\end{equation}
Consequently, $f_n^{\otimes n}:=T(P_{\theta}^{\otimes n}, P_{\theta + h_n/\sqrt{n}}^{\otimes n}) \to f_{\infty}= T(P_{\infty}, Q_{\infty})$, where $dQ_{\infty}(x)= e^{x} dP_{\infty}(x)$.
\end{theorem}

A large class of statistical models \cite{vandervaart1998asymptotic}, including most exponential families with density $p_{\theta}(x)= d(\theta) h(x) \exp\left( Q(\theta)^Tt(x)\right)$ for appropriate choices of $d(\cdot), h(\cdot), Q(\cdot), t(\cdot)$, location models with density $p_{\theta}(x)= f(x-\theta)$ for appropriate $f(\cdot)$ satisfies the conclusion of the theorem.

Now, we describe the Poisson case, which is beyond Gaussian, but  still belongs to the infinitely divisible class $\cI_{\cT}$. It is one of the key motivating examples for unifying with the Gaussian theory.

\begin{corollary}[Poisson limits]\label{cor:Poissonlimit}
Consider a Bernoulli model $(P_{\theta}\equiv\mathrm{Ber}(\theta): 0\leq \theta \leq 1)$ on  $(\Omega, \cF)= (\{0,1\}, 2^{\{0,1\}})$ so that $P_{\theta}(1)= \theta$ and $P_{\theta}(0)= 1- \theta$. Now, consider $P_{n}= P_{\frac{\lambda_1}{n}}$ and $Q_{n}= P_{\frac{\lambda_2}{n}}$ for some $\lambda_1, \lambda_2 >0$\footnote{We are specializing to the case of our more general model \eqref{eq:statmodel} with $\theta_n=\lambda_1/n$ and $\delta_n = (\lambda_2 -\lambda_1)/n$.} and denote $f_n=T(P_n, Q_n)$. Then $f_n^{\otimes n} \to f_{\infty}= T(P(\lambda_1), P(\lambda_2))$, where $P(\lambda)$ denotes the Poisson distribution  on $(\mathbb{Z}_{\geq 0}, 2^{\mathbb{Z}_{\geq 0}}), P(\lambda) (n) := e^{-\lambda} \frac{\lambda^n}{n!}$ for all $n\in \mathbb{Z}_{\geq 0}$.

\begin{proof}
    The proof follows from a \textit{sufficiency} argument. First, observe from the likelihood ratio test  
    \begin{equation}
    f_n^{\otimes n} =T (P_n^{\otimes n}, Q_n^{\otimes n})= T(B(n, \lambda_1/n), B(n, \lambda_2/n)), \text{ where } B(n, \theta) \text{ for } \theta \in [0,1]
    \end{equation}
  is the binomial distribution on $(\{0, \cdots, n\}, 2^{\{0, \cdots, n\}})$, $B(n,\theta)(k)= \binom{n}{k}\theta^{k} (1-\theta)^{k}$ for $k\in  \{0, \cdots, n\}$. This equality follows from lemma \ref{lem: TOFB}  $ T(c_1X + c_2, c_1Y+c_2)= T(X,Y)$ and
  \begin{equation}
      T(P,Q)= T(L_P,L_Q) 
      \text{ where }
      L_P= P\circ \left( \log \frac{dQ}{dP}\right)^{-1} 
      \text{ and }
      L_Q= Q\circ \left( \log \frac{dQ}{dP}\right)^{-1}.
  \end{equation}
  \begin{equation}
      P_n^{\otimes n}\circ \left(\log(\frac{dQ_n^{\otimes n}}{dP_n^{\otimes n}}\right)^{-1} = n \log \left(\frac{1-\lambda_2/n}{1-\lambda_1/n}\right) + B(n, \lambda_1/n) \log \left(\frac{\lambda_2(1-\lambda_1/n)}{(1-\lambda_2/n)\lambda_1}\right)
  \end{equation}
  \begin{equation}
      Q_n^{\otimes n}\circ \left(\log(\frac{dQ_n^{\otimes n}}{dP_n^{\otimes n}}\right)^{-1} = n \log \left(\frac{1-\lambda_2/n}{1-\lambda_1/n}\right) + B(n, \lambda_2/n) \log \left(\frac{\lambda_2(1-\lambda_1/n)}{(1-\lambda_2/n)\lambda_1}\right)
  \end{equation}
 From $B(n, \lambda/n) \overset{d}{\to} P(\lambda)$ \cite{Kallenberg2021FMP3} and lemma \ref{lem:pollard} we  have for  $N_1 \overset{d}{=} P(\lambda_1), N_2\overset{d}{=} P(\lambda_2)$
  \begin{equation}
      T(P_n^{\otimes n}, Q_n^{\otimes n}) \to T\left(\lambda_1 -\lambda_2 + N_1 \log\left(\frac{\lambda_2}{\lambda_1}\right), \lambda_1 -\lambda_2 + N_2 \log\left(\frac{\lambda_2}{\lambda_1}\right)\right) = T(N_1, N_2)
  \end{equation}
\end{proof}
\textbf{The triangular case.}  Consider the inhomogeneous products of Bernoulli model $P^{(n)} = \otimes_{i=1}^n P_{n,i},$  with  $P_{n,i}= \mathrm{Ber}(p_{n,i})$ and $Q^{(n)} = \otimes_{i=1}^n Q_{n,i}$ with $Q_{n,i}= \mathrm{Ber}(q_{n,i})$ where  
\begin{equation}
0\leq p_{n,i}, q_{n,i} \leq 1,
\quad 
\max_{1\leq i \leq n}(p_{n,i} \vee q_{n,i}) \to 0, 
\text{ and }
\quad 
\sum_{i=1}^{k_n} p_{n,i} \to \lambda_1, \quad \sum_{i=1}^{k_n} q_{n,i} \to \lambda_2. \text{ Then} 
\end{equation}
a similar computation yields $T(P^{(n)}, Q^{(n)}) \to T(P(\lambda_1), P(\lambda_2))$ \cite{Durrett2019}[Theorem 3.6.1].
\end{corollary}

\begin{proposition}[More infinitely divisible limits]
Given trade-off functions  $\{f^{(n)}, g^{(n)}: n\in \mathbb{Z}_{\geq 1}\}$ \label{prop:moreidp}
\begin{equation}
\text{If } f^{(n)} \to f_{\infty} \in \cT \text{ and } g^{(n)} \to g_{\infty} \in \cT \text{ poinwise } \text{ then }f^{(n)} \otimes g^{(n)} \to f_{\infty}\otimes g_{\infty} \in \cT.
\end{equation}
One can specialize to the product case $f^{(n)}:=f_n^{\otimes n}, g^{(n)}= g_n^{\otimes n}$ for some $f_n, g_n \in \cT$. 
\end{proposition}
\begin{proof}
    First, recall that the pointwise convergence of monotone functions $\{f^{(n)}\}_n (\{g^{(n)}\}_n)$ mapping $[0,1] \to [0,1]$ to a limiting monotone continuous function $f_{\infty} (g_{\infty})$ mapping  $[0,1] \to [0,1]$ implies uniform convergence on $[0,1]$, and so we analyze under  $\|h_1-h_2\cdot \|_{\infty}^{[0,1]} = \sup\{h_1(x)-h_2(x): x \in [0,1]\}$. Now, the proof follows from a simple \textit{Lipschitz} argument, where we observe that
    \begin{equation}
        \|f^{(n)}\otimes g^{(n)}  - f_{\infty} \otimes g_{\infty}\|_{\infty}^{[0,1]}
        \leq 
        \|f^{(n)}- f_{\infty}\|_{\infty}^{[0,1]} 
        +
        \|g^{(n)}  - g_{\infty}\|_{\infty}^{[0,1]}
        \to 0
    \end{equation}
   This follows from  the `trick' $|G(a_n,b_n)-G(a,b)| = |G(a_n,b_n)- G(a_n,b) + G(a_n,b) - G(a,b)|$
    \begin{align}
        &
        \text{applied to }
        \|g_1\otimes f_1 - g_2 \otimes f_2\|_{\infty}^{[0,1]}
        =
        \|g_1\otimes f_1 - g_1\otimes f_2 + g_1\otimes f_2- g_2 \otimes f_2\|_{\infty}^{[0,1]}
        \\
        &
        \leq 
        \|g_1\otimes f_1 - g_1\otimes f_2\|_{\infty}^{[0,1]} 
        + 
        \|g_1\otimes f_2 - g_2\otimes f_2\|_{\infty}^{[0,1]} 
        \overset{(a)}{\leq}
        \|f_1 - f_2\|_{\infty}^{[0,1]}
        +
        \|g_1 - g_2\|_{\infty}^{[0,1]} 
    \end{align}
    
    where the inequality $(a)$ $\|g\otimes f_1 - g\otimes f_2\|_{\infty}^{[0,1]}  \leq  \|f_1-f_2\|_{\infty}^{[0,1]}$ is proven below.
\end{proof}
\begin{lemma}[Le-Cam distance under product] \label{lem:LeCdef}
    Consider trade off functions $g,f_1,f_2 \in \cT$. Then,
    \begin{equation}
        \|g\otimes f_1 - g\otimes f_2\|_{\infty}^{[0,1]}  \leq  \|f_1-f_2\|_{\infty}^{[0,1]}.
    \end{equation}
\end{lemma}
\begin{proof}
We need to show that for any given $\varepsilon >0$ if
for all $x\in [0,1], |f_1(x)-f_2(x)|\leq \varepsilon$, then 
\begin{equation}
-\varepsilon \leq g\otimes f_1 (\alpha) - g\otimes f_2(\alpha) \leq \varepsilon
\text{ for all } \alpha \in [0,1] 
\end{equation}
  Now, denote $g = T(P,Q),  f_1 = T(P_1,Q_1), f_2=T(P_2,Q_2)$, corresponding to binary Blackwell-Le-Cam experiments  $(\Omega, \cF, P, Q) \leftrightarrow g, (\Omega_1, \cF_1, P_1, Q_1) \leftrightarrow f_1,  (\Omega_2, \cF_2, P_2,Q_2) \leftrightarrow f_2\footnote{First, the exact representation of the probability measures such as $(P,Q)$  and the associated sample spaces  $(\Omega, \cF)$ on which they are defined do not appear in the conclusion of the result. Second, this proof relies on the fact that one can construct a jointly measurable map $f(w,w_2): \Omega \times \Omega_2 \to [0,1]$ from a collection $\{f^{w}(\cdot): w\in \Omega\}$ of marginally measurable maps $f^{w}(w_2) : \Omega_2 \to [0,1]$, we restrict ourselves to the canonical choice of the standard Borel space $ \Omega= \Omega_1= \Omega_2= [0,1]$ with $\cF=\cF_1=\cF_2= \cB_{[0,1]}$}.$

So, we need to show $\|T(P\otimes P_1, Q \otimes Q_1)- T(P\otimes P_2, Q\otimes Q_2)\| \leq  \|T(P_1, Q_1) - T(P_2,Q_2)\|$

Now, consider $\phi'_1: \Omega\times \Omega_1 \to [0,1]$ as the (Borel measurable) Neyman-Pearson optimizer at level $\alpha =P\otimes P_1(\phi'_1)$, with the value of the minimum type II error $g\otimes f_1(\alpha) = 1- Q\otimes Q_1(\phi'_1)$. 
It suffices to construct  a Borel measurable $\phi'_2: \Omega\times \Omega_2 \to [0,1]$  such that $\alpha = P\otimes P_2(\phi'_2)$ and $g\otimes f_2(\alpha) \leq   1-Q\otimes Q_2(\phi'_2) \leq g\otimes f_1(\alpha) + \varepsilon$. For a fixed $w \in \Omega $ define  $\phi^{'w}_1: \Omega_1 \to [0,1]$ be the slice of $\phi'_1$ at $w \in \Omega$ as $\phi^{'w}_1(w_1) = \phi'_1(w,w_1)$, and therefore Borel  measurable. Then $\phi^{'w}_1$ is a suboptimal test between $P_1$ and $Q_1$ with type I error $\alpha^w= P_1(\phi^{'w}_1)$ and type II error $Q_1(1- \phi^{'w}_1) \geq f_1(\alpha^w)$. 

Now, consider the optimal (Borel measurable) Neyman-Pearson test $\phi^{'w}_2: \Omega_2 \to [0,1] $  between $P_2$ and $Q_2$ at level $\alpha^w = P_2(\phi_2^{'w})$ so that $f_2(\alpha^w)= Q_2(1-\phi_2^{'w})$. Moreover, by patching these marginally measurable maps\footnote{Although this is outside the scope of this paper, we believe some abstract results of the kind required to define joint Radon-Nikodym derivatives out of marginal ones \cite{polyanskiy_wu_2025}[Theorem 2.12],  \cite{Doob1953StochasticProcesses}, \cite{Cinlar2011ProbabilityStochastics} would be helpful to conclude, in our context, that indeed such a jointly measurable map $\phi_2'$ can be constructed. given that all the intermediate choices we made in the proof are well behaved.}, we obtain a jointly measurable mapping $\phi'_2: \Omega \times \Omega_2 \to [0,1]$ defined as $\phi'_2(w,w_2) = \phi_2^{'w}(w_2)$.   Observe that the type I error of the map $\phi_2'$  is given by 
\begin{equation}
P\otimes P_2 (\phi'_2)= P(\alpha^w) =P \otimes P_1(\phi'_1) =\alpha.
\end{equation}
Now, consider the type II error 
$g\otimes f_2(\alpha) \leq Q\otimes Q_2 (1- \phi_2')= Q(f_2(\alpha^{w})) \leq Q(f_1(\alpha^w)) +\varepsilon) \leq Q\otimes Q_1 (1-\phi_1') +\varepsilon  = g\otimes f_1(\alpha) +\varepsilon$.
As a consequence we have $g\otimes f_2 (\alpha) \leq g\otimes f_1(\alpha) + \varepsilon$, and since $\alpha\in [0,1]$ was arbitrary, and moreover one can switch $f_1$ and $f_2$ in the conclusion above, since the assumption $\|f_1-f_2\|_{\infty} \leq \varepsilon$ is symmetric, and therefore the result holds immediately.
\end{proof}

\textbf{Blackwell Le-Cam perspective.}  The proof of Lemma \ref{lem:LeCdef} is inspired from  the proof of \cite{dong2022gdp}[Lemma C.1] which proves the comparison version:  for trade-off functions $g,f_1,f_2 \in \cT$ if $f_1\leq f_2$ then $g \otimes f_1 \leq g\otimes f_2$. However, we comment that this quantitative $\varepsilon$ generalization
$\|f_1- f_2\|_{\infty}^{[0,1]}\leq  \varepsilon$ implies $\|g\otimes f_1- g\otimes f_2\|_{\infty}^{[0,1]}\leq  \varepsilon$ closely matches with the generalization of Le-Cam \cite{LeCam1986} going from Blackwell-comparison ($0$-deficiency) to $\varepsilon$-deficiency. More precisely, a similar result holds for the Le-Cam distance as well $\Delta (\cE \otimes \cE_1, \cE\otimes \cE_2) \leq \Delta(\cE_1, \cE_2)$ for binary experiments $\cE, \cE_1, \cE_2$ \cite{Torgersen1991Comparison}, \cite{LeCam1986}. 

\textbf{A different proof of Lemma \ref{eq:Lecamdef}}  Observe that the proof of Lemma \ref{eq:Lecamdef}
 was conducted by examining various Neyman-Pearson optimizers. This is because no explicit description of $g\otimes f$ was given purely in terms of trade off functions $g$ and $f$.  Now, by a reformulation of \cite{Torgersen1991Comparison}[Complement 2, Chapter 10] we observe that with $Q_g'[0, \alpha] = 1- g(\alpha)$ for $\alpha \in (0,1]$ and $Q_f'[\{0\}] = 1-g(0)$
 \begin{equation}
     g\otimes f (\alpha)=  \inf_{\varphi} \{ Q_g'(f(\varphi(\cdot))): ([0,1], \cB_{[0,1]}) \overset{\varphi}{\to} ([0,1], \cB_{[0,1]}) \text{ decreasing }, \int_{0}^{1} \varphi(x) =\alpha\}
 \end{equation}
 Now, using the variational representation above, the proof of the Lemma \ref{eq:Lecamdef} is immediate, since $\|f_1-f_2\|_{\infty} \leq \varepsilon$ implies $Q_g'(f_2(\varphi(\cdot))) \leq Q_g'(f_1(\varphi(\cdot))) +\varepsilon$ for any given (appropriate) $\varphi$, and therefore $g\otimes f_1(\alpha) \leq  g \otimes f_2(\alpha)  +\varepsilon$. Now, $\alpha \in [0,1]$ is arbitrary, and one can switch $f_1$ and $f_2$  as the assumption $\|f_1-f_2\|_{\infty} \leq \varepsilon$  is symmetric. Hence, the result of the lemma \ref{eq:Lecamdef} holds. 

\textbf{More infinitely divisible limits.} A consequence of Proposition \ref{prop:moreidp}, the closure of the infinitely divisible trade-off functions $\cI_{\cT}$ under tensor products \ref{lem:TOF1} and weak convergence\footnote{Weak limits of infinitely divisible distributions are infinitely divisible \cite{Kallenberg2021FMP3}[Theorem 7.7].}, and given that one can explicitly derive Gaussian \ref{cor:Gausslimit} and Poisson limits \ref{cor:Poissonlimit}, one can allow $P_{\infty} \circ \log\left(\frac{dQ_{\infty}}{dP_{\infty}}\right)^{-1}= G+P$ to have a Gaussian $G$ and a Poisson $P$ component. Theoretically, given any infinitely divisible random variable $X$ on $\mathbb{R}$ with $\e[\exp (\mathbb{X})]=1$, one can construct $P_{\infty} \overset{d}{=} X$ and define $dQ_\infty (x)= e^{x}dP_{\infty}(x)$. Then, since $X$ is infinitely divisible, by lemma \ref{lem:TOF1} there exists $f_n=T(P_n, Q_n) \in \cI_{\cT}$ with $dQ_n =e^x P_n$ such that $T(P_n^{\otimes n}, Q_n^{\otimes n}) = T(P_{\infty}, Q_{\infty})$. Now, as $n\to  \infty$,  $f_n^{\otimes n} \to f_{\infty} = T(P_{\infty}, Q_{\infty})$.

 \textbf{Beyond infinite divisibility (mixtures)}
Here, we demonstrate how to go beyond the infinitely divisible privacy framework $\cI_{\cT}$. We achieve this by working under a random (but growing) number of differentially private operations\footnote{The intuition behind this is that in practice, often one does not know based on the same dataset $S$, how many differentially private statistics $(M_1(S), M_2(S, M_1(S)), \cdots)$  one has to release over time. So, we capture it as a time change sequence of integer-valued random variables $\{N_n:n\in \mathbb{Z}\}$. }.  Consider the canonical probability space $([0,1]^{\infty}, \cB_{\infty}, U^{\infty})$ as the countable product of the standard uniform distribution $U$ on $[0,1]$. Now consider a sequence of trade-off functions $\{f_n: n \in \mathbb{Z}_{\geq 1}\} \subset \cT$ with  $f_n(0)=1$ so that $f_n=T(U,Q_{f_n})$ with $Q_{n}:=Q_{f_n} \ll U$.  Define the likelihood ratio process  $\{g_n(t):=L_n(\lfloor t n \rfloor): t\geq 0\}$ under $U^{\infty}$ 
\begin{equation}
L_n(k)(\ul{x}) 
:=
\log\left(\frac{dQ_n^{\otimes k}}{dU^{\otimes k}}\right) (\ul{x})
= 
\sum_{i=1}^{k} \log |f_n'(1-x_i)| \text{ for } \ul{x} \in  [0,1]^{\infty} 
.\end{equation}  
Assume that there exists a Lévy process  $\{L(t): t\ge 0\}$ so that the entire  likelihood ratio process $\{g_n (t) :=L_n(\lfloor t n \rfloor): t\geq  0\}$, under $U^{\infty}$  converges in distribution  to the process $\{L(t): t\geq 0\}$ 
\begin{equation}
\text{ under }
U^{\infty}
\quad
\{g_n (t): t\geq  0\} \overset{d}{\to} 
\{L(t): t\geq 0\},
\text{ and }
\mathbb{E}[ \exp (L(t))] =1  \text{ for all } t\geq 0.
\end{equation}
Consider a sequence  $\{N_n: n\in \mathbb{Z}_{\geq 1}\}$ be integer-valued random variables on some probability space $(\Omega',\cF', \mathbb{P}')$, and consider the joint product space  $(\Omega' \times [0,1]^{\infty},\cF \otimes \cB_{\infty}, \mathbb{P}' \otimes U^{\infty})$, so that\footnote{This independence assumption precisely captures the intuition that we do not know in practice (and hence model it as random) how many rounds of differentially private operations one has to run on a given dataset $S$.} 
\begin{equation}
\text{ for all } 
n \in \mathbb{Z}_{\geq 1}
\text{ we have independence }
N_n \perp  g_n.
\end{equation}
We assume $\mathbb{P}'$ we have  $N_n/n \overset{d}{\to}\tau$ for some probability measure  $\tau$  on $([0,\infty), \cB_{[0,0\infty)}$\footnote{Under independence assumption $\{X_n\}_n \perp \{Y_n\}_n$ marginal convergence $X_n \overset{d}{\to} X$, $Y_n \overset{d}{\to} Y$ implies joint convergence $(X_n, Y_n) \to (X,Y)$ \cite{vandervaart1998asymptotic}.  We apply this to $(g_n, N_n/n) \overset{d}{\to} (g,\tau)$ }. Then
\begin{equation}
\text{under } \mathbb{P}'\otimes U^{\infty} \quad g_n\left(\frac{N_n}{n}\right)
=
L_n(N_n) \overset{d}{\to} L(\tau), \text{ and } \mathbb{E}[\exp L(\tau)] =1 
\end{equation}
Moreover, by \ref{lem:pollard} we have pointwise convergence of  $f_n^{\otimes N_n} = T(U^{\otimes N_n}, Q_{n}^{N_n}) \to f_{\infty}$ on $[0,1]$ where
\begin{equation} \label{eq:Levylimit}
    f_{\infty} = T(P_{\infty}, Q_{\infty}), \text{ with }
    P_{\infty} \overset{d}{=} \Law(L(\tau)) \text{ and }
    dQ_{\infty}(x) = e^{x} dP_{\infty}(x).
\end{equation}

\textbf{Basics of Levy processes.} A Levy process $L$ with $L(0)\equiv 0$ is a stochastic process $\{L(t): t\geq 0\}$, a collection of random variables on some probability space  (assumed to be right-continuous function with left limits almost surely) with independent\footnote{It means $(L(t_1), L(t_2)- L(t_1), \cdots, L(t_k)- L(t_{k-1}))$ are jointly independent for all $0\leq t_1 < \cdots < t_k$.} and stationary\footnote{Stationarity $L(t)-L(s) \overset{d}{=} L(t-s)$ is dropped, when one analyses triangular systems $f_{n1} \otimes \cdots \otimes f_{nN_n}$.} increments. 

Since we are interested only in weak convergence\footnote{We mean weak convergence of the entire process, which implies convergence of finite dimensional distributions and much more. We refer to \cite{JohnstoneOnatski2020SpikedModels}[Conclusions] for more details on this.}, we recall that marginal distributions of a Levy process $\{\mu_t: t\geq 0\}$ are all infinitely divisible satisfying $\mu_{t+s} = \mu_t *\mu_s$ for all $s,t \geq 0$ with $\mu_0= \delta_0$. Moreover, the joint distribution of the process $\{L(t): t\ge 0\}$ is determined by the $\mu_1 \overset{d}{=} L(1)$ \cite{Kallenberg2021FMP3}[Lemma 16.1] and can be written explicitly as a sum of drifted Brownian motion 
\begin{equation}
\{mt + \sigma B_t:t \geq 0\} \text{ for a standard Brownian motion } \{B_t:t \geq 0\}, \text{ and}
\end{equation} 
a Poisson process (see \cite{Kallenberg2021FMP3}[Theorem 16.2] for details on the Poisson process part).

\textbf{Gaussian mixtures example.} For demonstrations of our results of \eqref{eq:Levylimit} going beyond the infinitely divisible privacy framework of $\cI_{\cT}$, we restrict ourselves to limiting Levy processes $\{L(t): t\geq 0\}$ with almost surely continuous paths, and hence reduce to the case of drifted Brownian motion \cite{Kallenberg2021FMP3}[Proposition 16.9] $\{L(t) = mt + \sigma B_t: t\geq 0\}$ for some real number  $m= -\frac{\sigma^2}{2}$\footnote{This constraint comes from the restriction $\e[\exp L(t)]  =\e[\exp (mt +\sigma B_t)] =1$ for all $t \geq 0$.} 

For simplicity, consider the case where $\tau = \lambda_1 \delta_{t_1} + \cdots \lambda_k \delta_{t_k}$ for some $ 0 <t_1 < \cdots <t_k < \infty$ with $\lambda_i \geq 0$ and $\sum_{i=1^k} \lambda_i =1$. Then with the notation $\mu_{t} = N( -\frac{\sigma^2}{2} t, \sigma^2 t)$ for $t\geq 0$ we have

\begin{equation} \label{eq:GMM}
P_{\infty}= \lambda_1  \mu_{t_1} + \cdots  +\lambda_k \mu_{t_k},
\text{ and }
dQ_{\infty}=   \lambda_1  e^xd\mu_{t_1} + \cdots + \lambda_k e^x d\mu_{t_k},
\text{ so }
f_{\infty}= T(P_{\infty}, Q_{\infty}).
\end{equation}

It is known that nontrivial mixtures of Gaussians are not infinitely divisible \cite{BoseDasGuptaRubin2002}. So, such a trade-off function $f_{\infty}= T(P_{\infty}, Q_{\infty})$ is not just a departure from the asymptotic Gaussian differential privacy  \cite{dong2022gdp}, but also from our enlarged infinitely divisible privacy $\cI_{\cT}$, which we established under a large (but deterministic) number of differentially private operations.

\begin{figure}[h]
    \centering
    \includegraphics[width=0.7\linewidth]{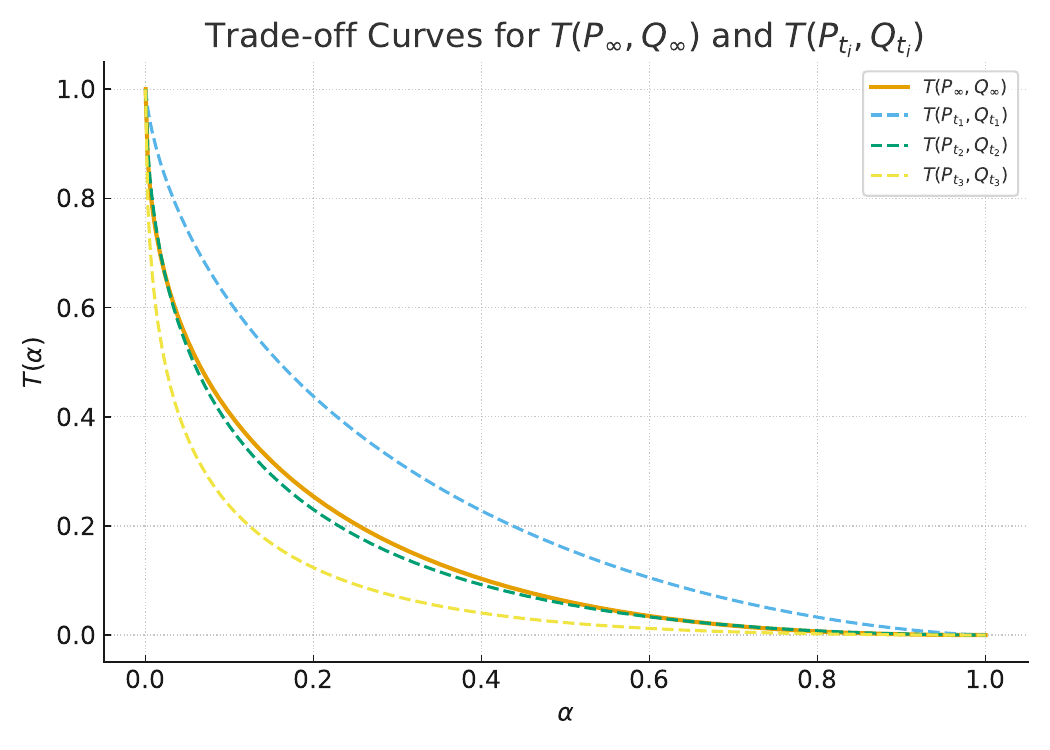}
    \caption{
    Trade-off curves for the Gaussian mixture experiment
    $T(P_{\infty}, Q_{\infty})$ and for the component Gaussian experiments
    $T(P_{t_i}, Q_{t_i})$, $i=1,2,3$.
    Here
    $P_{\infty} = \sum_{i=1}^3 \lambda_i \mu_{t_i}$ and
    $dQ_{\infty}(x) = \sum_{i=1}^3 \lambda_i e^x \, d\mu_{t_i}(x)$,
    with $\mu_{t_i} = \mathcal N(-\tfrac{\sigma^2}{2} t_i, \sigma^2 t_i)$.
    The solid curve shows $T(P_{\infty}, Q_{\infty})$, while the dashed curves
    show the individual components $T(P_{t_1}, Q_{t_1})$,
    $T(P_{t_2}, Q_{t_2})$, and $T(P_{t_3}, Q_{t_3})$.
    }
    \label{fig:mixture-vs-components-tradeoff}
\end{figure}

\textbf{Blackwell-Le-Cam perspective.} The  Gaussian mixture example $(P_{\infty}, Q_{\infty})$ with the mixing distribution  $\tau = \lambda_1 \delta_{t_1} + \cdots \lambda_k \delta_{t_k}$ can be seen as a \textit{mixture} of binary experiments $\cE_i = (P_{t_i},Q_{t_i})$ where $P_{t_i} = \mu_{t_i} = N(-\frac{\sigma^2}{2}t_i, \frac{\sigma^2}{2}t_i)$ and $dQ_{t_i} (x)= e^{x} dP_{t_i} (x)$, all of which are on  $\Omega_i=\mathbb{R}, \cF_i =\cB_{\mathbb{R}}$.  $P_{\infty}$ can also be seen to fit within the framework of locally asymptotic mixed normality (`Normal' distributions with random variances), as introduced in \cite{Jeganathan1982AsymptoticEstimation}, and has recently also appeared in the asymptotics of (local) differentially private statistical analysis \cite{Steinberger2024EfficiencyLDP}.

Moreover, a reformulation of \cite{Torgersen1991Comparison}[Complement 2, Chapter 10] is that for a binary  mixture experiment $(P_{\infty}, Q_{\infty})$ by mixing  experiments $\{(P_i,Q_i) :1 \leq i \leq k\}$  so that $P_{\infty} = \sum_{i=1}^k \lambda_i P_{i}$, and $Q_{\infty}= \sum_{i=1}^k \lambda_i Q_{i}$. Then we have for $f_{\infty}= T(P_{\infty}, Q_{\infty})$ and $f_{i}= T(P_i, Q_i)$
\begin{equation} \label{eq:tofmixture}
    f_{\infty}(\alpha) 
    =
    \inf_{0 \leq \alpha_i \leq 1} \{\sum_{i=1}^k \lambda_i f_i (\alpha_i) : \sum_{i=1}^k \alpha_i = \alpha\}
    \text{ for all } \alpha \in [0,1]
\end{equation}

\textbf{Wald-Wolfowitz extension.} In this `mixture' framework above, we assumed that the random (but growing) number $\{N_{n}\}_{n \geq 1}$ of differentially private operations are independent of the `outcomes' of the differentially private operations $M^{(n,k)}(\cdot) = (M_{n1}(\cdot),  \cdots, M_{nk}(\cdot))$ on a given dataset $S$. 

Under the sequential analysis, one might consider stopping at time $N_n$ based on the release of the previous outputs $M^{n, N_n-1}(S)$,  if enough information about the data has already been released by that time. Although this is beyond the scope of this paper, we think that, from a practical point of view, it is an extremely interesting direction to pursue, where we believe a different analysis, more in line with Wald-Wolfowitz's framework of the sequential probability ratio test \cite{Siegmund1985SequentialAnalysis}, along with the techniques of  \cite{PenaGine1999Decoupling}, 
\cite{deLaPenaLaiShao2009SelfNormalizedProcesses}, 
\cite{Silvestrov2004RandomlyStopped}.

We believe that the result \ref{eq:Levylimit} would still hold, where $\tau$ would be a stopping time  with respect to the Filtration $\{\cF_t= \sigma(L(s): 0\leq s\leq t):t\geq 0\}$ generated by the Levy process $\{L(t):t\geq 0\}$
\footnote{To be mathematically precise, in concluding, $\e[\exp(L(\tau))] =1$ we need some uniform integrability conditions on the sequence $\{\exp(L_n(N_n)): n \in \mathbb{Z}_{\geq 1}\}$ to apply Martingale optional stopping theorems.}.
\begin{figure}[h]
    \centering
    \includegraphics[width=0.7\linewidth]{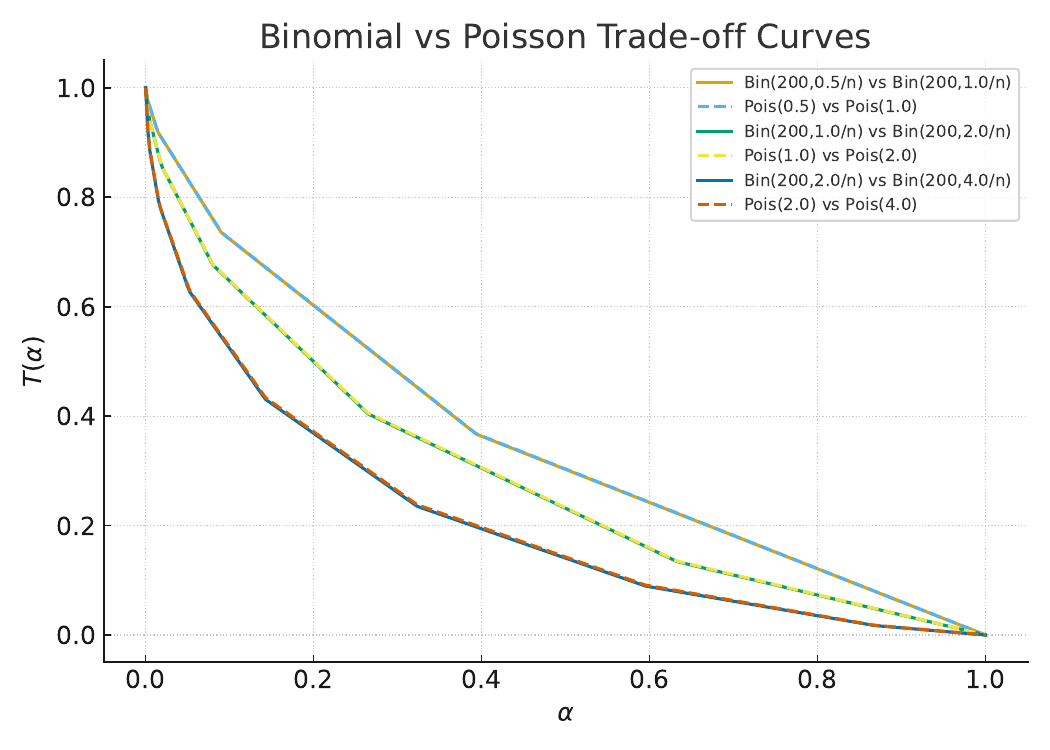}
    \caption{
    Trade-off curves for testing between $\mathrm{Bin}(n,\lambda_1/n)$ and
    $\mathrm{Bin}(n,\lambda_2/n)$ (solid lines, with $n=200$) compared to
    the Poisson limits $T(\mathrm{Pois}(\lambda_1),\mathrm{Pois}(\lambda_2))$
    (dashed lines), for several pairs $(\lambda_1,\lambda_2)$.
    The curves are nearly indistinguishable, illustrating the Poisson approximation at the level.
    }
    \label{fig:binom-poisson-tradeoff}
\end{figure}

\section{How to be Poisson differentially private in practice} \label{sec:BIDP}

Now, we describe a mechanism for achieving a baseline notion of Poisson differential privacy $f_{\mu}= T(P(1), P(\mu))$ for Poissons  $P(\cdot)$ with mean $\mu>1$ \ref{def:fDP}.  The purpose is similar to the noise-adding mechanism results for achieving Gaussian differential privacy $f_{\mu}= T(Z, Z+\mu)$ for $Z\sim N(0,1)$ or, more generally, $f_{\mu}= T(X, X+\mu)$ for a symmetric log-concave random variable $X$ \ref{lem:SLCTOF}. However, as we will demonstrate, for baseline privacy families such as $\{P(\lambda): \lambda >1\}$ with respect to $P(1)$, that are not inherently a shift family including $\{X+ \mu: \mu >0\}$ with respect to $X$ (distribution of $X$ is symmetric and log-concave), noise adding mechanisms is not the natural way to achieve privacy.

Consider a metric space $(\cX,d)$\footnote{The definition of a metric space includes the product space example $\cX= \cX_1 \times \cdots \times \cX_n$ under the Hamming metric $d_H(\ul{x}, \ul{y}) = \sum_{i=1}^n \mathbf{1}(x_i \neq y_i)$, capturing how many coordinates $\ul{x}$ and $\ul{y}$ differ.} It consists of the allowable datasets $S$ and a real-valued statistic $g: (\cX, d) \to \mathbb{R}$.  Denote the sensitivity or the modulus of continuity at 1 as\footnote{The modulus of continuity $w_{g}(\delta)$ is an important  measure of continuity of $g$. In the literature on statistical estimation \cite{DonohoLiu1991GeometrizingIII} under differential privacy constraints \cite{RohdeSteinberger2020GeometrizingLDP}, \cite{DuchiRuan2024RightComplexityLDP}, one usually considers the metric space $\cX$ as an appropriate subset of Borel probability measures $\cP(\mathbb{R}^d, \cB_{\mathbb{R}^d})$ with an appropriate choice of the metric $d$, including total variation and Hellinger distance.} 
\begin{equation}
    w_g(1):= \sup_{d(x,y)\leq 1} |g(x) -g(y)| 
    \quad 
    \text{ satisfying }
    \quad 
    w_{cg}(1) = c w_{g}(1) \text{ for any } c>0.
\end{equation} 
\begin{lemma}[Modulus of continuity] \label{lem:MOC}
For a map $(X,d_X) \overset{g}{\to} (Y,d_Y)$ between metric spaces, define
\begin{equation}
w_g: [0, \infty] \to [0, \infty]
\quad 
w_g(\delta)= \sup_{d_X(x_1,x_2)\leq \delta} d_Y(g(x_1),g(x_2)).
\end{equation}
Consider the composition $(X,d_X) \overset{g}{\to} (Y,d_Y) \overset{f}{\to} (Z,d_Z)$ then $w_{f\circ g} (\delta) \leq w_f(w_g(\delta))$ for all $\delta \geq 0$.
\end{lemma}
\begin{proof}
Define $L=w_{f \circ g} (\delta)= \sup_{d_X(x_1,x_2)\leq \delta } d_Z(f(g(x_1)),f(g(x_2)))$. Now  for all $x_1,x_2\in X$
\begin{equation}
d_X(x_1,x_2)\leq \delta \rightarrow d_Y(g(x_1),g(x_2)) \leq w_g(\delta).
\end{equation}
Consequently, $\{(g(x_1),g(x_2)):  d_X(x_1,x_2)\leq \delta\} \subset \{(y_1, y_2):  d_Y(y_1,y_2)\leq w_g(\delta)\}$. Therefore 
\begin{equation}
L
\leq  
\sup_{d_Y(y_1,y_2)\leq w_g(\delta)}  d_Z(f(y_1),f(y_2))
\leq 
w_f(w_g(\delta))
\end{equation}
\end{proof}
\begin{theorem}\cite{dong2022gdp}[Theorem 1, Proposition A.3] \label{thm:GDPdong} Given a metric space $(\cX,d)$ of datasets and a real-valued statistic $g: (\cX,d) \to \R$, consider the noise-adding mechanism (Markov Kernel)   $M:(\mathbb{R}, \cB_{\mathbb{R}}) \to (\mathbb{R}, \cB_{\mathbb{R}}) $  as $M(g(x)) = g(x) + X$ for a symmetric-log-concave  $X$ \footnote{With an abuse of notation, we denote the mechanism $M$ and the corresponding Markov kernel $M(z,A) = \mathbb{P}(z+ X \in A)= P_X(A-z)$ for $A \subset \cB_{\mathbb{R}}$ and $z \in \mathbb{R}$, where $A-z$ is the Minkowski translate $\{A-z\}:= \{a-z: a \in A\}$, and $P_X$ is a symmetric ($P_X(A)= P_X(-A), \text{ with } -A:= \{-a: a\in A\})$, log-concave ($x\ \to \log\frac{dP_X}{d\lambda}(x)$ is a concave function on the convex support $(\text{Range}(X))$) measure, the distribution of $X$.}. Then 
\begin{equation} \label{eq:GDPdong}
  \inf_{d(x,y)\leq 1} T(M(g(x)), M(g(y))) \geq  T(X, X + w_g(1))
\end{equation}
Consequently, the mechanism $M$ satisfies $f^X_{w_{g}(1)}$-DP with $f^X_{\mu}= T(X, X+\mu)$ for any $\mu >0$.
\end{theorem}
\textbf{Noise level.} Observe that having fixed the noise distribution $X$, the noise level $\mu>0$ and the statistic  $g$ before implementing the mechanism $M$, one can scale $g \to cg$  so that $w_{cg}(1)=c w_g(1)= \mu$ (if the mechanism designer is allowed to scale the original statistic $g$) or equivalently scale the noise $X \to cX$ with $c$ so that $T(cX, cX + w_g(1))= T(X, X + c^{-1} w_g(1))$ with $c^{-1} w_g(1)= \mu$. 

Now, we describe a mechanism that achieves an asymmetric version of Poisson differential privacy.
\begin{theorem}[Poisson differential privacy] \label{thm:PDP}
Given a metric space $(\cX,d)$ of datasets and a count  statistic $g: (\cX,d) \to \mathbb{Z}_{\geq 0}$\footnote{We are trying to achieve Poisson differential privacy, which is suitable when we expect the released outputs to be count statistics (elements of $\mathbb{Z}_{\geq 0}$). So, its natural to consider the original statistic $g$ to also take values in $\mathbb{Z}_{\geq 0}$. However, a quick inspection of the proof reveals that our results hold for any real valued statistic $g$, since nowhere in the proof do we use the fact that $g$ takes only integer values or non-negative values.} with the baseline Poisson trade-off function $f_{\infty}= T(P(\mu_1), P(\mu_2))$ for some $\mu_2 >\mu_1 >0$. Consider the mechanism $M: (\mathbb{Z}_{\geq 0}, 2^{\mathbb{Z}_{\geq 0}}) \to (\mathbb{Z}_{\geq 0}, 2^{\mathbb{Z}_{\geq 0}})$ defined as

\begin{equation} \label{eq:NN}
   M(g(x)) \sim  P(N_2 e^{N_1 g(x)}), \text{ with } N_1=  \frac{\log  \left(\frac{\mu_2}{\mu_1}\right)}{w_g(1)} \text{ and }
   N_2= \frac{|\mu_2 -\mu_1|}{w_{h\circ g}(1))}, \text{ where }
\end{equation}
\begin{equation} \label{eq:h}
    h: (\mathbb{R}, |\cdot|)  \to (\mathbb{R}, |\cdot|) \quad 
    h(y) = e^{N_1 y}= \left(\frac{\mu_2}{\mu_1}\right)^{\frac{y}{w_g(1)}}.
\end{equation}
Then $M$ satisfies the following (asymmetric) version of Poisson differential privacy\footnote{For  $f= T(P,Q)$ implies the generalized inverse $f^{-1} = T(Q,P)$ \cite{dong2022gdp}[Lemma A.2].}.
\begin{equation} \label{eq:part1}
    \inf_{d_{\cX}(x_1,x_2) \leq 1: g(x_2) \geq  g(x_1)} T(M(g(x_1)), M(g(x_2))) \geq  T(P(\mu_1), P(\mu_2)) = f_{\infty}
\end{equation}
\begin{equation} \label{eq:part2}
    \inf_{d_{\cX}(x_1,x_2) \leq 1: g(x_2) \leq  g(x_1)} T(M(g(x_1)), M(g(x_2))) \geq  T(P(\mu_2), P(\mu_1)) 
    = f_{\infty}^{-1}
\end{equation}
Consequently, $M$ satisfy (symmetric) $f_{\mu_1, \mu_2}$-DP \ref{def:fDP} where  $f_{\mu_1, \mu_2}= \min(f_{\infty}, f_{\infty}^{-1})$. More precisely,

\begin{equation} \label{eq:final}
     \inf_{d_{\cX}(x_1,x_2) \leq 1} T(M(g(x_1)), M(g(x_2))) \geq  f_{\mu_1, \mu_2}
\end{equation}

\end{theorem}

\begin{proof}
For parameters $\lambda_1', \lambda_2', \mu_1', \mu_2' >0$,  if  $(\lambda_1',\lambda_2')$ and $(\mu_1',\mu_2')$ are similarly ordered (this means either $\lambda_1' < \lambda_2'$ and $\mu_1' <\mu_2'$ simultaneously  or $\lambda_2' <\lambda_1'$ and $\mu_2' <\mu_1'$  simultaneously)  and satisfy
\begin{equation} \label{eq:suffcon}
  \textbf{sufficient conditions:}
  \quad 
  |\lambda_2'- \lambda_1'| \leq |\mu_2'- \mu_1'| 
    \text{ and  }
    \frac{\max(\lambda_1', \lambda_2')}{\min(\lambda_1', \lambda_2')} \leq 
    \frac{\max(\mu_1', \mu_2')}{\min(\mu_1', \mu_2')} 
\end{equation}

   Then by taking $0< c:= \frac{|\lambda_2'- \lambda_1'|}{ |\mu_2'- \mu_1'| } \leq 1$  and a $\lambda \geq 0$ defined in the following way\footnote{$(c,d)$ captures the (slope, intercept) pair of the line $y= cx + d$ joining the points $(\mu_1', \lambda_1')$ and $(\mu_2', \lambda_2')$.} 

    \[
\lambda 
:= \lambda_1' - c\,\mu_1'
=
\begin{cases}
\displaystyle 
\frac{\lambda_1' \mu_2' - \lambda_2' \mu_1'}{\mu_2' - \mu_1'} 
& \text{when } \lambda_2' > \lambda_1' \text{ and } \mu_2' > \mu_1' \text{ satisfying } \eqref{eq:suffcon}, \\[1.2em]
\displaystyle
\frac{\lambda_2' \mu_1'- \lambda_1' \mu_2'}{\mu_1' - \mu_2'} 
& \text{when } \lambda_1' > \lambda_2' \text{ and } \mu_1' > \mu_2' \text{ satisfying } \eqref{eq:suffcon}.
\end{cases}
\]

 we have  $\lambda_1'= c \mu_1' + \lambda$ (by definition of $\lambda$), and  $\lambda_2' = c \mu_2' + \lambda$ holds because of the following.  

    \[
\lambda_2' - c\,\mu_2'
=
\begin{cases}
\displaystyle 
\frac{\lambda_1' \mu_2' - \lambda_2' \mu_1'}{\mu_2' - \mu_1'}  =\lambda 
& \text{when } \lambda_2' > \lambda_1' \text{ and } \mu_2' > \mu_1' \text{ satisfying } \eqref{eq:suffcon}, \\[1.2em]
\displaystyle
\frac{\lambda_2' \mu_1'- \lambda_1' \mu_2'}{\mu_1' - \mu_2'}  = \lambda
& \text{when } \lambda_1' > \lambda_2' \text{ and } \mu_1' > \mu_2' \text{ satisfying } \eqref{eq:suffcon}.
\end{cases}
\]

Then applying \eqref{eq:PS2} (Thinning) and  \eqref{eq:PT1} (Superposition) of Lemma \ref{lem:PTOF} we have the following chain of inequalities for any similarly ordered (all of them $>0$) $(\lambda_1',\lambda_2')$ and $(\mu_1',\mu_2')$ satisfying \eqref{eq:suffcon}
\begin{equation} \label{eq:keyineq}
    T(P(\mu_1'), P(\mu_2')) \overset{\eqref{eq:PS2}}{\leq} T(P(c\mu_1'), P(c\mu_2')) \overset{\eqref{eq:PT1}}{\leq} T(P(c\mu_1' + \lambda), P(c\mu_2' + \lambda)) = T(P(\lambda_1'), P(\lambda_2')).
\end{equation}
Now, it remains to show  with the choice of $N_1,N_2$ as mentioned in the statement of Theorem \ref{thm:PDP}
\begin{equation}
N_1=  \frac{\log  \left(\frac{\mu_2}{\mu_1}\right)}{w_g(1)} \text{ and }
   N_2= \frac{|\mu_2 -\mu_1|}{w_{h\circ g}(1))} \text{ with }  h(y) = e^{N_1 y}= \left(\frac{\mu_2}{\mu_1}\right)^{\frac{y}{w_g(1)}} \text{ for a given statistic}
\end{equation}
 $g: (\cX, d) \to \mathbb{Z}_{\geq 0}$ and a baseline trade-off function $f_{\infty}= T(P(\mu_1), P(\mu_2))$ for some $\mu_2 > \mu_1 >0$ 
  that   \eqref{eq:suffcon} holds for any $(\lambda_1', \lambda_2') = (N_2e^{N_1 g(x_1)}, N_2 e^{N_1 g(x_2)})$ against $(\mu_1', \mu_2')= (\mu_1, \mu_2)$ as 

 \begin{equation}
 \sup_{d_{\cX}(x_1,x_2) \leq 1: g(x_2) \geq  g(x_1) } \frac{\exp(N_1 g(x_2))}{\exp(N_1 g(x_1))} \leq 
    \frac{\mu_2}{\mu_1} = \frac{\max(\mu_1, \mu_2)}{\min(\mu_1, \mu_2)}
\leftrightarrow N_1  w_g(1) \leq  \log \left(\frac{\mu_2}{ \mu_1}\right) \text{ and }
    \end{equation}
    \begin{equation}
   \sup_{d_{\cX}(x_1,x_2) \leq 1: g(x_2) \geq  g(x_1) } |N_2 e^{N_1 g(x_1)} - N_2 e^{N_1 g(x_2)} | \leq |\mu_2- \mu_1| \leftrightarrow
   N_2 w_{h\circ g}(1)) \leq |\mu_2- \mu_1| 
\end{equation}

    Now, having verified the sufficient conditions for every pair $(x_1, x_2)$ such that $d_{\cX}(x_1, x_2) \leq  1$ and $g(x_2) \geq g(x_1)$, we have the following by applying the  conclusion of  \eqref{eq:keyineq}
\begin{equation} \label{eq:conc1}
\inf_{d_{\cX}(x_1,x_2) \leq 1: g(x_2) \geq  g(x_1)} T\left(P(N_2 e^{N_1 g(x_1)}), P(N_2e^{N_1g(x_2)})\right)  \geq T(P(\mu_1), P(\mu_2))=f_{\infty}.
\end{equation}

This concludes the proof of \eqref{eq:part1}, and by symmetry\footnote{The required sufficient conditions \eqref{eq:suffcon} to establish for the choice of $N_1,N_2$ as well as its conclusions \eqref{eq:keyineq} are symmetric with respect to order reversing map $(\lambda_1', \lambda_2'), (\mu_1', \mu_2') \to (\lambda_2', \lambda_1'), (\mu_2', \mu_1') $} we have  \eqref{eq:part2} (same choice of $N_1, N_2$ as above)

\begin{equation} \label{eq:conc2}
\inf_{d_{\cX}(x_1,x_2) \leq 1: g(x_2) \leq  g(x_1)} T\left(P(N_2 e^{N_1 g(x_1)}), P(N_2e^{N_1g(x_2)})\right) \geq T(P(\mu_2), P(\mu_1)) =f_{\infty}^{-1}.
\end{equation}
Now, combining the two above \eqref{eq:conc1}, \eqref{eq:conc2} we have the following \eqref{eq:final}.
\begin{equation}
    \inf_{d_{\cX}(x_1,x_2) \leq 1} T(M(g(x_1)), M(g(x_2))) \geq  \min (f_{\infty}, f_{\infty}^{-1})=:f_{\mu_1, \mu_2}.
\end{equation}
\end{proof}

Now, we remark one several aspects of this mechanism $M$ along with extensions of the result.

\begin{figure}[h]
    \centering
    \includegraphics[width=0.5\linewidth]{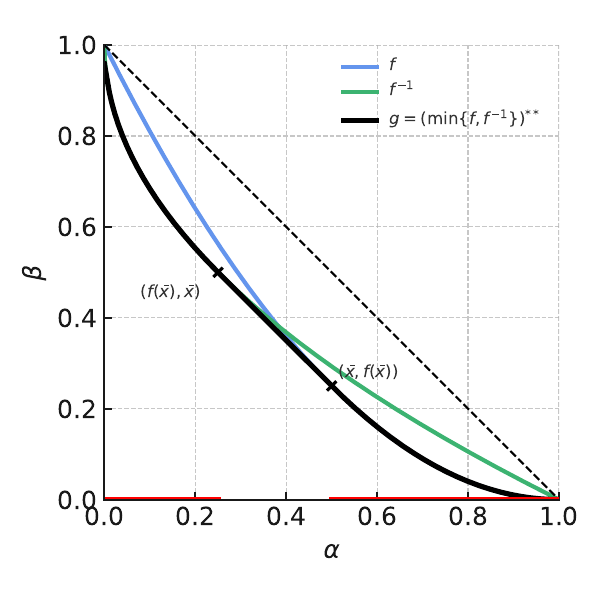}
    \caption{
    Double convex conjugate.  Here, the base trade-off function is
    \(f:[0,1] \to [0,1]\) given by $ f(x) = 1 - 2x + x^2,$ for $ x \in [0,1]$ and 
    and its functional inverse on $[0,1]$ is $f^{-1}(y) = 1 - \sqrt{y},$ for $ y \in [0,1].$  The black curve shows
    $g = \bigl(\min\{f,f^{-1}\}\bigr)^{**}$,
    the double convex conjugate of $\min\{f,f^{-1}\}$.
    The two marked points are
    $(f(\bar{x}), \bar{x})$ and $(\bar{x}, f(\bar{x}))$,
    with $\bar{x} = 1/2$ and $f(\bar{x}) = 1/4$.}
    \label{fig:symm-f-example}
\end{figure}

\begin{figure}[h]
    \centering
    \includegraphics[width=0.6\linewidth]{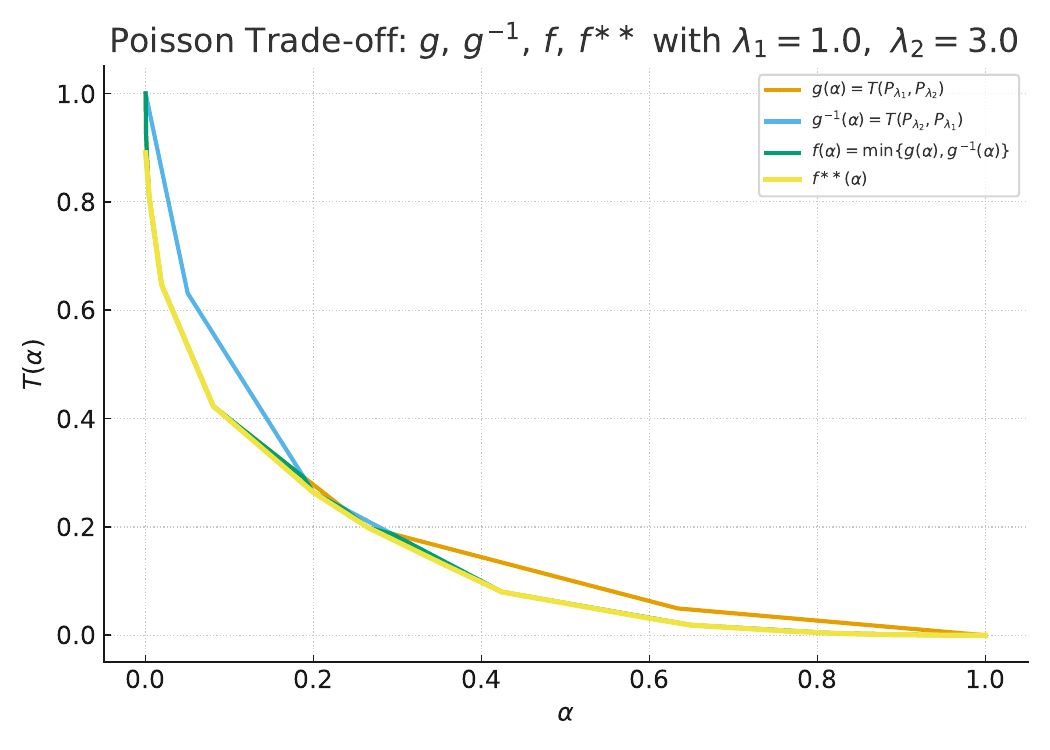}
    \caption{
    Trade-off functions for a Poisson experiment.
    Let $P_{\lambda_1} = \mathrm{Pois}(\lambda_1)$ and
    $P_{\lambda_2} = \mathrm{Pois}(\lambda_2)$ with
    $\lambda_1 = 1$ and $\lambda_2 = 3$.
    The blue curve shows
    $g(\alpha) = T(P_{\lambda_1}, P_{\lambda_2})(\alpha),$
    and the orange curve shows $g^{-1}(\alpha) = T(P_{\lambda_2}, P_{\lambda_1})(\alpha).$
    The green curve corresponds to  $ f(\alpha) = \min\{g(\alpha), g^{-1}(\alpha)\},$
    and the black curve shows the double convex conjugate $
        f^{**}(\alpha) = \bigl(\min\{g, g^{-1}\}\bigr)^{**}(\alpha).$
    }
    \label{fig:poisson-tradeoff-g-ginv-f-fss}
\end{figure}

\begin{remark}[Double conjugate] $f_{\mu_1, \mu_2} = \min(T(P(\mu_1), P(\mu_2)), T(P(\mu_2), P(\mu_1)))$ for $\mu_1, \mu_2 >0$ is not a trade-off function, since the minimum of two convex functions is not a convex function in general, and it does not hold  in the Poisson trade-off functions (see figure \ref{fig:poisson-tradeoff-g-ginv-f-fss} for a counterexample).

Therefore,  to fit with the $f$-DP framework of \cite{dong2022gdp} for  a trade-off function $f\in \cT$ \ref{def:fDP}, we replace $f_{\mu_1, \mu_2}$ with its double conjugate $g=f_{\mu_1, \mu_2}^{**}: [0,1] \to [0,1]$, which is the (convex envelope from below) largest lower-semi-continuous convex function $g$ lying below $f_{\mu_1, \mu_2} =\min(f_{\infty}, f_{\infty}^{-1})$ \cite{Rockafellar1970ConvexAnalysis}. However, for a trade-off function $f\in \cT$ we have an explicit expression for the double conjugate $g=\min(f,f^{-1})^{**} \in \cT$, following \cite{dong2022gdp}[Prop E.1]. 

For $f\in \cT$, define the first time the slope (derivative) of $f$ matches that of the function $y(x)= 1-x$ (ideal privacy curve) as $\bar{x}:= \bar{x}_f = \inf\{x \in [0,1] |-1 \in \partial f(x)\}$, where $\partial f(x)$ is the sub-gradient of the trade-off function $f$ at $x$ \cite{Rockafellar1970ConvexAnalysis}. Then in the case $\bar{x} \leq f(\bar{x})$ we have
\begin{equation}
g= \min\{f,f^{-1}\}^{**}(x)
=
\begin{cases}
f(x), & x \in [0,\bar{x}],\\[3pt]
\bar{x} + f(\bar{x}) - x, & x \in [\bar{x}, f(\bar{x})],\\[3pt]
f^{-1}(x), & x \in [f(\bar{x}),1].
\end{cases}
\end{equation}
When $\bar{x} > f(\bar{x})$ below, see  figure \ref{fig:symm-f-example} (see \cite{dong2022gdp}[Figure 6, Appendix] for $\bar{x} \leq f(\bar{x})$ case)
\begin{equation}
g
=
\min\{f,f^{-1}\}^{**}(x)
=
\begin{cases}
f^{-1}(x), & x \in [0,f(\bar{x})],\\[3pt]
\bar{x} + f(\bar{x}) - x, & x \in [f(\bar{x}), \bar{x}],\\[3pt]
f(x), & x \in [\bar{x},1].
\end{cases}
\end{equation}
\end{remark}

\begin{figure}[h]
    \centering
    \includegraphics[width=0.55\linewidth]{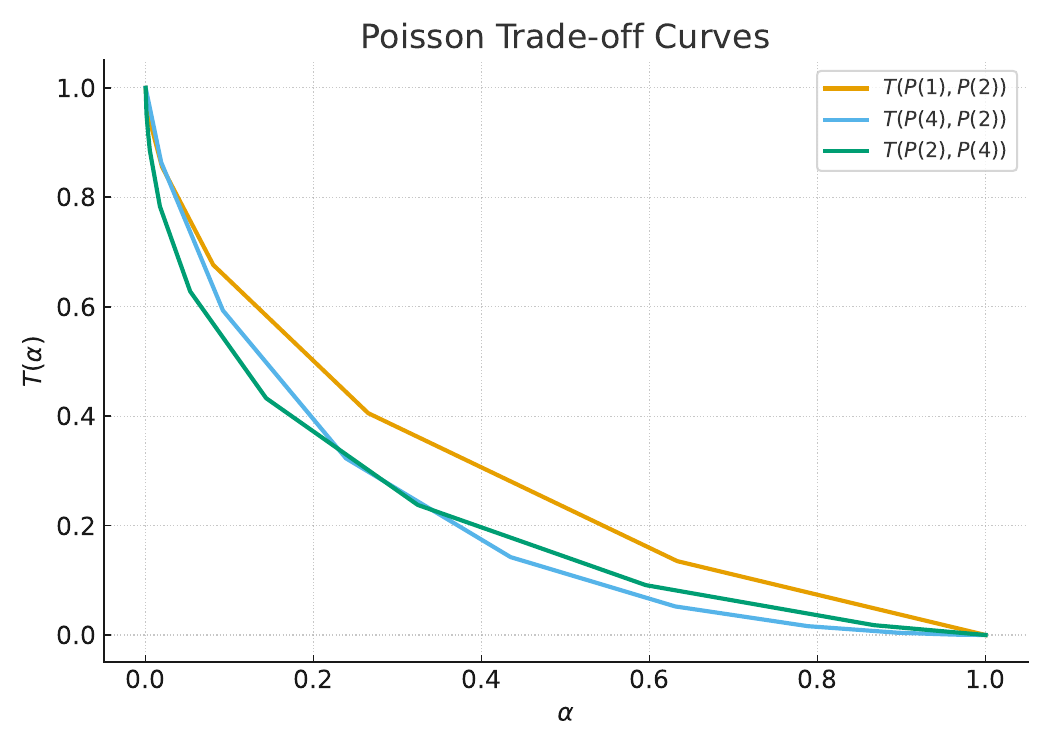}
    \caption{
    Trade-off curves for Poissons:
    $T(P(1),P(2))$, $T(P(4),P(2))$, and $T(P(2),P(4))$.
    }
    \label{fig:poisson-tradeoffs-1-2-4}
\end{figure}

\begin{remark}[Tightness of our mechanism]
    It has been proved in \cite{dong2022gdp}[Theorem 1] (see also Theorem \ref{thm:GDPdong}), that for  a Gaussian baseline trade-off function $f_{\mu} = T(N(0,1),N(\mu,1))$\footnote{This holds more generally for a  shift family $f_{\mu}= T(X, X+\mu)$ with symmetric log concave $X$.},  Gaussian noise mechanism $M(g(x)) =  g(x) + cZ$ with $c^{-1} w_g(1) =\mu$ tightly captures $f_{\mu}$-DP. 
    
    This means that the inequality in Equation \eqref{eq:GDPdong} is an equality and it happens as soon as there exist two datapoints $x_1, x_2$ with $d_{\cX}(x_1, x_2) \leq 1$  achieve the supremum $w_g(1)= $ $|g(x_2)- g(x_1)|$.

    Similarly, if there exist two data points $x_1, x_2$ with $d_{\cX}(x_1, x_2) \leq 1$ achieve the supremum  jointly,
    \begin{equation} \label{eq:jointeq}
    w_g(1)=  |g(x_2)- g(x_1)| \text{ and } w_{h\circ g}(1) =|h(g(x_2))-h(g(x_1))|.
    \end{equation}
Then we have simultaneous equality in Equations \eqref{eq:part1} and \eqref{eq:part2}\footnote{We will have individual equality in Equations  \eqref{eq:part1} and \eqref{eq:part2} according to which equality we require among the supremas  $w_g(1)= $ $|g(x_2)- g(x_1)|$ and $w_{h\circ g}(1) = |h(g(x_2))-h(g(x_1))|$ for some $d_{\cX}(x_1,x_2) \leq 1$. }, since the  pair of constants $(c, \lambda)$ used in the proof becomes $(0, 0)$ causing equality of trade-off functions in Equations \eqref{eq:part1} and \eqref{eq:part2}. 
    
    Therefore,  the mechanism $M(g(x)) \sim P(N_2 e^{N_1 g(x)})$ achieves tightness  \eqref{eq:part1} and \eqref{eq:part2} together (or individually based on what we require on Equation \eqref{eq:jointeq}). However, because of the inherent asymmetry of the Poisson trade-off curves $T(P(\lambda_1), P(\lambda_2))\neq  T (P(\lambda_2), P(\lambda_1))$, and the fact that Blackwell ordering \eqref{eq:Blackwellordering} is a partial order (see figure \ref{fig:poisson-tradeoffs-1-2-4} for a proof of this incomparability)
    \begin{equation}
    T(P(\lambda_1), P(\lambda_2)) \ngeq  T (P(\lambda_2), P(\lambda_1)) \text{ and } 
    T(P(\lambda_2), P(\lambda_1)) \ngeq  T (P(\lambda_1), P(\lambda_2)).
    \end{equation}
We have that the mechanism $M(g(x)) = P(N_2e^{N_1g(x)})$ proposed in Theorem \ref{thm:PDP} overall achieves a $\min(f_{\infty}, f_{\infty}^{-1})^{**}$- differential  privacy \ref{def:fDP} ,  where $f_{\infty}= T(P(\lambda_1), P(\lambda_2)) $.
 \end{remark}

\begin{remark}[Asymmetry in Poisson] Observe that $T(P(\lambda_1), P(\lambda_2)) \neq T(P(\lambda_2), P(\lambda_1))$ (see figure \ref{fig:poisson-tradeoff-g-ginv-f-fss}), as opposed to the symmetry of trade-off curves corresponding to the shift (location) family  $ T(X,X+\mu)= T(X+\mu, X)$ for symmetric-log concave $X$ and $\mu>0$. This asymmetry also happens because of the one-sided nature of the support of $P(\lambda)$ being $\mathbb{Z}_{\geq 0}$. Therefore, we believe it is not natural to consider a symmetric (around zero)  noise-adding mechanism (discrete or continuous) for privacy concerning count statistics $g$ that are inherently non-negative integer valued.
\end{remark}
\begin{remark}[Utility of the framework]
The utility of the $f$-differential privacy framework of \cite{dong2022gdp} lies in the ability to choose the baseline $f$ based on the class of problems at hand. For count statistics, it is natural to have a comparable Poisson trade-off function as the baseline rather than two Gaussians, which are inherently continuous.  In the Poisson case (see Theorem \ref{thm:PDP}), instead of adding Poisson noise, having received the value of the statistic $g(x)$, the mechanism $M$ outputs a sample from the Poisson distribution with parameter $N_2e^{N_1 g(x)}$. 
\end{remark}
\begin{remark}[Almosus sure Implementation in practice]
First, given the privacy level $\mu_2 > \mu_1 >0$, and the statistic $g$ the computation of the constants $N_1, N_2$ (see Equation \eqref{eq:NN}) is immediate, once the modulus of continuity values $w_g(1), w_{h\circ g}(1)$ have been computed\footnote{The values $w_g(1), w_{h\circ g}(1)$ can be unbounded, unless the range of $g$ is bounded, which it will be in many use case scenarios of count statistics, such as degree of a vertex in a sparse network \cite{KarwaSlavkovic2016}.} (see Equation \eqref{eq:h} for the definition of $h$). In Theorem  \ref{thm:PDP}, we have given only the distributional information that for all $x \in \cX$ we have $M(g(x)) \sim P(\lambda(g(x)))$, where $\lambda(n)= N_2e^{N_1n}$, but one can prescribe a sample by sample (almost sure) implementation using the standard conversion between the cumulative distribution function $F$ and the quantile function $F^{-1}$. 
In the Poisson case, more precisely, let  
\begin{equation}
M(n) \sim  P(\lambda(n)), \text{ and }
p_n(k) := P(\lambda(n))(\{k\}) = e^{-\lambda(n)}\frac{\lambda(n)^k}{k!},
\text{ for }k\in\mathbb{Z}_{\ge 0},
\end{equation}
Denote the Cumulative distribution function $F_n: \mathbb{R} \to [0,1]$ for  $P(\lambda(n))$, for $k\in\mathbb{Z}_{\ge 0}$ as\footnote{$F_n$ remains constant between points of jump $F_n(x)= F_n(k)$ for all $k\leq x < k+1$ for all $k \in \mathbb{Z}_{\geq -1}$.}
\begin{equation} \label{CDF}
F_n(k)
:=
P(\lambda(n))(\{0, \cdots, k\})
=
\sum_{j=0}^{k} p_j(n), 
\text{ with } F_n(-1):=0, \text{ and}
\end{equation}
\begin{equation}
0=F_n(-1) < F_n(0) < F_n(1) < \cdots \leq 1, 
\text{ and }\qquad \lim_{k\to\infty} F_n(k)=1.
\end{equation}

Now, almost surely set $M(g(x))=k$ if and only if $F_n(k-1)\leq U < F_n(k).$ for $n=g(x)$, and $U$ is a uniformly random sample from $[0,1]$ (independent of $g(x)$)\footnote{This is the step where we are using external randomness. This is very similar to how one injects independent centered Gaussian noise $G$ (of appropriate variance) into a released statistic $g(x)$ to generate a sample from a Gaussian with mean $g(x)$ to achieve Gaussian differential privacy \ref{thm:GDPdong}.}. Moreover, for each $x \in \cX$ 
\begin{equation}
    \mathbb{P}\bigl(M(g(x))=k\bigr)= U (\{u\in [0,1]: F_n(k-1)\leq u < F_n(k)\})= P(\lambda(n)) (\{k\}) 
    \quad k\in\mathbb{Z}_{\ge 0}.
\end{equation}
So $M(g(x))\sim P(\lambda(g(x)))$ satisfies the distributional requirement\footnote{Observe that, since in practice, one only outputs one noisy output $M(g(x))$ given a query request $g(x)$, and hence there is no requirement on the  joint distributions  of the output $\{M(g(x)): x\in \cX\}$.}. As is immediate, any Markov kernel $M$ with one-dimensional output $M(g(x))$ and a different sampling distribution $F_{g(x)}$ from $P(\lambda(g(x)))$ can also be implemented for a released statistic $g(x)$ in a similar way, almost surely as $F_{g(x)}^{-1}(U)$\footnote{Assume for simplicity, $F_{g(x)}$ has continuous distribution, otherwise, one has to randomize .} where
$F_{g(x)}^{-1}: [0,1] \to \mathbb{R}$ is the quantile function defined as 
\begin{equation}
    F_{g(x)}^{-1} (u) = \inf\{  t \in \mathbb{R}: F_{g(x)} (t)\geq u\}
\end{equation}
\end{remark}

\begin{remark}[Infinitely divisible extension]
    It is an interesting future direction to generalize the result of Theorem \ref{thm:PDP} to describe differentially private mechanisms covering the cases of other baseline trade-off functions $f= T(P, Q)$ with infinitely divisible distributions, such as Geometric, discrete Laplace, and Negative-binomial, more generally for a exponential family of the following form
    \begin{equation}
    f_{\theta_1, \theta_2} =T(P_{\theta_1}, P_{\theta_2}) \text{ where }    \{P_{\theta}: dP_{\theta} (x)= \frac{e^{\theta x} d P(x)}{P(e^{\theta x})}\}_{\theta \in [-1,1]} \text{  with } T(P_0,P_1) \in \cI_{\cT}.
    \end{equation}
 It would also be interesting  to analyse how distribution-specific ($T(P,Q)$) mechanisms \footnote{We believe that except for shift baseline trade-off functions  $f=T(X,X+\mu)$ for some $\mu>0$  mechanisms of the noise adding kind would not be tight in general such as in the Poisson case \ref{thm:PDP}.}   would compare with noise addition mechanisms such as Skellam \cite{AgarwalKairouzLiu2021Skellam},  discrete Gaussian \cite{CanonneKamathSteinke2020DiscreteGaussianDP}\footnote{Although the Gaussian distribution is infinitely divisible, but discrete Gaussian is not \cite{BoseDasGuptaRubin2002}.}, and canonical noise-distribution framework \cite{awan2023canonical} . 
\end{remark}

\begin{remark}[Natural apperance of Poisson trade-off functions] \label{poissonnatural}
We demonstrate how the Poisson differentially private framework with a baseline trade-off function $f_{\infty}= T(P(\lambda_1), P(\lambda_2))$   would be an appropriate benchmark for privacy on a graph-based relational private dataset. 

More precisely, consider a graph-based private dataset $G_n= ([n],E_{n})$, where we would like to  protect the community label vector $\beta \in \{0,1\}^n$ of the vertices $\{1,\cdots, n\}$ upon sequential queries on the labels of the vertices. Since a single coordinate of the entire mechanism $M = (M_{n1}, \cdots, M_{nn})$, where $M_{ni}$ releases a noisy label for vertex $i$, is binary, it is natural to consider that for each $M_{ni}$ (for each $1\leq i \leq n$) satisfies $f_n$-DP \ref{def:fDP} for a baseline trade-off function $f_n= T(\text{Ber}(\lambda_1/n), \text{Ber}(\lambda_2/n))$. Then, as $n \to \infty$ or, equivalently, as the number of vertices within the network grows, upon releasing the (noisy) labels of the vertices, the joint privacy level is given by $f_n^{\otimes n}$ and converges to a Poisson trade-off curve $f_{\infty}= T(P(\lambda_1), P(\lambda_2))$ \ref{cor:Poissonlimit}. 

Now, observe that the natural limits obtained here deviate from Gaussian trade-off curves \cite{dong2022gdp} (see figure \ref{fig:binom-poisson-tradeoff}), and the sequential queries for coordinates of a high-dimensional vector $\beta \in \{0,1\}^n$ can be asked without its connection to a graph-based dataset. For example, the coordinates of $\beta$ might correspond to the characteristics of a single person  on different topics. 

We end by referring to \cite{KarwaSlavkovic2016}, which discusses the $(\varepsilon,\delta=0)$ differentially private framework on an exponential random graph model $G$ and achieves this by adding discrete Laplace noise on the  degree sequence $d=(d_1, \cdots,d_n)$ of the graph $G$. We believe that it would be interesting to work under the Poisson differential privacy framework \ref{thm:PDP} and analyze how the results vary when  we apply our Poisson optimal mechanism $M$ instead.
\end{remark}

Now, we  prove a Poisson Thinning and superposition result required in the proof of Theorem \ref{thm:PDP}. 
\begin{lemma}[Montonicity properties of Poisson trade-off functions under translation and scaling] \label{lem:PTOF}
    \begin{equation} \label{eq:PS1}
        T(P(\lambda_1), P(\lambda_2)) \geq  T(P(c\lambda_1), P(c\lambda_2)) \text{ if } 1 \leq  c <\infty, \lambda_1, \lambda_2 > 0,
    \end{equation}
    \begin{equation} \label{eq:PS2}
        T(P(\lambda_1), P(\lambda_2)) \leq  T(P(c\lambda_1), P(c\lambda_2)) \text{ if } 0 <c \leq 1, \lambda_1, \lambda_2 > 0,
    \end{equation}
    \begin{equation} \label{eq:PT1}
        T(P(\lambda_1), P(\lambda_2)) \leq T(P(\lambda_1 +\lambda), P(\lambda_1 + \lambda)) \text{ if } \lambda, \lambda_1, \lambda_2 > 0,
    \end{equation}
    \begin{equation} \label{eq:PT2}
        T(P(\lambda_1), P(\lambda_2)) 
        \geq T(P(\lambda_1 -\lambda), P(\lambda_1 -\lambda)) \text{ if } \lambda, \lambda_1- \lambda, \lambda_2- \lambda > 0.
    \end{equation}
\end{lemma}
\begin{proof}
    From  Blackwell's theorem \cite{dong2022gdp}[Theorem 2] $T(P,Q) \leq T(KP, KQ)$ for any pairs of probability measures $P,Q$ and Markov kernel $K$. So,  its enough to establish that there exists  Markov kernels $K_c, K_{\lambda} $ defined from  $(\Omega,\mathcal{F})=(\mathbb{Z}_{\geq 0},2^{\mathbb{Z}_{\geq 0}})$ to itself so that  for $0< c <1$ and $\lambda, \lambda_1 , \lambda_2 >0$ we have (by symmetry $(c \to c^{-1}, \lambda \to -\lambda)$, its enough to establish \eqref{eq:PS2} and \eqref{eq:PT1})
     \begin{equation} \label{eq:thinning}
         K_c P(\lambda_1) = P(c\lambda_1) \text{ and }  K_c P(\lambda_2) = P(c\lambda_2) 
     \end{equation}
     \begin{equation} \label{eq:adding}
          K_{\lambda} P(\lambda_1) = P(\lambda_1 +\lambda) \text{ and} 
          K_{\lambda} P(\lambda_2) = P(\lambda_2 + \lambda) 
     \end{equation}

The construction of $K_{\lambda}$ captures the intuition that we are adding the same independent Poisson noise of mean $\lambda$ to both sides. More precisely, define $K_{\lambda}: \mathbb{Z}_{\geq 0} \times 2^{\mathbb{Z}_{\geq 0}} \to [0,1]$ as\footnote{By definition of $K_{\lambda}(n, A)$, we have  $A \to K_{\lambda}(n, A)$ is  a probability measure on $(\mathbb{Z}_{\geq 0},2^{\mathbb{Z}_{\geq 0}})$ for every $n \in \mathbb{Z}_{\geq 0}$, and for every $A \subset \mathbb{Z}_{\geq 0}$, the map   $\mathbb{Z}_{\geq 0} \ni n \to  K_{\lambda}(n, A) \in [0,1]$ is $2^{\mathbb{Z}_{\geq 0}}$ to $\cB_{[0,1]}$ measurable .} 
\begin{equation}
K_{\lambda}(n,A):= \mathbb{P}(n+P\in A), \text{ for } n\in\mathbb{Z}_{\geq 0},\ A\subseteq\mathbb{Z}_{\geq 0}, P \sim P(\lambda).
\end{equation}

Then by the definition of  Markov kernel $K : (\Omega, \cF) \to (\Omega', \cF')$   we have from \cite{Kallenberg2021FMP3}
\begin{equation} \label{eq:MK}
(KP)(A):=\int_{\Omega} K(x,A)dP(dx), 
\text{ for }
A\in\mathcal{F}  
\text{ and a  probability measure } P \text{ on }(\Omega, \cF).
\end{equation} 
So, for all $k\in \mathbb{Z}_{\geq 0}$ with $P(\lambda)(\{k\}) = e^{-\lambda} \frac{\lambda^k}{k!}$ supported on $\mathbb{Z}_{\geq 0}$ we have the following
\begin{equation}
    (K_{\lambda}P(\lambda_1))(\{k\})
  =
  \sum_{0 \leq i \leq k} e^{-\lambda} \frac{\lambda^{k-i}}{(k-i)!} e^{-\lambda_1} 
  \frac{\lambda_1^i}{i!}
  =
  e^{-(\lambda +\lambda_1)} \frac{(\lambda_1 +\lambda)^{k}}{k!}
  =
  P(\lambda + \lambda_1)(\{k\}) 
\end{equation}

Therefore, $K_{\lambda}P(\lambda_1)= P(\lambda+ \lambda_1)$, and  by a similar computation, we have  $K_{\lambda}P(\lambda_2) \equiv P(\lambda+ \lambda_2)$.

The construction of $K_c$ is based on Poisson thinning. For $c\in(0,1]$, define $K_c:\mathbb{Z}_{\geq 0}\times 2^{\mathbb{Z}_{\geq 0}}\to[0,1]$
\[
K_c(n,\{k\})
=
\mathbb{P}\bigl(\mathrm{Bin}(n,c)
:=k\bigr)=\binom{n}{k}c^k(1-c)^{n-k} \mathbf{1}( 0\leq k\leq n)\text{ for }  n\in \mathbb{Z}_{\geq 0}.
\]
Then $K_c$ is a Markov kernel on $(\mathbb{Z}_{\geq 0},2^{\mathbb{Z}_{\geq 0}})$, since $\sum_{k \in \mathbb{Z}_{\geq 0}} K_c(n,\{k\}) =1$ for all $n\in \mathbb{Z}_{\geq 0}$, and for every $A \subset \mathbb{Z}_{\geq 0}$, the map   $\mathbb{Z}_{\geq 0} \ni n \to  K_c(n, A) \in [0,1]$ is $2^{\mathbb{Z}_{\geq 0}}$ to $\cB_{[0,1]}$ measurable.  Now,
\begin{equation}
L= (K_{1}P(\lambda_1))(\{k\})
=\sum_{n=k}^{\infty} 
 \binom{n}{k}c^k(1-c)^{n-k}\,e^{-\lambda_1}\frac{\lambda_1^n}{n!}.
\end{equation}
 Now, applying the change of variable  $n \to k +n$ to the sum above, we have that 
\begin{equation}
L
=
\sum_{n=0}^{\infty} \frac{(n+k)!}{n! k!}\,c^k(1-c)^n\,e^{-\lambda_1}\frac{\lambda_1^{n+k}}{(n+k)!}
=
e^{-\lambda_1}\frac{(c\lambda_1)^k}{k!}\sum_{n=0}^{\infty}\frac{\bigl((1-c)\lambda_1\bigr)^n}{n!}
\end{equation}
which  equals  $e^{-c\lambda_1}\frac{(c\lambda_1)^k}{k!}= P(c\lambda_1)(\{k\})$. Since this is true for any $k\in \mathbb{Z}_{\geq 0}$, we have $K_cP(\lambda_1)= P(c\lambda_1)$. A similar computation shows  $K_cP(\lambda_2)= P(c\lambda_2)$, thereby establishing \eqref{eq:thinning}.
\end{proof}

\section{Neyman-Pearson lemma under a coarser $\sigma$-algebra}

The Neyman-Pearson optimizer  $\varphi^*$ of the trade-off function $T(P,Q)$ for probability measures $(P,Q)$ on a measurable space $(\Omega, \cF)$ involves the likelihood ratio (Radon-Nikodym derivative) random variable $F= \frac{dQ}{dP}: (\Omega, \cF) \to [0,\infty]$ \eqref{eq:LRT}, but we pose a natural question: what happens when we shrink the original $\sigma$-algebra $\cF$ to a smaller one $\cG$\footnote{Shrinking the $\sigma$-algebra to a smaller $\sigma$ algebra captures the intuition of blurring the information in subsets of the sample space $\Omega$ by reducing the number of allowed subsets from $\cF$ to $\cG$. A complementary way would be to say that coarsening information is the same as testing with fewer allowed observables $\varphi$, which are available at our disposal. See \cite{RevuzYor1999Continuous} for more details on the intuition behind \textit{filtration} of $\sigma$-algebras.}, so that the Radon-Nikodym derivative $H= \frac{dQ}{dP}$ is not a $\cG$ measurable random variable; therefore, $\varphi^*$ cannot be the Neyman-Pearson optimizer for $(P,Q)$ on $(\Omega, \cG)$. In this section, we answer this question in full generality. Our motivation comes from a Bayesian perspective of differential privacy as proposed in \cite{StrackYang2024PrivacyPreserving}, whose connections with the $f$-DP framework will also be an interesting future direction.

\begin{theorem} \label{thm:NP-coarsen}
Consider a measurable space $(\Omega,\mathcal F)$, and let $P,Q$ be probability measures on $(\Omega,\mathcal F)$\footnote{We do not need this assumption of absolute continuity of $Q\ll P$, and the result holds more generally when we extend our definition of the Radon-Nikodym derivative $\frac{dQ}{dP}:= \frac{dQ}{d\mu} / \frac{dP}{d\mu}$ for a dominating  $2\mu= P+Q$.}. Denote the likelihood ratio
$F:=\frac{dQ}{dP}: (\Omega, \cF) \to [0,\infty]$ (defined only $P$-almost surely). Let $\mathcal G\subseteq \mathcal F$ be a sub-$\sigma$-algebra (a coarsening), and $P^{\mathcal G}$ and $Q^{\mathcal G}$ are the restrictions of $P,Q$ to $\mathcal G$.  Define the \emph{conditional likelihood ratio} $G\;:=\; \mathbb{E}_P\!\left[F \,\middle|\, \mathcal G\right]: (\Omega, \cG) \to [0,\infty]$.  Then,
\begin{itemize}
    \item $G$ is (a $P^{\mathcal G}$-almost sure version of) the Radon–Nikodym derivative $\frac{dQ^{\mathcal G}}{dP^{\mathcal G}}$ on $(\Omega, \cG)$.
\end{itemize}

Now, consider the coarsened trade-off function $T_{\cG}(P,Q)\equiv T(P^{\mathcal G},Q^{\mathcal G}): [0,1] \to [0,1]$ defined as
\begin{equation}
T(P^{\mathcal G},Q^{\mathcal G})(\alpha)
=
\inf_{\varphi}
\Bigl\{\e_{Q}[1-\vp] 
        \,\Bigm\vert\,
       \e_{P}[\vp]\le \alpha,\;
        \vp:(\Omega, \cG)\!\to[0,1]\text{ measurable}
      \Bigr\}.
\end{equation}
 \begin{itemize}
     \item The Neyman-Pearson optimizer $\varphi_{\cG}^*$ achieving $T_{\cG}(P,Q)$ at level $\alpha \in [0,1]$ is given by \eqref{eq:LRT}
     \end{itemize}
\begin{equation}\label{eq:LRTC}
            \varphi^*(w) :=\mathbf{1}\left( \frac{dQ^{\mathcal{G}}}{dP^{\mathcal{G}}} >\tau^* \right) + \lambda^* \mathbf{1}\left(\frac{dQ^{\mathcal{G}}}{dP^{\mathcal{G}}} =\tau^* \right), \quad \frac{dQ^{\mathcal{G}}}{dP^{\mathcal{G}}}: = \frac{dQ^{\mathcal{G}}}{d\mu^{G}}/\frac{dP^{\mathcal{G}}}{d\mu^{\cG}}:(\Omega, \cG) \to [0, \infty]
        \end{equation}
  for $2\mu^{\mathcal{G}}=P^{\mathcal{G}}+Q^{\cG}$, the ratios of Radon–Nikodym derivatives, and $\tau^*, \lambda^*$ are determined  by choosing $\tau^* \in [0,\infty]$ as the unique number 

  such that $P^{\cG}\left(\frac{dQ^{\cG}}{dP^{\cG}} \geq \tau^{*}\right) \geq \alpha \geq  P^{\cG}\left(\frac{dQ^{\cG}}{dP^{\cG}} > \tau^{*}\right)$, 
  \[\lambda^*:= \frac{\alpha - P^{\cG}\left(\frac{dQ^{\cG}}{dP^{\cG}} > \tau^{*}\right) }{P^{\cG}\left(\frac{dQ^{\cG}}{dP^{\cG}} \geq \tau^{*}\right)  - P^{\cG}\left(\frac{dQ^{\cG}}{dP^{\cG}} > \tau^{*}\right) } \mathbf{1}\left(\alpha - P^{\cG}\left(\frac{dQ^{\cG}}{dP^{\cG}} > \tau^{*}\right) >0\right) \text{ so } \mathbb{E}_{P}[\varphi_{\cG}^*]=\alpha.\] 

\begin{itemize}
    \item Moreover, the trade-off function satisfies  data-processing inequality (Blackwell ordering)
\end{itemize}
\begin{equation}
T(P^{\mathcal G},Q^{\mathcal G})(\alpha)\;\ge\;T(P,Q)(\alpha) \text{ for all }\alpha\in[0,1], \text{ with equality for all } \alpha \text{ if and only if }
\end{equation}

 $F$ itself is $\mathcal G$-measurable $P$ almost surely or equivalently, $\mathbb{E}_{P}[F|\cG]=F$, $P$ almost surely\footnote{A complementary way is to say that the information $\cG$ is sufficient for the binary experiment $(\Omega, \cF, P,Q)$.}.
 
\end{theorem}

\begin{proof}
First, we show that $G$ is the Radon-Nikodym derivative $\frac{dQ^{\mathcal G}}{dP^{\mathcal G}}$. Now, observe that for any measurable function $(\Omega,\mathcal G) \overset{\varphi}{\to} [0,1]$ , we have  $\mathbb{E}_Q[\varphi]
= \mathbb{E}_P[F\,\varphi]
= \mathbb{E}_P\!\big[\,\mathbb{E}_P(F\mid \mathcal G)\,\varphi\,\big]
= \mathbb{E}_P[G\,\varphi].$

Second, the constrained minimization of $\mathbb{E}_Q[1-\varphi]$ over
$\mathcal G$-measurable $\varphi$ with $\mathbb{E}_P[\varphi]\le\alpha$ reduces to the standard Neyman-Pearson problem for the reduced experiment
$(\Omega,\mathcal G,P^{\mathcal G},Q^{\mathcal G})$ with likelihood ratio
$\frac{dQ^{\mathcal G}}{dP^{\mathcal G}}$. The Neyman-Pearson lemma \eqref{eq:LRTC} gives the threshold form of $\varphi_{\cG}^*$.

Finally, the inequality
$T(P^{\mathcal G},Q^{\mathcal G})\ge T(P,Q)$ follows immediately because  $\varphi: (\Omega, \cG) \to [0,1]$ is measurable implies $\varphi: (\Omega, \cF) \to [0,1] $ is measurable. Moreover, if $\mathbb{E}_{P}[F|\cG]=F$, then $\varphi_{\cG}^{*}=\varphi^*$ in \eqref{eq:LRTC} and hence equality holds for all $\alpha$ holds. Conversely, equality holds at all $\alpha \in [0,1]$, by the uniqueness of Neyman-Pearson optimizer, $\varphi_{\cG}^{*}= \varphi^*$ at all $\alpha$, and therefore $F=\frac{dQ}{dP}$ itself is $\cG$ measurable\footnote{By the restriction of $\varphi^*_{\cG}= \varphi^*$ at all $\alpha \in [0,1]$, we need $\{F >\tau\}_{\tau \geq 0}$ and $\{F=\tau\}_{\tau \geq 0}$ to be $\cG$ measurable, and therefore $F$ itself is $G$ measurable, since it is by definition non-negative.}. So, $G=\mathbb{E}_{P}[F|\cG]= F, P$ almost surely.
\end{proof}
\subsection{A Gaussian shift example on a sub-$\sigma$ algebra}
\textbf{Some consequences in practice.} Now, we describe the theory developed above with a Gaussian shift example of how privacy (trade-off function) degrades upon coarsening the $\sigma$ algebra\footnote{The results extend verbatim with appropriate substitution for any shift family $\{P_{\theta}(\cdot)= P(\cdot - \theta): \theta \in \}$  where $P_{0}$ is a symmetric-log concave probability measure on $(\mathbb{R}, \cB_{\mathbb{R}})$ \cite{dong2022gdp}.}.

Consider a Gaussian shift family $\{P_{\theta}(\cdot)= P(\cdot - \theta): \theta \in \mathbb{R}\}$ with $P=P_0= N(0,1)$ and for $\mu>0$ let $Q=P_{\mu}= N(\mu,1)$ on $(\Omega, \cF)= (\mathbb{R}, \cB_{\mathbb{R}})$. Now, coarse to the $\sigma$-algebra to define $\cG$ as\footnote{$\sigma$-algebra generated by a countable partition $\{E_n\}_{n\in \mathbb{Z}}$ of  $\Omega$ can be  written as $\{\{\bigcup_{A \subset \mathbb{Z}} E_n: A\subset \mathbb{Z}\}\}.$}
\begin{equation}
\mathcal G
=
\sigma\big([n,n{+}1): n\in\mathbb Z\big) 
=
\{\bigcup_{n \in A} [n,n+1): A\subset \mathbb{Z}\},
\end{equation}
So this $\sigma$-algebra is naturally isomorphic to the $\sigma$-algebra $2^{\mathbb{Z}}$ (power set) under the isomorphism $A \to \bigcup_{n \in A} [n, n+1)$, and the  test functions $\varphi$ for discriminating between $P$ and  $Q$ on $(\Omega, \cG)$ only allowed to depend on the bin index \(N=\lfloor X\rfloor\in\mathbb Z\) of the observed value\footnote{The intuition is that the practitioner does not observe the actual output $X$, but only see the processed output $\lfloor X\rfloor$ taking integer values, and is supposed to base his discriminating procedure based on this information.}.

\begin{proposition} \label{prop:GTOFC}
Let $P_0= N(0,1), P_{\mu}= N(\mu,1)$ for $\mu>0$ on $(\Omega, \cF) \equiv (\mathbb{R}, \cB_{\mathbb{R}})$ with a sub-$\sigma$ algebra $\cG:= \sigma([n, n+1): n \in \mathbb{Z})$. Then, $
        T(P_0, P_{\mu})(\alpha) = \Phi(\Phi^{-1}(1-\alpha)- \mu) \forall  \alpha \in [0,1],$ \text{ and}
    \begin{equation}
        T_{\cG}(P_0, P_{\mu})(\alpha_k) =  \Phi(\Phi^{-1}(1-\alpha_k)- \mu) \text{ for } \alpha_k = 1- \Phi(k) \text{ for all } k\in \mathbb{Z}, \text{ and}
    \end{equation}
    \begin{equation}
        T_{\cG}(P_0, P_{\mu})(\alpha)=  \; T_{\cG}(P_0,P_{\mu})(\alpha_{k+1}) - \lambda^*\,q_k, \text{ for } \alpha \in [\alpha_{k+1}, \alpha_k] \text{ for all } k \in \mathbb{Z}
        \end{equation}
        \begin{equation}
        \text{and }\lambda^* = \frac{\alpha-\alpha_{k+1}}{\alpha_k - \alpha_{k+1}}, \text{ with }
        q_k= \Phi(k+1-\mu)- \Phi(k-\mu)
    \end{equation}
\end{proposition}
\begin{proof}
The proof follows from \ref{thm:NP-coarsen}. First, observe that the  laws $\mathbb{P}^{\cG}$  and $\mathbb{Q}^{\cG}$ on $(\Omega, \cG) \equiv (\mathbb Z, 2^{\mathbb{Z}})$ are
\begin{equation}
p_n:=P(N=n)=\Phi(n{+}1)-\Phi(n),\qquad
q_n:=Q(N=n)=\Phi(n{+}1-\mu)-\Phi(n-\mu),
\end{equation}
where \(\Phi\) is the standard normal CDF. For \(\mu>0\) the likelihood ratios $F_n$  are increasing in $n \in \mathbb{Z}$\footnote{A similar computation would show that $F_n$ would be decreasing in $n$  when the shift  \(\mu<0\)).}.
\begin{equation}
F_n:=\frac{q_n}{p_n}
=
\frac{\int_{n}^{n+1} \exp( -\frac{(x-\mu)^2}{2}) dx}{\int_{n}^{n+1} \exp(-\frac{(x)^2}{2}) dx} 
=
\frac{\int_{n}^{n+1} F_{\mu}(x) dP_0(x)) dx}{\int_{n}^{n+1} dP_0(x)}
=
\mathbb{E}_P[F_{\mu}|\cG](x) \text{ for } x\in [n, n+1),
\end{equation}
where we denote $F_{\mu}(x)= \frac{dP_{\mu}}{dP_0}(x)= \frac{\exp( -\frac{(x-\mu)^2}{2})}{\exp(-\frac{x^2}{2})}$ as the Radon-Nikodym derivative, and a $\cG$ measurability requirement of $G=\mathbb{E}_P[F_{\mu}|\cG]$ reveals that $\mathbb{E}_P[F_{\mu}|\cG]$ is constant on each of those partitions $E_n= [n,n+1)$ generating the $\sigma$-algebra $\cG$ \cite{Durrett2019}. Moreover, from the  integrability requirement $\mathbb{E}_P[F_{\mu} \varphi]= \mathbb{E}_P[G\varphi]$ for bounded $\cG$ measurable functions $\varphi$ that constant value is the  weighted average as computed above (take $\varphi \equiv 1$). Now, to show that $F_n$ is increasing in $n$, observe that for any $\mu> 0$, $\frac{dP_{\mu}}{dP_0}(x) = \exp\left(\mu x- \frac{\mu^2}{2}\right)$ is monotone in $x$ whenever $\mu>0$\footnote{As shown in \cite{dong2022gdp}, this condition of monotonicity (increasing) of $\frac{dP_{\mu}}{dP_0}(x)$ in $x$ for all (but fixed) $\mu>0$ is equivalent to the condition that $P_0$ is log concave. This is the reason this Gaussian  result extends verbatim to symmetric-log concave distribution $P_0$ with the shift family  $\{P_{\theta}(\cdot)= P(\cdot - \theta): \theta \in \mathbb{R} \}$.}. As a consequence, we have $F_{\mu}(n) \leq F_{n} \leq F_{\mu}(n+1) \leq F_{n+1} \leq  F_{\mu}(n+2)$ for every $n \in \mathbb{Z}$.

So the Neyman-Pearson optimizer on \((\Omega,\mathcal G) \equiv (\mathbb{Z}, 2^{\mathbb{Z}})\) is a threshold in \(N\)\footnote{To follow the theorem \ref{thm:NP-coarsen}, one could equivalently threshold at the value based on $F_n$ as well.} (along with some randomization as required). For the non-randomized part, consider an  threshold \(k \in \mathbb{Z}\),
\begin{equation}
\text{type I error}= \alpha_k := P(N\ge k)=\sum_{n\ge k} p_n = 1-\Phi(k),
\text{ and }
\end{equation}
\begin{equation}
\text{type II error}=T_{\cG}(P,Q)(\alpha_k) := Q(N< k)= \sum_{n< k} q_n= \Phi(k-\mu).
\end{equation}
Thus, the trade-off function $T_{\cG}(P^{\cG},Q^{\cG})$ at the points with type I errors $\{\alpha_k\}_{k\in \mathbb{Z}} \subset [0,1]$\footnote{These are the values of type I errors which do not require randomization in Neyman-Pearson lemma.},   have the respective minimum value of type II errors as $\{\Phi(k-\mu)\big)\}_{k\in \mathbb{Z}}$  or equivalently as the subset of points $\{\alpha_k, \Phi(\Phi^{-1}(1-\alpha_k)- \mu)\}_{k\in \mathbb{Z}} \subset \{\alpha, T_{\cG}(P,Q)(\alpha)\}_{\alpha \in [0,1]}\subset [0,1]^2$. Now, to have the entire  $T_{\cG}(P,Q)$, one linearly interpolates the points above on the curve. 
More precisely, for any \(\alpha\in[\alpha_{k+1},\,\alpha_k]\), define $\varphi_{\cG}^*(N) = \mathbf{1}(N >k) + \lambda^* \mathbf{1}(N =k),$ where 
\(\lambda^*\in[0,1]\) is given by
\begin{equation}
\lambda^* = \frac{\alpha-\alpha_{k+1}}{\alpha_k - \alpha_{k+1}}, \text{ and } 
T_{\cG}(P,Q)(\alpha)
\;=\; T_{\cG}(P,Q)(\alpha_{k+1}) - \lambda^*\,q_k.
\end{equation}
\end{proof}
This finishes the description of the computation of the  trade-off function $T_{\cG}(P,Q)$ for a Gaussian shift family, which is clearly equal to $T(P,Q)$ at the points of type $I$ errors $\{\alpha_k\}_{k\in \mathbb{Z}}$ without requiring randomization, and linear interpolated otherwise. Therefore, it always lies above the original Gaussian trade-off curve $T(P,Q)$. Moreover, the \(\mu<0\), case be dealt in exactly similar way.

\begin{figure}[h]
  \centering
  \includegraphics[width=0.7\linewidth]{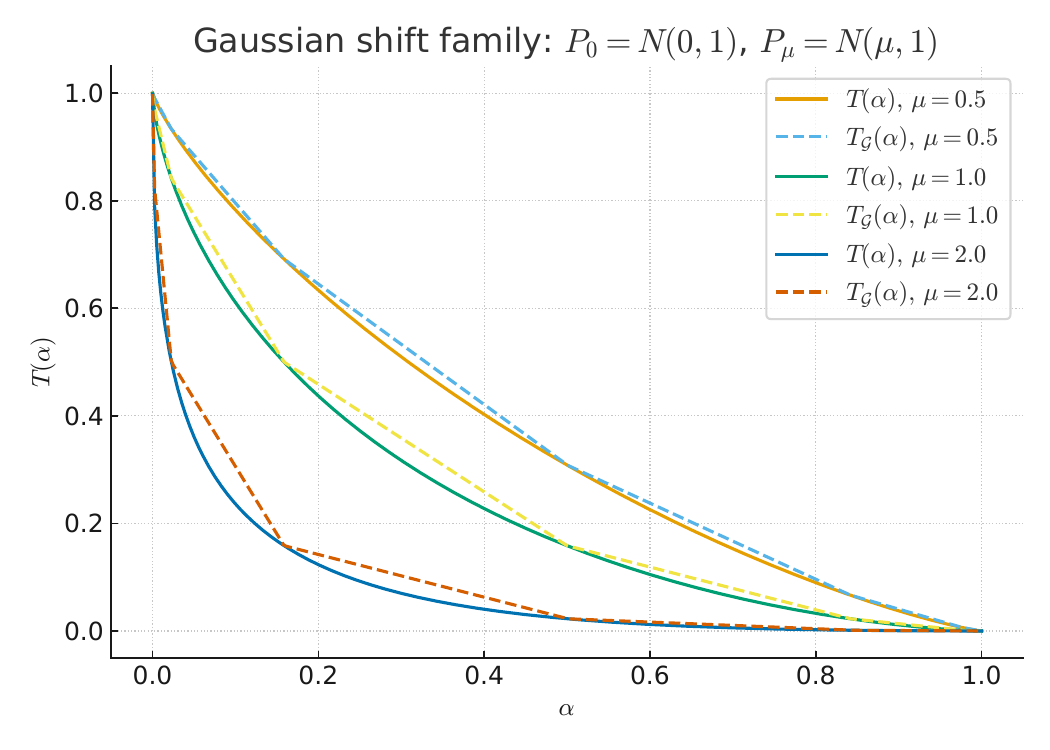}
  \caption{Trade-off function $T(\alpha)$ and its coarsened version
    $T_{\mathcal G}(\alpha)$ for the Gaussian shift family
    $P_0 = N(0,1)$, $P_\mu = N(\mu,1)$.}
  \label{fig:gaussian-tradeoff}
\end{figure}
\begin{figure}[h]
  \centering
  \includegraphics[width=0.7\linewidth]{fig_laplace_tradeoff.pdf}
  \caption{Trade-off function $T(\alpha)$ and its coarsened version
    $T_{\mathcal G}(\alpha)$ for the Laplace shift family
    $P_0 = \mathrm{Lap}(0,1)$, $P_\mu = \mathrm{Lap}(\mu,1)$.}
  \label{fig:laplace-tradeoff}
\end{figure}

\section{Conclusion and Future Work}
\label{sec:conclusion}

In this paper, we introduce the hypothesis testing-based framework of infinitely divisible privacy \ref{thm:IDP}, encompassing both the Gaussian differential privacy of \citet{dong2022gdp} and the newly proposed Poisson differential privacy \ref{cor:Poissonlimit}, which are not just natural baselines for count statistics but also appear under repeated compositions of (nearly perfect) differentially private outputs. Along the way, we resolve the $s^2=2k$ conjecture of \cite{dong2022gdp} as a consequence of contiguity  \ref{thm:GDP}. 

In practice, the number of private operations made on a dataset is uncertain. We show that under randomization of the number of operations, the limiting trade-off can even \emph{escape} the universal class of infinitely divisible trade-off functions \ref{eq:Levylimit}  We illustrate this with a locally asymptotically mixed normal (LAMN) example \eqref{eq:GMM}, following the spirit of \citet{LeCam1986}.

Then, we construct a mechanism in Theorem~\ref{thm:PDP} that \emph{optimally} attains (asymmetric) Poisson differential privacy for real-valued statistics, including count statistics, thus going beyond standard additive-noise mechanisms. Afterwards, we demonstrate that Poisson baselines naturally appear for graph-based datasets, where released statistics are inherently discretized \ref{poissonnatural}.

Finally, inspired by a Bayesian perspective on DP \citep{StrackYang2024PrivacyPreserving}, we quantify exactly how trade-off functions degrade under \textit{information coarsening}   and derive a coarsened Neyman-Pearson lemma \ref{thm:NP-coarsen} with equality characterizations via conditional likelihood ratios ~\ref{thm:NP-coarsen}. We demonstrate this by explicitly working out the Gaussian case ~\ref{prop:GTOFC} that extends to symmetric-log concave family.

Several interesting directions follow, and many have been detailed in the remarks within the text. First, our algorithmic result suggests the construction of private mechanisms that are optimal with respect to a given natural, infinitely divisible trade-off curve. Second, it would be extremely interesting to illustrate the tightness of the proposed Poisson differential privacy framework on graph-based problems of modern interest \cite{KarwaSlavkovic2016}.   Third, it would be interesting to analyze how privacy degrades when one is allowed to stop adaptively. Finally, it would be fascinating to find an equivalence between the $f$-differential privacy framework of \cite{dong2022gdp}  based on comparison of trade-off functions and the  (Bayesian)  framework of privacy-preserving signals \cite{StrackYang2024PrivacyPreserving}, since coarsening the underlying $\sigma$ algebra clarifies how information restrictions lead to an increase in indistinguishability.

\newpage 

\bibliographystyle{iclr2026_conference}
\bibliography{references}

\appendix
\section{Supplementary material }
\label{app:proofs}
We first define  a Markov kernel (randomized mechanism) following \cite{Kallenberg2021FMP3}.
\begin{definition}[Markov kernel] \label{def:MK}
Let $(\cX, \cF_{\cX})$ and $(\cY, \cF_{\cY})$ be measurable spaces.
A \emph{Markov kernel} (also called a \emph{stochastic kernel} or \emph{transition kernel})
from $(\cX, \cF_{\cX})$ to $(\cY, \cF_{\cY})$ is a map
\begin{equation}
K : \cX \times \cF_{\cY} \to [0,1], \text{ usually written } K(x, B) \text{ for } x \in \cX, B \in \cF_{\cY},  
\end{equation}
  such that for every fixed $x \in X$, the set function $B \mapsto K(x, B)$ is a probability measure on $(\cY, \cF_{\cY})$. Moreover, for every fixed $B \in \cF_{\cY}$, the function $x \mapsto K(x, B)$ is $\cF_{\cX}$-measurable.
\end{definition}

We state the following Berry Esseen CLT, where $\mathbf{kl}$ denotes the vector
$(\operatorname{kl}(f_1), \ldots, \operatorname{kl}(f_n))$ and
$\boldsymbol{\kappa}_2$, $\boldsymbol{\kappa}_3$, $\overline{\boldsymbol{\kappa}}_3$
are defined similarly; in addition, $\|\cdot\|_1$ and $\|\cdot\|_2$ are the $\ell_1$ and $\ell_2$ norms, respectively.

Our \textbf{first} lemma summarizes some of the basic but extremely fruitful properties of the Gaussian trade-off function $f_{G, \ve}(x)=T(N(0,1), N(\ve,1))(x)= \Phi(\Phi^{-1}(1-x) - \ve)$ for all $x \in [0,1]$, and any $\ve \geq 0$. Among others, it proves that $T\left(N\left(-\frac{\mu^2}{2}, \mu^2\right), N\left(\frac{\mu^2}{2}, \mu^2\right)\right)= f_{G,|\mu|}$ for any $\mu \in \mathbb{R}$.

\begin{lemma} (\cite{dong2022gdp})\label{lem: GTOF} The Gaussian trade off functions satisfy the following properties:
\begin{itemize}
\item Monotonicity:  For any pair  $\ve_1, \ve_2  \geq0$, $\ve_1 \leq \ve_2$ if and only if  
\begin{equation} \label{eq:GMontone}
   f_{ G,\ve_1}(x)  \geq f_{G,\ve_2}(x) \text{ for all } x \in [0,1]. 
\end{equation}
\item Closure under suprema:  For any collection of $(\ve_i)_{i\in I} \subset 
\R$ with index set $I$
\begin{equation} \label{eq:Gsuprema}
   \inf_{i\in I} f_{ G,\ve_i}(x) = f_{G,\sup_{i\in I}\ve_i}(x) \text{ for all } x \in [0,1]
\end{equation}
\item Symmetry: For any  $\mu_1, \nu_2 \in \R^p$ and $\sigma >0$ 
\begin{equation} \label{eq:symmetry}
   T(\nu_2 + \sigma N(0, \II_p), \mu_1+ \sigma N(0, \II_p)) =   T(\mu_1+ \sigma N(0, I_p), \nu_2 + \sigma N(0, I_p))
\end{equation}
\item Dimension freeness: For any  $\mu_1, \nu_2 \in \R^p$ and $\Sigma  \succ 0$ let  $\ve^2 := \langle (\mu_1- \nu_2), \Sigma^{-1} (\mu_1- \nu_2)\rangle $. Then
\begin{equation} \label{eq:dto1a}
    T(\mu_1+ \sqrt{\Sigma} N(0, \II_p), \nu_2 + \sqrt{\Sigma}N(0, \II_p)) \equiv T(N(0,1), N(\ve, 1)), 
\end{equation}
\end{itemize}
 \end{lemma}
\textbf{Intuition and importance:} The proof of this lemma is immediate from the explicit description of $f_{G, \ve}(x)=\Phi(\Phi^{-1}(1-x) - \ve)$ where $\Phi(\cdot)$ is the CDF of a standard Gaussian variable $\Phi(t)= \mathbb{P}[N(0,1)\leq t]$ for all $t\in \mathbb{R}$. It also requires applying the Neyman-Pearson lemma or the likelihood ratio test \cite{polyanskiy_wu_2025}. But, the conclusions that they imply are extremely powerful. 
\par
\textbf{Monotonicity:} The \textit{monotonicity} condition \eqref{eq:GMontone} reduces an apriori difficult functional comparison between two functions $f$ and $g$ at uncountably many points to a comparison of just one parameter $\varepsilon$.
\par
\textbf{Closure under suprema:} The closure under suprema property \eqref{eq:Gsuprema} essentially says two very important things. First, it makes it easy to identify what the suprema of an apriori arbitrary collection of functions $\{f_i\}_i$ is (in fact explicitly). Second, the limiting object is a function of the same kind: it is again a Gaussian trade-off function with a different choice of parameter $\ve$.
\par
\textbf{Symmetry:} The symmetry property  \eqref{eq:symmetry}  (requires Neyman-Pearson lemma) is very interesting because the definition of a trade-off function $T(P,Q)$ ~\ref{def:TOF} is asymmetric in general, between its first and second arguments.  But, for a pair of shifted isotropic Gaussians, they match because of the spherical symmetry (orthogonal invariance) of the standard Gaussian density in any dimension.
\par
\textbf{Dimension freeness:} The dimension freeness property \eqref{eq:dto1a}  (requires Neyman-Pearson lemma)  makes the case for Gaussian certifiability in high dimensions stronger than any other notion of certifiability. This is because, almost all the results of classical statistics that are true in low dimensions, fail to hold in high dimensions, because many of the quantities involved in controlling the errors are dimension dependent and blows up when $p \uparrow \infty$. It is often the case that finding a dimension-free quantity or even an inequality that \textit{tensorizes} \footnote{We are referring dimension free Poincare and Logarithmic Sobolev inequalities of high-dimensional statistics that are extremely important in obtaining concentration bounds in high dimensions \cite{vanhandelAPC550}.} help us resolve high dimensional issues.

Our \textbf{second} lemma says that the some aspects of the Gaussian trade-off functions generalize. 
\begin{lemma}\cite{dong2022gdp}[Lemma A.2, Proposition A.3] \label{lem:SLCTOF}
Consider a \textit{symmetric} random variable $X\overset{d}{=} -X$ with CDF $F$ having a log-concave Lebesgue density. Consider its  baseline trade-off function $f_{X,\ve}(x)=T(X,\ve+X)(x)$ for all $x \in [0,1]$. 
\begin{enumerate}
    \item \textit{Symmetry:} For any $t_1,t_2 \in \mathbb{R}$ and $\delta = |t_2- t_1|$ the following holds  for all $x \in [0,1]$.
    \begin{equation}
T(t_1+X, t_2+X)(x) = F\left(F^{-1}(1-x)-\delta \right)=T(t_2+X,t_1+X)(x). 
\end{equation}
\begin{equation}
\text{So, }T(t_1+X, t_2+X)(x)= T(X,\delta+X)(x)=T(\delta+X,X)(x) 
\end{equation}
\item \textit{Monotonicity:}
For $\ve_1 ,\ve_2 \geq 0$. Then $\ve_1 \leq \ve_2$ if and only if 
\begin{equation}
    f_{X,\ve_1}(x) \geq f_{X,\ve_2}(x) \text{ for all } x \in [0,1].
\end{equation}
\end{enumerate}
\end{lemma}
\textbf{Intuition and importance:}  The proof of this lemma is along the same lines of the proof given in \cite{dong2022gdp}[Lemma A.2, Proposition A.3]. Using log-concavity of the density of $X$ and symmetry $X \overset{d}{=}-X$, there is an explicit description of $f_{X, \ve}(x)=F(F^{-1}(1-x) - \ve)$ where $F(\cdot)$ is the CDF of the random variable $X$, $F(t)= \mathbb{P}[X\leq t]$ for all $t\in \mathbb{R}$, and $\ve \geq 0$.

\end{document}